%% file: gen2.tex
\documentstyle[txmac,a4,
amssymb,%
case,%
twoside,%
nocaphead,%
epsf,%
varthm,%
rotate,%
myrot%
,multicol%
,mypic,times,mathptm]{article}

\advance\oddsidemargin by -1.9cm
\advance\evensidemargin by -1.9cm
\advance\textwidth by 3.8cm

\def\mynewtheo#1#2{%
\newtheorem{@#1}{#2}[section]%
\newenvironment{#1}{\begin{@#1}\rm}{\end{@#1}}}

\mynewtheo{lemma}{Lemma}
\mynewtheo{exer}{Exercise}
\mynewtheo{theo}{Theorem}
\mynewtheo{rem}{Remark}
\mynewtheo{defi}{Definition}
\mynewtheo{conj}{Conjecture}
\mynewtheo{corr}{Corollary}
\mynewtheo{prop}{Proposition}
\mynewtheo{question}{Question}
\mynewtheo{exam}{Example}

\newenvironment{theorem}{\begin{theo}}{\end{theo}}
\newenvironment{conjecture}{\begin{conj}}{\end{conj}}

\newenvironment{eqn}{\begin{equation}}{\end{equation}}

\begin{document}

\makeatletter

\input{myeqn.tex}

\newcommand{\mybin}[2]{\text{$\Bigl(\begin{array}{@{}c@{}}#1\\#2%
\end{array}\Bigr)$}}
\newcommand{\mybinn}[2]{\text{$\biggl(\begin{array}{@{}c@{}}%
#1\\#2\end{array}\biggr)$}}

\def\overtwo#1{\mbox{\small$\mybin{#1}{2}$}}
\newcommand{\mybr}[2]{\text{$\Bigl\lfloor\mbox{%
\small$\displaystyle\frac{#1}{#2}$}\Bigr\rfloor$}}
\def\mybrtwo#1{\mbox{\mybr{#1}{2}}}

\def\myfrac#1#2{\raisebox{0.2em}{\small$#1$}\!/\!\raisebox{-0.2em}{\small$#2$}}
\def\ffrac#1#2{\mbox{\small$\ds\frac{#1}{#2}$}}

\def\myeqnlabel{\bgroup\@ifnextchar[{\@maketheeq}{\immediate
\stepcounter{equation}\@myeqnlabel}}

\def\@maketheeq[#1]{\def\theequation{#1}\@myeqnlabel}

\def\@myeqnlabel#1{%
{\edef\@currentlabel{\theequation}
\label{#1}\enspace\eqref{#1}}\egroup}

{\catcode`\^^M=\active%
\gdef\@entry{\catcode`\^^M=\active\let^^M\@kentry\frenchspacing%
\catcode`\ 11\relax%
\begin{center}\begin{eqnarray*}}%
{\catcode`\ =11\relax%
\gdef\@@kentry#1 #2:{\@vobeyspaces\@@@kentry{#1}{#2}:}
}%
\gdef\@@@kentry#1#2:#3^^M{
\hbox to 0.5cm{$#1_{#2}$\hss} & & {\verb|#3|}\\[-2pt]^^M}%
}
\newenvironment{kentry}{\@entry}{\end{eqnarray*}\end{center}}
\def\@kentry{\@ifnextchar\end{}{\@@kentry}}


\def\rato#1{\hbox to #1{\rightarrowfill}}
\def\hookrato#1{\hbox to #1{$\lhook\joinrel$\rightarrowfill}}

\def\arrowname#1{{\enspace
\setbox7=\hbox{F}\setbox6=\hbox{%
\setbox0=\hbox{\footnotesize $#1$}\setbox1=\hbox{$\to$}%
\dimen@\wd0\advance\dimen@ by 0.66\wd1\relax
$\stackrel{\rato{\dimen@}}{\copy0}$}%
\ifdim\ht6>\ht7\dimen@\ht7\advance\dimen@ by -\ht6\else
\dimen@\z@\fi\raise\dimen@\box6\enspace}}
 
\def\leftrightharpoonfill{$ \mathord \leftharpoondown\mkern -17mu 
\m@th \smash -\mkern -7mu\cleaders \hbox {$\mkern -2mu\smash -\mkern
-2mu$}\hfill \mkern -7mu\mathord \rightharpoonup$}
\def\lrhto#1{\hbox to #1{\leftrightharpoonfill}}

\def\opp#1{\mathrel{%
\setbox7=\hbox{F}\setbox6=\hbox{%
\setbox0=\hbox{\footnotesize $#1$}\setbox1=\hbox{$\to$}%
\dimen@\wd0\advance\dimen@ by 0.66\wd1\relax
$\stackrel{\lrhto{\dimen@}}{\copy0}$}%
\ifdim\ht6>\ht7\dimen@\ht7\advance\dimen@ by -\ht6\else
\dimen@\z@\fi\raise\dimen@\box6}}

\def\edge#1#2{\picline{#1}{#2}\pt{#1}\pt{#2}}
\def\chrd#1#2{\picline{1 #1 polar}{1 #2 polar}}
\def\arrow#1#2{\picvecline{1 #1 polar}{1 #2 polar}}

\def\labch#1#2#3{\chrd{#1}{#2}\picputtext{1.5 #2 polar}{$#3$}}
\def\labar#1#2#3{\arrow{#1}{#2}\picputtext{1.5 #2 polar}{$#3$}}
\def\labbr#1#2#3{\arrow{#1}{#2}\picputtext{1.5 #1 polar}{$#3$}}

\def\CD{\szCD{4mm}}
\def\szCD#1#2{{\let\@nomath\@gobble\small\diag{#1}{2.4}{2.4}{
  \pictranslate{1.2 1.2}{
    \piccircle{0 0}{1}{}
    \picveclength{0.4}
    \picvecwidth{0.1}
    #2
}}}}

\def\labline#1#2#3#4{\picvecline{#1}{#2}\pictranslate{#2}{
  \picputtext{#3}{$#4$}}}
\def\lablineb#1#2#3#4{\picvecline{#1}{#2}\pictranslate{#1}{
  \picputtext{#3}{$#4$}}}

\def\pt#1{{\picfillgraycol{0}\picfilledcircle{#1}{0.06}{}}}
\def\labpt#1#2#3{\pictranslate{#1}{\pt{0 0}\picputtext{#2}{$#3$}}}

\def\GD#1{{\let\@nomath\@gobble\scriptsize\diag{6mm}{3.0}{3.0}{
  \pictranslate{1.5 1.5}{
    \piccircle{0 0}{1}{}
    #1
}}}}

\def\rottab#1#2{%
\expandafter\advance\csname c@table\endcsname by -1\relax
\centerline{%
\rbox{\centerline{\vbox{\setbox1=\hbox{#1}%
\hbox to \wd1{\hfill\vbox{{%
\caption{#2}}}\hfill}%
\vskip6mm
\copy1}}%
}%
}%
}


\author{A. Stoimenow\footnotemark[1]\\[2mm]
\small Department of Mathematics, \\
\small University of Toronto,\\
\small Canada M5S 3G3\\
\small e-mail: {\tt stoimeno@math.toronto.edu}\\
\small WWW: {\hbox{\tt http://www.math.toronto.edu/stoimeno/}}
}

{\def\thefootnote{\fnsymbol{footnote}}
\footnotetext[1]{Supported by a DFG postdoc grant. 
}
}

\title{\large\bf \uppercase{Knots of genus two}\\[4mm]
{\small\it This is a preprint. I would be grateful
for any comments and corrections!}}

\date{\large Current version: \curv\ \ \ First version:
\makedate{29}{6}{1998}}

\maketitle

\makeatletter

\def\ul{\underline}
\let\ol\overline
\let\point\pt
\let\ay\asymp
\let\pa\partial
\let\al\alpha
\let\ap\alpha
\let\be\beta
\let\Dl\Delta
\let\Gm\Gamma
\let\gm\gamma
\let\de\delta
\let\dl\delta
\let\eps\epsilon
\let\lm\lambda
\let\Lm\Lambda
\let\sg\sigma
\let\vp\varphi
\let\om\omega
\let\diagram\diag
\let\es\enspace

\let\sm\setminus
\let\wt\widetilde
\let\tl\tilde
\def\ncap{\not\mathrel{\cap}}
\def\vol{\text{\rm vol}\,}
\def\sgn{\text{\rm sgn}\,}
\def\cf{\text{\rm cf}\,}
\def\md{\max\deg}
\def\mc{\max\cf}
\def\Lra{\Longrightarrow}
\def\lra{\longrightarrow}
\def\so{\Rightarrow}
\def\So{\Longrightarrow}
\let\ds\displaystyle
\def\is{\big|\,S\cap [-k,k]^{n}\,\big|}
\let\ndiv\nmid

\def\lb{\linebreak[0]}
\def\lz{\lb\verb}

\let\reference\ref

\long\def\@makecaption#1#2{%
   \vskip 10pt
   {\let\label\@gobble
   \let\ignorespaces\@empty
   \xdef\@tempt{#2}%
   }%
   \ea\@ifempty\ea{\@tempt}{%
   \setbox\@tempboxa\hbox{%
      \fignr#1#2}%
      }{%
   \setbox\@tempboxc\hbox{%
      {\fignr#1:}\capt\ #2}%
      }%
   \ifdim \wd\@tempboxc >\captionwidth {%
      \hbox{
      \@tempdima=\captionwidth\relax
      \advance\@tempdima by 0\@captionmargin\relax
      \parbox[t]{\@tempdima}{%
      \leftskip=\@captionmargin
      \rightskip=\@captionmargin
      \unhbox\@tempboxc\par}%
      }%
      }%
   \else
      \hbox to\captionwidth{\hfil\box\@tempboxc\hfil}%
   \fi}%
\def\fignr{\small\sffamily\bfseries}%
\def\capt{\small\sffamily}%
\newbox\@tempboxc

\newdimen\@captionmargin\@captionmargin2cm\relax
\newdimen\captionwidth\captionwidth\hsize\relax

\def\eqref#1{(\protect\ref{#1})}

\def\proof{\@ifnextchar[{\@proof}{\@proof[\unskip]}}
\def\@proof[#1]{\noindent{\bf Proof #1.}\enspace}

\def\hint{\noindent Hint: }
\def\problem{\noindent{\bf Problem.} }

\def\@mt#1{\ifmmode#1\else$#1$\fi}
\def\qed{\hfill\@mt{\Box}}
\def\qqed{\hfill\@mt{\Box\enspace\Box}}

\def\cM{{\cal M}}
\def\cU{{\cal U}}
\def\cB{{\cal B}}
\def\cC{{\cal C}}
\def\cV{{\cal V}}
\def\cP{{\cal P}}
\def\tP{{\tilde P}}
\def\tZ{{\tilde Z}}
\def\fg{{\frak g}}
\def\tr{\text{tr}}
\def\cZ{{\cal Z}}
\def\cD{{\cal D}}
\def\bR{{\Bbb R}}
\def\cE{{\cal E}}
\def\bZ{{\Bbb Z}}
\def\bN{{\Bbb N}}

\def\bysame{\same[\kern2cm]\,}

\def\br#1{\left\lfloor#1\right\rfloor}
\def\BR#1{\left\lceil#1\right\rceil}

\def\abstractname{}

\@addtoreset {footnote}{page}

\renewcommand{\section}{%
   \@startsection
         {section}{1}{\z@}{-1.5ex \@plus -1ex \@minus -.2ex}%
               {1ex \@plus.2ex}{\large\bf}%
}
\renewcommand{\@seccntformat}[1]{\csname the#1\endcsname .
\quad}

\def\bC{{\Bbb C}}
\def\bP{{\Bbb P}}



\def\@test#1#2#3#4{%
  \let\@tempa\go@
  \@tempdima#1\relax\@tempdimb#3\@tempdima\relax\@tempdima#4\unitxsize\relax
  \ifdim \@tempdimb>\z@\relax
    \ifdim \@tempdimb<#2%
      \def\@tempa{\@test{#1}{#2}}%
    \fi
  \fi
  \@tempa
}

\def\go@#1\@end{}
\newdimen\unitxsize
\newif\ifautoepsf\autoepsftrue

\unitxsize4cm\relax
\def\epsfsize#1#2{\epsfxsize\relax\ifautoepsf
  {\@test{#1}{#2}{0.1 }{4   }
		{0.2 }{3   }
		{0.3 }{2   }
		{0.4 }{1.7 }
		{0.5 }{1.5 }
		{0.6 }{1.4 }
		{0.7 }{1.3 }
		{0.8 }{1.2 }
		{0.9 }{1.1 }
		{1.1 }{1.  }
		{1.2 }{0.9 }
		{1.4 }{0.8 }
		{1.6 }{0.75}
		{2.  }{0.7 }
		{2.25}{0.6 }
		{3   }{0.55}
		{5   }{0.5 }
		{10  }{0.33}
		{-1  }{0.25}\@end
		\ea}\ea\epsfxsize\the\@tempdima\relax
		\fi
		}

\def\rr#1{{\catcode`\_=11\relax\epsffile{#1.eps}}}

\def\vis#1#2{\hbox{\begin{tabular}{c}\rr{k-#1-#2} \\ $#1_{#2}$\end{tabular} }}
\def\vis@#1#2{\hbox{\begin{tabular}{c}\rr{k-#1-#2} \end{tabular} }}


\makeatletter
\def\mybrace#1#2{\@tempdima#1em\relax
\advance\@tempdima by -1em\relax
\setbox\@tempboxa=\hbox{\raisebox{-0.5\@tempdima}{$\ds 
\left.\rule[0.5\@tempdima]{\z@
}{0.5\@tempdima}\right\} #2$}}\dp\@tempboxa=\z@
\box\@tempboxa}

\def\namedarrow#1{{\enspace
\setbox7=\hbox{F}\setbox6=\hbox{%
\setbox0=\hbox{\footnotesize $#1$}\setbox1=\hbox{$\to$}%
\dimen@\wd0\advance\dimen@ by 0.66\wd1\relax
$\stackrel{\rato{\dimen@}}{\copy0}$}%
\ifdim\ht6>\ht7\dimen@\ht7\advance\dimen@ by -\ht6\else
\dimen@\z@\fi\raise\dimen@\box6\enspace}}
  
\def\ssim{\stackrel{\ds \sim}{\vbox{\vskip-0.2em\hbox{$\scriptstyle *$}}}}

\let\old@tl\~\def\~{\raisebox{-0.8ex}{\tt\old@tl{}}}
\let\lra\longrightarrow
\def\bt{\bar t'_{2}}
\let\sm\setminus
\let\eps\varepsilon
\let\ex\exists
\let\x\ex
\let\fa\forall
\let\ps\supset

\def\rs#1{\raisebox{-0.4em}{$\big|_{#1}$}}


\def\@test#1#2#3#4{%
  \let\@tempa\go@
  \@tempdima#1\relax\@tempdimb#3\@tempdima\relax\@tempdima#4\unitxsize\relax
  \ifdim \@tempdimb>\z@\relax
    \ifdim \@tempdimb<#2%
      \def\@tempa{\@test{#1}{#2}}%
    \fi
  \fi
  \@tempa
}

%
	
\def\uu#1{\setbox\@tempboxa=\hbox{#1}\@tempdima=-0.5\ht\@tempboxa
\advance\@tempdima by 0.5em\raise\@tempdima\box\@tempboxa}

\def\epsfs{\@ifnextchar[{\@epsfs}{\@@epsfs}}
\def\@epsfs[#1]#2{{\ifautoepsf\unitxsize#1\relax\else
\epsfxsize#1\relax\fi\epsffile{#2.eps}}}
\def\@@epsfs#1{{\epsffile{#1.eps}}}
	
\def\eepsfs{\@ifnextchar[{\@eepsfs}{\@@eepsfs}}
\def\@eepsfs[#1]#2{\uu{\ifautoepsf\unitxsize#1\relax\else
\epsfxsize#1\relax\fi\epsffile{#2.eps}}}
\def\@@eepsfs#1{\uu{\epsffile{#1.eps}}}

\def\xepsfs#1#2{{\catcode`\_=11\relax\ifautoepsf\unitxsize#1\relax\else
\epsfxsize#1\relax\fi\epsffile{#2.eps}}}

{\let\@noitemerr\relax
\vskip-2.7em\kern0pt\begin{abstract}
\noindent{\bf Abstract.}\enspace
We classify all knot diagrams of genus two and three, and
give applications to positive, alternating and homogeneous
knots, including a classification of achiral genus 2 alternating
knots, slice or achiral 2-almost positive knots, a proof of
the 3- and 4-move conjectures, and the calculation of the maximal
hyperbolic volume for weak genus two knots. We also study the values
of the link polynomials at roots of unity, extending denseness results
of Jones. Using these values, examples of knots with unsharp
Morton (weak genus) inequality are found. Several results are
generalized to arbitrary weak genus.
\end{abstract}
}

{\parskip0.2mm\tableofcontents}
\vspace{7mm}

\section{Introduction}

The notion of a Seifert surface of a knot is classical
\cite{Seifert}. Seifert proved the existence of these surfaces
by an algorithm 
how to construct such a surface out of some diagram of the
knot. Briefly, the procedure is as follows (see \cite[\S 4.3]{Adams}
or \cite{Rolfsen}): smooth out all crossings of the diagram, plug in
discs into the resulting set of disjoint (Seifert)
circles and connect the circles along the crossings
by half-twisted bands. We will call the resulting surface
canonical Seifert surface (of this diagram) and its genus the
genus of the diagram. The canonical (or weak) genus
of a knot we call the minimal genus of all
its diagrams.

The weak genus appears in previous work of several authors,
mainly in the context of showing it being equal to the
classical Seifert genus for large classes of knots, see \cite{%
Cromwell} and \em{loc.\ cit.} However, in \cite{Morton}, Morton
showed that this is not true in general. Later, further 
examples have been constructed \cite{Moriah,Kobayashi}.

Motivated by Morton's striking observation, in \cite{gen1}
we started the study of the weak genus in its own right.
We gave a description of knot diagrams of genus one and
made some statements about the general case.

The present paper is a continuation of our work in \cite{gen1}, and
relies on similar ideas. 
Its motivation was the quest for more
interesting phenomena occurring with knot diagrams of (canonical)
genus higher than one. The genus one diagrams, examined
in \cite{gen1}, revealed to be a too narrow class for such phenomena.
In this paper we will study the weak genus in greater generality.
We will prove several new results about properties of knots with
arbitrary weak genus. In the cases of weak genus 2 and 3 we have
obtained a complete description of diagrams. Using this description,
we obtain computational examples and results, some of them solving
(at parts) several problems in previous papers of other authors.

For most practical applications,
it is useful to consider weak genus 2. We study it thus in detail.
All methods should work also for higher genera, but applying them in
practice seems hardly worthwhile, as the little qualitative renewment
this project promises is contraposed to an extremely rapid increase of
quantitative effort. Diagrams of genus two turned out to be attractive,
because their variety is on the one hand sufficient to exhibit
interesting phenomena and allowing to apply different types of
combinatorial arguments to prove properties of them and the knots
they represent, but on the other hand not too great to make
impossible argumentation by hand, or with a reasonable
amount of computer calculations. As we will see,
many of the theorems we will prove for weak genus 2 cannot
be any longer proved reasonably (at least with the same methods) for
weak genus 3, if they remain true at all.

We give a brief survey of the structure of the paper.

In \S\reference{S2} we prove our main result, theorem
\reference{thgen2}, the classification of diagrams of genus two.
It bases on a combination of computational and mathematical
arguments. The most subsequent sections are devoted to
applications of this classification.

In \S\reference{S3} we give asymptotical estimates for the number of
alternating and positive knots of genus two and given crossing number
and classify the achiral alternating ones.

In \S\reference{S4} we show non-homogeneity of 2 of the undecided
cases in \cite[appendix]{Cromwell}, following from the more
general fact that homogeneous genus two knots are positive or
alternating. 

In \S\reference{S5} and \S\reference{S6} we use the Gau\ss{} sum
inequalities of \cite{pos} in a combination with the
result of \S\reference{S2} to show how to classify all positive
diagrams of a positive genus two knot, on the simplest non-trivial
examples $7_3$ and $7_5$, and classify all $2$-almost positive
unknot diagrams, recovering a result announced by
Przytycki and Taniyama in \cite{PrzTan}, that the only non-trivial
achiral (resp.\ slice) 2-almost positive knot is $4_1$ (resp.\ $6_1$).

In \S\reference{S7} we prove that there is no almost positive
knot of genus one, and in \S\reference{S8} that any positive
knot of genus two has a positive diagram of minimal crossing number.
We also show an example of a knot of genus two which has a single
positive diagram.

Beside the results mentioned so far, which are direct applications
of the classification in theorem \reference{thgen2},
we will develop several new theoretical tools, valid
for arbitrary weak genus. Most of these tools can again be used
to study the genus two case in further detail. As such a tool,
%
most substantially we will deal with behaviour of the Jones and HOMFLY
polynomial in \S\reference{S9}. We show how unity root evaluations of
the polynomials give information on the weak genus, and use this tool
to exhibit the first examples of knots on which the weak genus
inequality of Morton \cite{Morton} is not sharp. We also give, as an
aside, using elementary complex analysis and complex Lie group theory,
generalizations of some denseness theorems of Jones in \cite{Jones2}
about the values at roots of unity of the Jones polynomial of
knots of small braid index. Unity root evaluations of the
Jones polynomial seem to have become recently of interest 
because of a variety of relations to quantum physics,
in particular the volume conjecture. See \cite{DLL}.

Since these unity root evaluations are closely related to the
Nakanishi-Przytycki $k$-moves, we give several applications to
these moves in \S\reference{SecQ}, in particular the proof of
the $3$- and $4$-move conjecture for weak genus two knots
in \S\reference{3move} and \S\reference{4move}. We also discuss
how the criteria using the Jones and HOMFLY polynomial, and the
examples they give rise to, can be complemented by the
Brandt-Lickorish-Millett-Ho polynomial $Q$.

A further theoretical result is an asymptotical estimate for the
quality of the Seifert algorithm in giving a minimal (genus)
surface in \S\reference{S11}.

In \S\reference{Shyp}, we consider the hyperbolic volume.
Brittenham \cite{Brittenham} used a similar approach to ours to
prove that the weak genus bounds the volume of a hyperbolic knot.
We will slightly improve Brittenham's estimate of the maximal
hyperbolic volume for given weak genus, and (numerically)
determine the exact maximum for weak genus $1$ and $2$.



At the end of the paper we present the description for knot diagrams
of genus three in \S\reference{S12}, solving completely the knots
undecided for homogeneity in Cromwell's tables \cite[appendix]
{Cromwell}.

In \S\reference{S13} we conclude with some questions, and
a counterexample to a conjecture of Cromwell \cite{Cromwell2}.

Although a part of the material presented here (in particular the
examples illustrating our theoretical results) uses some computer
calculations, we hope that it has been obtained (and hence is
verifiable) with reasonable effort. To facilitate this, we include
some details about the calculations.

Further applications of the classification of genus 2 diagrams
are given in several subsequent papers. For example, in \cite{2apos}
this classification is used to give a short proof of a result
announced in \cite{PrzTan}, that positive knots of genus at least
2 have $\sg\ge 4$ (which builds on the result for genus 2 stated here
in corollary \reference{crs4}), in \cite{apos} to give a specific
inequality between the Vassiliev invariant of degree 2 and the
crossing number of almost positive knots of genus 2, and in \cite{restr}
to generalize the classification of $k$-almost positive achiral knots
for the case $k=2$ (announced also in \cite{PrzTan} and given here as
proposition \reference{ppaa}) for alternating knots to $k\le 4$.
These results are hopefully sufficiently motivating our approach.

\subsection*{\bf Notation.} For a knot $K$ and a (knot) diagram $D$,
$c(D)$ denotes the crossing number of $D$, $c(K)$ the crossing 
number of $K$ (the minimal crossing number of all its diagrams), 
$w(D)$ the writhe of $D$, $w(K)$ the writhe of an alternating diagram 
of $K$, if $K$ is alternating (this is an invariant of $K$, see
\cite{Kauffman}), and $n(D)$ the number of Seifert circles of $D$.
$\sg$ denotes the signature of a knot, $u$ denotes its unknotting
number, $\tl g$ denotes its weak genus and $g$ its classical 
(Seifert) genus. $!K$ denotes the obverse (mirror image) of a
knot $K$. Often we will assume a diagram to be reduced without
each time pointing it out. It should be always clear from the
context, where this is the case.

$v_2$ denotes the Vassiliev knot invariant of degree $2$,
normalized to be zero on the unknot and one on the trefoil(s).
$v_3$ denotes the primitive Vassiliev invariant of degree $3$,
normalized to be $4$ on the positive (right-hand) trefoil.
As usual, $V$ denotes the Jones \cite{Jones}, $\Dl$ the Alexander
\cite{Alexander}, $\nabla$ the Conway \cite{Conway}, $Q$
the Brandt-Lickorish-Millett-Ho \cite{BLM,Ho}, and $P$ the HOMFLY
(or skein) \cite{HOMFLY} polynomial. For the HOMFLY polynomial, we
use the variable convention of \cite{LickMil}.

For some polynomial $Y$ and some integer $k$ we denote by
$[Y(x)]_{x^k}$ the coefficient of $x^k$ in $Y(x)$. The minimal
(resp.\ maximal) degree of $Y$ we call the minimal (resp.\ 
maximal) $k$ with $[Y(x)]_{x^k}\ne 0$ and denote it by
$\min\deg_x Y$ (resp. $\max\deg_x Y$). The span of $Y$ is
the difference between its maximal and minimal degrees. In case
$Y$ has only one variable, its indication in the notation will
be omitted. The encoded notation for polynomials we use is the one of
\cite{poly}: if the absolute term occurs between the minimal
and maximal degrees, then it is bracketed, else the minimal degree is
recorded in braces before the coefficient list.

We use the notation of \cite{Rolfsen} for knots with
up to 10 crossings, renumbering $10_{163}\,\dots
10_{166}$ by eliminating $10_{162}$, the Perko duplication
of $10_{161}$, as has been done in the tables of \cite{BurZie}.
The notation of \cite{KnotScape} is used 
for knots from 11 crossings on. (Note, that for 11
crossing knots this notation \em{differs} from this of
\cite{Conway} and \cite{Perko}.) We use the convention
of the Rolfsen pictures to distinguish between the knot
and its obverse whenever necessary.

For two sequences of positive integers $(a_n)_{n=1}^{\infty}$
and $(b_n)_{n=1}^{\infty}$ we say that $a_n$ is
$O(b_n)$ iff $\limsup\limits_{n\to\infty}\myfrac{a_n}{b_n}<\infty$,
$O^\succeq(b_n)$ iff $\liminf\limits_{n\to\infty}\myfrac{a_n}{b_n}>0$,
and $O^\asymp(b_n)$ iff it is both $O(b_n)$ and $O^\succeq(b_n)$.

$\bZ$, $\bN$, $\bN_+$, $\bR$ and $\bC$ denote the integer, natural,
positive natural, real and complex numbers respectively.

For a set $S$, the expressions $|S|$ and $\#S$ are equivalent and both
denote the cardinality of $S$. In the sequel the symbol '$\subset $'
denotes a not necessarily proper inclusion.

\section{Knot diagrams of canonical genus 2\label{S2}}

It is known that a Seifert surface obtained by applying
Seifert's algorithm on a knot diagram $D$ has genus
\[
g(D)\,=\,\frac{c(D)-n(D)+1}{2}\,.
\]
This formula is shown by homotopy retracting the surface to a graph and
determining its Euler characteristic by a simple vertex and edge count.
The \em{weak} (or \em{canonical}) \em{genus} $\tl g(K)$ of a knot $K$
is defined as
\[
\tl g(K)\,:=\,\min\,\{\,g(D)\,:\,\mbox{$D$ is a diagram of $K$}\,\}
\,.
\]

In the following we will describe all knot diagrams of genus $2$
and deduce consequences for knots of weak genus two from this
classification.

As a preparation, we (re)introduce some terminology, recalling
\em{inter alia} some of the definitions and facts of \cite{gen1};
more details may be found there.

First we need to introduce some transformations of diagrams
which we will crucially need later.

In 1992, Menasco and Thistlethwaite \cite{MenThis} proved the
(previously long conjectured) statement, that reduced alternating
diagrams of the same knot (or link) must be transformable
by \em{flypes}, where a flype is shown of figure \ref{fig1}.

\begin{figure}[htb]
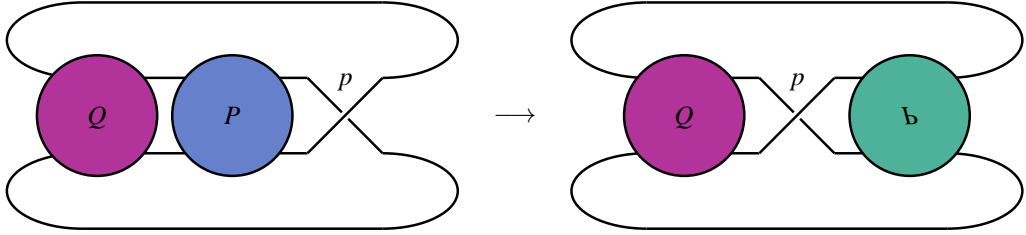

\[
\diag{1cm}{6}{3}{
    \pictranslate{3 1.5}{
       \picputtext[d]{1.5 0.4}{$p$}
       \picmultigraphics[S]{2}{1 -1}{
           \picmultiline{0.12 1 -1.0 0}{2 -0.5}{1 0.5}
           \picmultigraphics[S]{2}{-1 1}{
                \picellipsearc{-2 -1.0}{1 0.5}{90 270}
           }
           \picline{-2 -1.5}{2 -1.5}
           \picline{-2 -0.5}{1 -0.5}
      }
      \picfillgraycol{0.4 0.5 0.8}
      \picfilledcircle{0.0 0}{0.8}{$P$}
      \picfillgraycol{0.7 0.2 0.6}
      \picfilledcircle{-1.8 0}{0.8}{$Q$}
   }
}\quad\lra\quad
\diag{1cm}{6}{3}{
    \pictranslate{3 1.5}{
       \picputtext[d]{0 0.4}{$p$}
       \picmultigraphics[S]{2}{1 -1}{
           \picmultiline{0.12 1 -1.0 0}{0.5 -0.5}{-0.5 0.5}
           \picmultigraphics[S]{2}{-1 1}{
                \picellipsearc{-2 -1.0}{1 0.5}{90 270}
                \picline{-2 -0.5}{-0.5 -0.5}
           }
           \picline{-2 -1.5}{2 -1.5}
      }
      \picscale{1 -1}{
           \picfillgraycol{0.3 0.7 0.6}
           \picfilledcircle{1.5 0}{0.8}{$P$}
      }
      \picfillgraycol{0.7 0.2 0.6}
      \picfilledcircle{-1.5 0}{0.8}{$Q$}
   }
}
\]
\caption{A flype near the crossing $p$\label{fig1}}
\end{figure}

\begin{figure}[htb]
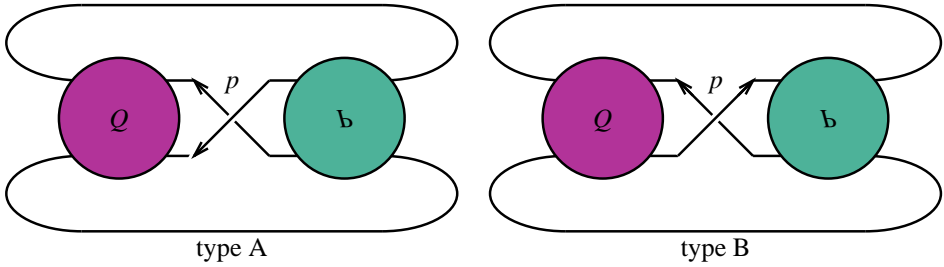

\[
\begin{array}{cc}
\diag{1cm}{6}{3}{
    \pictranslate{3 1.5}{
       \picputtext[d]{0 0.4}{$p$}
       \picmultigraphics[S]{2}{1 -1}{
           \picmultivecline{0.12 1 -1.0 0}{0.5 -0.5}{-0.5 0.5}
           \picmultigraphics[S]{2}{-1 1}{
                \picellipsearc{-2 -1.0}{1 0.5}{90 270}
                \picline{-2 -0.5}{-0.5 -0.5}
           }
           \picline{-2 -1.5}{2 -1.5}
      }
      \picscale{1 -1}{
           \picfillgraycol{0.3 0.7 0.6}
           \picfilledcircle{1.5 0}{0.8}{$P$}
      }
      \picfillgraycol{0.7 0.2 0.6}
      \picfilledcircle{-1.5 0}{0.8}{$Q$}
   }
}
&
\diag{1cm}{6}{3}{
    \pictranslate{3 1.5}{
       \picputtext[d]{0 0.4}{$p$}
       \picmultigraphics[S]{2}{1 -1}{
           \picmultiline{0.12 1 -1.0 0}{0.5 -0.5}{-0.5 0.5}
           \picmultigraphics[S]{2}{-1 1}{
                \picellipsearc{-2 -1.0}{1 0.5}{90 270}
                \picline{-2 -0.5}{-0.5 -0.5}
           }
           \picline{-2 -1.5}{2 -1.5}
       }
       \picvecline{-0.3 0.3}{-0.5 0.5}
       \picvecline{0.3 0.3}{0.5 0.5}
      \picscale{1 -1}{
           \picfillgraycol{0.3 0.7 0.6}
           \picfilledcircle{1.5 0}{0.8}{$P$}
      }
      \picfillgraycol{0.7 0.2 0.6}
      \picfilledcircle{-1.5 0}{0.8}{$Q$}
   }
}\\[3mm]
\mbox{type A} & \mbox{type B} 
\end{array}
\]
\caption{A flype of type A and B\label{ffl2}}
\end{figure}

The tangle $P$ on figure \reference{fig1} we call \em{flypable}, and we
say that the crossing $p$ \em{admits a flype} or that \em{the diagram
admits a flype at (or near) $p$}.

According to the orientation near $p$ we distinguish two
types of flypes, see figure \ref{ffl2}.

A \em{clasp} (we call it alternatively also a \em{matched
crossing pair}) is a tangle of the form
\[
\begin{tabular}{c@{\quad}c}
\diag{1cm}{2}{1}{
  \picrotate{-90}{
    \lbraid{-0.5 0.5}{1 1}
    \lbraid{-0.5 1.5}{1 1}
    \picvecline{-0.95 1.9}{-1 2}
    \picvecline{-0.05 0.1}{0 0}
  }
} &
\diag{1cm}{2}{1}{
  \picrotate{-90}{
    \lbraid{-0.5 0.5}{1 1}
    \lbraid{-0.5 1.5}{1 1}
    \picvecline{-0.95 1.9}{-1 2}
    \picvecline{-0.05 1.9}{0 2}
  }
} \\[8mm]
reverse clasp & parallel clasp
\end{tabular}\,,
\]
distinguished into \em{reverse} and \em{parallel} clasp
depending on the strand orientation.

By switching one of the crossings in a clasp and applying
a Reidemeister II move, one can eliminate both crossings.
This procedure is called \em{resolving} a clasp.
For the discussion below it is important to remark how
resolving a clasp affects the genus of the diagram. It reduces
the genus by one, if the clasp is parallel, or if it is reverse 
and the Seifert circles on which the two clasped strands lie
after the resolution are distinct. In this case we will
call the clasp \em{genus reducing}. Contrarily, a clasp resolution
preserves the the genus of the diagram if the
clasp is reverse and the strands obtained after the resolution
belong to the same Seifert circle (as for
example in the $\bt$ move we will just introduce). Then we call
the clasp \em{genus preserving}.

We will also need a class of diagram moves studied by
Przytycki and Nakanishi.

\begin{defi}(see \cite{Przytycki})
A $t_k$ move is a local diagram move replacing a parallel pair
of strands by $k$ parallel half-twists. Similarly, a $\bar t_{k}$
move for $k$ even is a replacement of a reversely oriented
pair of strands by $k$ reversely oriented half-twists. 
\end{defi}

A $\bar t_{2}$ move is thus replacing a reversely oriented
pair of strands by a reverse clasp. Of particular
importance will be, as in \cite{gen1}, a special
instance of a $\bar t_{2}$ move.

\begin{defi}\label{def2.2}
A \em{$\bt$ move} we call a $\bar t_2$ move \cite{Przytycki}
applied near a crossing
\[
\diag{7mm}{1}{1}{
    \picmultivecline{0.12 1 -1.0 0}{0 0}{1 1}
    \picmultivecline{0.12 1 -1.0 0}{1 0}{0 1}
}\namedarrow{\bt}
\diag{7mm}{3}{2}{
  \picPSgraphics{0 setlinecap}
  \pictranslate{1.0 1}{
    \picrotate{-90}{
      \lbraid{0 -0.5}{1 1}
      \lbraid{0 0.5}{1 1}
      \lbraid{0 1.5}{1 1}
      \pictranslate{-0.5 0}{
      \picvecline{0.03 1.95}{0 2}
      \picvecline{0.03 -.95}{0 -1}
    }
    }
    }
}
\,,
\]
and \em{reducing} $\bt$ move the reverse operation of a $\bt$ move.
We call a diagram \em{$\bt$ irreducible} if there is no sequence
of type B flypes transforming it into a diagram, on which a
reducing $\bt$ move can be applied. Let $c_g$ denote the \em{maximal
crossing number of an alternating $\bt$ irreducible genus $g$ diagram}.
\end{defi}

A flype of type A never creates or destroys a fragment obtained from
a crossing by a $\bt$ move and commutes with type B flypes, hence 
the applicability of a reducing $\bt$ move after type B flypes is
independent of type A flypes. In terms of the associated
Gau\ss{} diagram \cite{ninv,VirPol}, a knot
diagram is (modulo crossing changes) $\bt$ reducible after type B
flypes iff it has three chords, which do not mutually intersect and
all intersect the same set of other chords.
  
In order to discard uninteresting cases, we will consider mainly
only prime diagrams. 

\begin{defi}
A diagram $D$ is called \em{composite}, if there is a closed curve
$\gm$ intersecting (transversely) the curve of $D$ in two
points, such that both in- and exterior of $\gm$ contain crossings
of $D$. Else $D$ is called \em{prime} or \em{connected}.
\end{defi}

It is a simple observation that $c_0=0$. Two results of \cite{gen1} were
$c_1=4$ (independently observed by Lee Rudolph) and $c_g\le 8c_{g-1}+6$,
so that $c_g=O(8^g)$. However, it was evident, that this bound is
far from sharp, and later we showed in \cite{STV} that $c_g\le 12g-6$.
The starting point for a significant part of the material that follows
is to obtain for $g=2$ a more precise description.

\begin{figure}
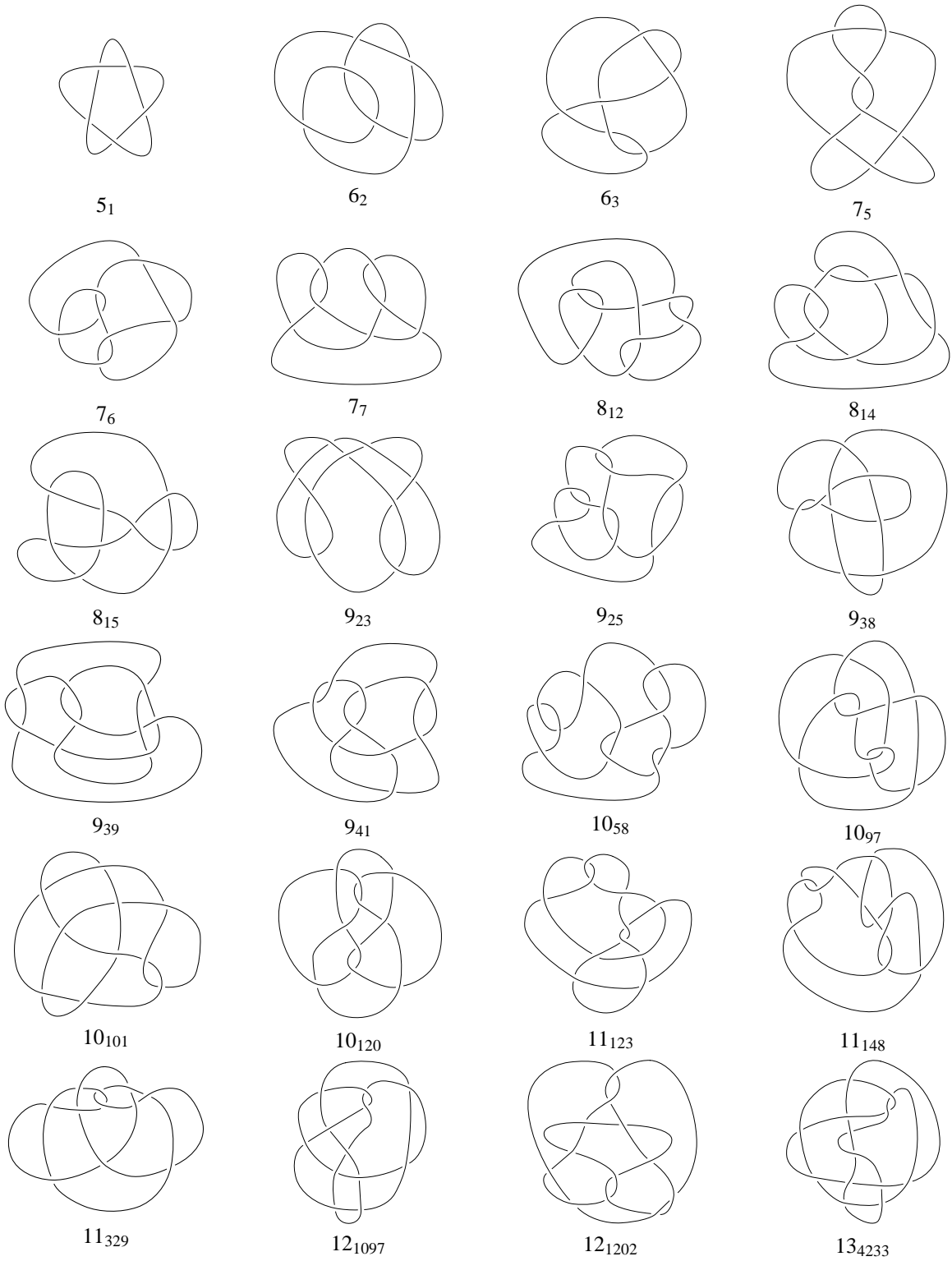

\[
\unitxsize2.8cm
\begin{array}{*4c}
\vis{5}{1} &
\vis{6}{2} &
\vis{6}{3} &
\vis{7}{5} \\
\vis{7}{6} &
\vis{7}{7} &
\vis{8}{12} &
\vis{8}{14} \\
\vis{8}{15} &
\vis{9}{23} &
\vis{9}{25} &
\vis{9}{38} \\
\vis{9}{39} &
\vis{9}{41} &
\vis{10}{58} &
\vis{10}{97} \\
\vis{10}{101} &
\vis{10}{120} &
\vis{11}{123} &
\vis{11}{148} \\
\vis{11}{329} &
\vis{12}{1097} &
\vis{12}{1202} &
\vis{13}{4233} \\
\end{array}
\]
\caption{\label{fig24}The 24 alternating genus 2 knots without an
alternating $\bt$ reducible diagram.}
\end{figure}

\begin{theo}\label{thgen2}
Let $K$ be a weak genus $2$ knot. Then any prime genus $2$ diagram
of $K$ is transformable by type B flypes into one which can be
obtained by crossing changes and $\bt$ moves from an alternating 
diagram of one of the 24 knots in figure \ref{fig24}.
\end{theo}

We will say that a diagram \em{generates} a \em{series} or
a \em{$\bt$ twist sequence} of diagrams by crossing changes and
$\bt$ moves (so that a $\bt$ twist sequence is a special case of
what was called in \cite{Bseq} a ``braiding sequence'').
In this terminology the classification result of genus one
diagrams in \cite{gen1} says that the only genus one generators
are (the alternating diagrams of) $3_1$ and $4_1$. 
Although we point out that some knots of figure
\ref{fig1} occur in multiple diagrams, it will be sometimes
possible and convenient to identify the series generated
by all its diagrams and call them a series \em{generated by the knot}.
For technical reasons (to have a numbering of the crossings)
it will turn out useful to record and fix a Dowker notation
\cite{DT} for each of these knots. (This is the one in the tables
of \cite{KnotScape}.)

\begin{kentry}
5 1:  6  8 10  2  4 
6 2:  4  8 10 12  2  6 
6 3:  4  8 10  2 12  6 
7 5:  4 10 12 14  2  8  6 
7 6:  4  8 12  2 14  6 10 
7 7:  4  8 10 12  2 14  6 
8 12:  4  8 14 10  2 16  6 12 
8 14:  4  8 10 14  2 16  6 12 
8 15:  4  8 12  2 14  6 16 10 
9 23:  4 10 12 16  2  8 18  6 14 
9 25:  4  8 12  2 16  6 18 10 14 
9 38:  6 10 14 18  4 16  2  8 12 
9 39:  6 10 14 18 16  2  8  4 12 
9 41:  6 10 14 12 16  2 18  4  8 
10 58:  4  8 14 10  2 18  6 20 12 16 
10 97:  4  8 12 18  2 16 20  6 10 14 
10 101:  4 10 14 18  2 16  6 20  8 12 
10 120:  6 10 18 12  4 16 20  8  2 14 
11 123:  4 10 14 20  2  8 18 22  6 12 16 
11 148:  4 10 16 20 12  2 18  6 22  8 14 
11 329:  6 12 18 22 14  4 20  8  2 10 16 
12 1097:  6 12 20 14 22  4 18 24  8  2 10 16 
12 1202:  6 20 10 24 14  4 18  8 22 12  2 16 
13 4233:  6 12 22 26 16  4 20 24  8 14  2 10 18
\end{kentry}

\proof[of theorem \reference{thgen2}]
By \cite{STV} any genus $2$ diagram of a 
weak genus $2$ knot can be obtained modulo type B flypes
by crossing changes and $\bt$ moves from an alternating 
diagram with at most $18$ crossings. Now the
24 knots in figure \ref{fig24} have been obtained 
by checking Thistlethwaite's tables of $\le 15$
crossing knots for $\bt$ irreducible alternating
genus $2$ diagrams.

It is important to note, that for each alternating
knot either all or no alternating diagrams are
$\bt$ irreducible modulo flypes. This follows
from the Menasco-Thistlethwaite flyping theorem
\cite{MenThis}, the fact that the applicability of
a reducing $\bt$ move is preserved by type A flypes
and that type A flypes and type B flypes commute
(i.e., if we can apply a type A flype and then a
type B flype, we can do so vice versa with the same
result). Hence it suffices to check the one specific
diagram included in the tables to figure out whether
the knot has a $\bt$ irreducible diagram.

It would be in principle possible to deal with the crossing numbers
$16$ to $18$ also by computer, but these tables are
not yet available to me (those of 16 crossings at least
at the time of the original writing), and to save a fair amount
of electronic capacity, it is preferable to
use mathematical arguments instead.
%
Let us give the following

\begin{lemma}\label{Lmm1}
If there is a $\bt$ irreducible alternating
genus $2$ diagram $D$ of $c$ crossings with a matched
crossing pair (clasp), then there is a $\bt$ irreducible
genus $2$ diagram of $c-2$ crossings, or $c\le 12$.
\end{lemma}

For the proof we need to make some definitions.

{
\let\Bt\bt
\let\bt\beta

\begin{defi}
A \em{region} of a knot diagram is a connected component of the complement
of its underlying curve in the plane. Every crossing $p$ is
bordered to four (not necessarily distinct) regions. We call two
of them $\al$ and $\bt$ \em{opposite} at $p$, notationally $\ap\opp{p} \bt$,
if they do not bound a common line segment (edge) in a neighborhood of $p$.
\[
\diag{1cm}{2}{2}{
  \picline{0 1}{2 1}
  \picmultiline{0.12 1 -1.0 0}{0 1 x}{2 1 x}
  \picputtext{0.8 1.25}{$p$}
  \picputtext{1.5 1.5}{$\bt$}
  \picputtext{0.5 0.5}{$\al$}
}
\]
One can see that if two of the four regions bordering a crossing
are equal, then they are opposite. In this case we call the
crossing \em{reducible} or \em{nugatory}, or an \em{isthmus}.
\end{defi}

\begin{defi}\label{dffr}
We call two crossings $p$ and $q$ of a knot diagram \em{linked},
notationally $p\cap q$, if the crossing strands are passed in
cyclic order $pqpq$ along the solid line, and unlinked if the cyclic
order is $ppqq$. Call two crossings $p$ and $q$ \em{equivalent},
if they are linked with the same set of other chords, that is if
$\fa c\ne p,q:\,c\cap p\iff c\cap q$. Call $p$ and $q$ \em{$\sim$-%
equivalent} $p\sim q$, if they are equivalent and unlinked and
\em{$\ssim$-equivalent} $p\ssim q$, if they are equivalent and linked.
\end{defi}

It is an exercise to check that $\sim$-equivalence and $\ssim$-%
equivalence are indeed equivalence relations and that two crossings
are $\sim$-\ (resp.\ $\ssim$-) equivalent if and only if after
a sequence of flypes they can be made to form a reverse (resp.\ 
parallel) clasp.

\proof[of lemma \reference{Lmm1}] Distinguish two cases for the
matched crossing pair in $D$.
{\def\theenumi{(\@roman\c@enumi)}
\def\labelenumi{\theenumi}
\begin{enumerate}
\item Strands are reverse and belong to distinct
Seifert circles. Then annihelating the matched crossing pair
gives a $c-2$ crossing alternating diagram $D'$. We claim
that (a) $D'$ is of genus $2$, and that (b)
it has no $\Bt$ reducible crossings.

The reason is that creating a situation of being able to perform a
$\Bt$ move after elimination of the matched pair always forces
the strands in the matched pair to belong to the same
Seifert circle (see figure \ref{fcl}). Namely, if after resolving
$3$ crossings $a$, $b$ and $c$ become $\sim$-equivalent, then there
are two regions $\al$ and $\bt$ of $D$, such that $\al\opp{p} \bt$
for an $p\in\{a,b,c\}$. Resolving the clasp joins two regions $\bt_1$
and $\bt_2$ of $D$ to one region $\bt$ of $D'$:
\[
\diag{1cm}{2}{3}{
  \picellipsearc{-0.1 1.5}{1.2 0.8}{-70 70}
  \picellipsearc{2.1 1.5}{1.2 0.8}{110 250}
  \picputtext{1 2.5}{$\bt_1$}
  \picputtext{1 0.5}{$\bt_2$}
}\quad\lra\quad
\diag{1cm}{2}{3}{
  \piccirclearc{-0.1 1.5}{0.8}{-60 60}
  \piccirclearc{2.1 1.5}{0.8}{120 240}
  \picputtext{1 1.5}{$\bt$}
}\,.
\]
Therefore, as $a$, $b$, and $c$ are not all $\sim$-equivalent in $D$,
w.l.o.g. $\al\opp{a} \bt_1$ and $\al\opp{b} \bt_2$ in $D$. But then
there exists in $D$ a dashed arc $\gm$ as in figure \ref{fcl}.
Then all Seifert circles on $D$ different form $k$, the Seifert circle
in the clasp, intersect the dashed curve $\gm$ totally only twice.
Thus both these crossings must belong to the same Seifert circle, and
hence resolving the clasp would be genus reducing.

\begin{figure}[htb]
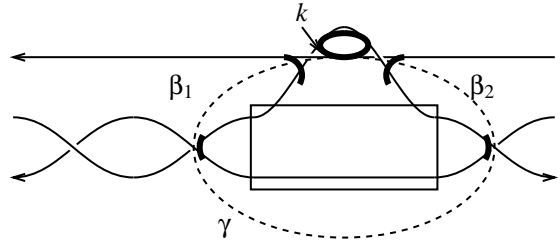

\[
\diag{8mm}{9}{4}{
 \pictranslate{9 1}{\picrotate{90}{
   \lbraid{0.5 6}{1 2}
   \lbraid{0.5 8}{1 2}
   \lbraid{0.5 1}{1 2}
   \picline{0 2}{0 5}
   \picline{2.0 0}{2.0 3.5}
   \picmulticurve{0.12 1 -1.0 0}{1 2}{1 2.5}{2.5 3}{2.5 3.5}
   \piccurve{2.5 3.5}{2.5 4}{1 4.5}{1 5}
   \picmultivecline{0.12 1 -1.0 0}{2 3.5}{2 9}
   \picvecline{0 8.95}{0 9}
   \picvecline{0 0.05}{0 0}
   {\piclinewidth{10}
    \piccirclearc{0.5 1.5}{0.4}{-120 -60}
    \piccirclearc{0.5 5.5}{0.4}{60 120}
    \piccirclearc{1.7 4.5}{0.3}{-120 0}
    \piccirclearc{1.7 2.5}{0.3}{0 120}
    \picellipse{2.2 3.5}{0.2 0.4}{}
   }
   \picvecline{2.5 4.2}{2.1 3.9}
   {\piclinedash{0.1}{0.05}
    \picellipse{0.5 3.5}{1.5 2.5}{}
   }
   \picbox{0.5 3.5}{1.4 3.1}{}
  }
 }
 \picputtext{3.5 0.2}{$\gm$}
 \picputtext{4.8 3.75}{$k$}
 \picputtext{2.8 2.55}{$\bt_1$}
 \picputtext{7.8 2.55}{$\bt_2$}
}
\]
\caption{When resolving a clasp makes a reducing $\Bt$ move 
applicable, the segments of the resolved clasp always belong to 
the same Seifert circle.\label{fcl}}
\end{figure}

Moreover, $D'$ has no reducible crossings. Assume that $p$ were such.
Then for some region $\al$ of $D'$ we have $\al\opp{p} \al$. But
then either $p$ is reducible in $D$, or $\al=\bt$ and
$\bt_1\opp{p} \bt_2$. Then we have a dashed curve $\gm$ like
\begin{eqn}\label{fgg}
\diag{1cm}{7}{4}{
  \pictranslate{0 4}{\picrotate{-90}{
    \lbraid{2 1}{1 2}
    \picline{2.5 2}{2.5 5.5}
    \picline{1.5 2}{1.5 5.5}
    \piccirclearc{2 5.5}{0.5}{0 180}
    \piccirclearc{2 6.2}{0.5}{180 360}
    \picline{2.5 6.2}{2.5 7}
    \picline{1.5 6.2}{1.5 7}
    { \piclinedash{0.1}{0.05}
      \picellipse{2 3.45}{1.3 2.45}{}
    }
    \picbox{2 3.5}{1.3 3.1}{}
    { \piclinewidth{10}
      \piccirclearc{2 1.5}{0.4}{-120 -60}
    }
 } 
 \picputtext{3.5 -1.0}{$\gm$}
 \picputtext{0.8 -2.5}{$p$}
}
 \picputtext{5.9 1.05}{$\bt_1$}
 \picputtext{5.9 2.95}{$\bt_2$}
}\,.
\end{eqn}
Then consider the Seifert circle in $D$ intersecting $\gm$ and apply
exactly the same argument as before to see that the clasp resolution
must be genus reducing.

\item Strands are parallel or belong to the same Seifert circle and are
reverse. Then annihelating the matched crossing pair reduces the
canonical genus of the diagram and we obtain a genus $1$ diagram $D'$.

We will show now that by $\Bt$ reducedness of $D$,
$D'$ has at most $4$ $\Bt$ reducible crossings. 
Thus $D'$ has at most $8$ crossings, and $D$ has at most $12$. 
To explain our argument again in more detail, we first need some definitions.

\begin{defi}
If $(a_1,\dots,a_n)$ is a finite sequence of objects, then 
$(a_{k_1},\dots,a_{k_l})$ is a subsequence if $k_i\ge k_{i-1}+1$, $k_1
\ge 1$ and $k_l\le n$, that is, the $a_{k_l}$s do not need to appear
immediately one after the other in $(a_1,\dots,a_n)$.
\end{defi}

\begin{defi}
Let $\al$ be a region of $D$, i.e. a connected component of the
complement of the plane curve of $D$ in the plane. Then consider
the sequence of regions
opposite to $\al$ at the crossings $\al$ borders in counterclockwise
order modulo cyclic permutation and call this \em{bordering sequence}
for $\al$ in $D$.
\[
\diag{5mm}{6}{6}{
  \pictranslate{3 3}{
    \picmultigraphics[rt]{5}{360 5 :}{
      \picline{3.4 -19 polar}{3.4 95 polar}
    }
    \picputtext{2.9 72 polar}{$\bt$}
    \picputtext{2.9 144 polar}{$\gm$}
    \picputtext{2.9 216 polar}{$\gm$}
    \picputtext{2.9 288 polar}{$\dl$}
    \picputtext{2.9 0 polar}{$\eps$}
    \picputtext{0 0}{$\al$}
  }
}\quad\lra\quad (\bt\gm\gm\dl\eps)
\]
\end{defi}

Note, that by connecting crossings with the same region $\gm$ opposite
to $\al$ by arcs in $\gm$ we see that the bordering sequence for $\al$
has no subsequence of the kind $\bt\gm\bt\gm$.

\begin{defi}
Call a set of crossings $\al_1,\dots,\al_n$ \em{mutually enclosed}
with respect to $\al$, if $\al_1,\dots,\al_n$ belong to the bordering sequence
for $\al$ and this bordering sequence can be cyclically permuted so that
the sequence $\al_1\al_2\dots\al_n\al_n\dots\al_2\al_1$ is a subsequence
of it.

The \em{enclosing index} $\eps_{\al,D}$ of $\al$ in $D$ is the
maximal size of a mutually enclosed set of crossings with respect to $\al$.
The enclosing index $\eps_D$ of $D$ is the maximal enclosing index
of all its regions.
\end{defi}

\begin{lemma}\label{Lmm2}
If we have a genus reducing clasp resolution $D\to D'$, joining
regions $\bt_1$ and $\bt_2$ of $D$ to $\bt$ of $D'$, and reduce by
Reidemeister I moves, flypes and reverse $\Bt$ moves
$D'$ to $D''$, then 
\[
c(D)-c(D'')\,\le\,4+4\eps_{\tl D}\,,
\]
where $\tl D$ is a diagram obtained by flypes from $D$.
\end{lemma}

\proof The absolute term `$4$', two of the crossings come from the clasp, and
two from the (Reidemeister I) reducible crossings in $D'$.

If there were three reducible crossings $a,b,c$ in $D'$ not
reducible in $D$, then $\bt_1\opp{p} \bt_2$ in $D$  for any
$p\in\{a,b,c\}$, and $a\sim b\sim c$ in $D$ (and not $a\ssim b
\ssim c$, as we can see form \eqref{fgg}), a contradiction
to its $\Bt$ irreducibility (see the remark after definition \ref{%
dffr}).

Separating $\bt$ in $D'$ into $\bt_1$ and $\bt_2$ in $D$ by reversing
the clasp resolution, enables us to add one $\Bt$ twist
to crossings participating in two mutually enclosed sets with respect to
$\bt$ in $D'$, leading to the term involving $\eps_{\tl D}$. \qed

Now, for any genus one diagram $D'$, $\eps_{D'}=1$, and using $c(D'')\le 4$
we obtain from lemma \ref{Lmm2} $c(D)\le 12$, concluding the
second case of the proof of lemma \ref{Lmm1}. \qed
\end{enumerate}
}

}

To show that there are no $\bt$ irreducible genus two diagrams
with $>13$ crossings we proceed by induction on the crossing number.

The cases of $14$ and $15$ crossings were excluded using
Thistlethwaite's tables (as I mentioned above). Then the
cases of $16$ and $17$ crossings can be (significantly)
reduced, applying lemma \reference{Lmm1}, to the cases with no
matched pair.

These cases we exclude as follows. Let $D$ be such a diagram
(that is, a genus 2 diagram with no matched pair).
A smoothing out of a crossing augments the number of $2$-gon
components of the diagram complement in the plane (or equivalently
the number of matched crossing pairs) by at most 2. Thus after
smoothing out a linked pair of crossings in $D$ we obtain a
diagram $D'$ of genus 1 with at most 4 matched pairs. Then
$D'$ is modulo its reducible crossings either a diagram obtained
from $3_1$ by at most two $\bt$ moves or one diagram obtained
from $4_1$ by at most one $\bt$ move.

Thus $D'$ has at most $7$ non-reducible crossings.
Now we count the reducible crossings of $D'$. 
(Compare to the proof of \cite[theorem 3.1]{gen1} or of
lemma \reference{Lmm2} above.) The smoothing out
of two crossings in $D$ identifies either two pairs or one
triple of regions. If $p$ is reducible in $D'$, then
$\beta_1\opp{p} \beta_2$ in $D$, where $\beta_{1,2}$ are among the
identified regions. There are two or three possible (unordered)
pairs $(\beta_1,\beta_2)$ of identified regions in $D$, and so there
are at most 4 or 6 crossings $p$ as above. Since two of these
crossings must be those smoothed out, $D'$ cannot have more than
$4$ reducible crossings.

We conclude that $D'$ must have at most $11$ crossings,
so $D$ has at most $13$ crossings. 

The same argument inductively excludes all higher crossing numbers. \qed


\begin{corr}
With $c_g$ as in definition \reference{def2.2},
we have $c_2=13$. 
\qed
\end{corr}

\begin{rem}
Note, that some of the 24 knots may have alternating
diagrams differing by a type A flype and twisting
at them gives mutated diagrams, the mutations
being type A ``flypes'' at a $\bt$ twisted crossing as shown
on figure \ref{figfl}. However, we can often ignore
these mutations, since for what we will do in the sequel
they will be mostly irrelevant.

For example, whenever we involve the Vassiliev invariants,
signature and knot polynomials in our proofs, the arguments
apply for all mutated diagrams as well, as these invariants
are preserved under mutation. (This is relevant for 
sections \reference{S4} to \reference{S8}.)
Also mutations do not occur in $\le 10$ crossings
(relevant for \S\reference{7_3_5} and \reference{S7}; the
few cases remaining can be checked directly), and 
rational knots and the unknot have no mutants \cite{HodgesonRubinstein}
(relevant for \S\reference{s5}).

Thus we consider only the one diagram for each of the 24 knots
given in figure \ref{fig24}. 
\end{rem}

\begin{figure}[htb]
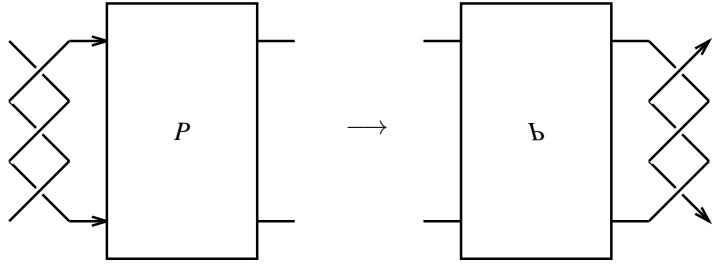

\[
\diag{1cm}{4}{3.5}{
   \picvecline{0.8 0.5}{1.3 0.5}
   \picvecline{0.8 2.9}{1.3 2.9}
   \picline{3.3 0.5}{3.8 0.5}
   \picline{3.3 2.9}{3.8 2.9}
   \picmultigraphics{3}{0 0.8}{
      \picline{0.8 0.5}{0 1.3}
      \picmultiline{0.12 1 -1.0 0}{0 0.5}{0.8 1.3}
   }
   \picfilledbox{2.3 1.7}{2 3.4}{$P$}
}\quad\lra\quad
\diag{1cm}{4}{3.5}{
   \pictranslate{1.5 1.7}{
      \picscale{1 -1}{
          \picfilledbox{0 0}{2 3.4}{$P$}
      }
   }
   \picline{0 0.5}{0.5 0.5}
   \picline{0 2.9}{0.5 2.9}
   \picline{2.5 0.5}{3 0.5}
   \picline{2.5 2.9}{3 2.9}
   \picmultigraphics{3}{0 0.8}{
      \picline{3.8 0.5}{3 1.3}
      \picmultiline{0.12 1 -1.0 0}{3 0.5}{3.8 1.3}
   }
   \picvecline{3.7 2.8}{3.8 2.9}
   \picvecline{3.7 0.6}{3.8 0.5}
}
\]
\caption{A ``flype'' near a $\bt$ twisted crossing is an iterated
mutation.\label{figfl}}
\end{figure}

\begin{rem}\label{rem1}
The present classification will later
be used to prove non-existence of minimal
canonical Seifert surfaces for some knots
of genus 2 when the obstruction of Morton
\cite{Morton} fails. (See remark \reference{remHS}.)
An explicit computer check gave
minimal canonical Seifert surfaces for all
knots up to 12 crossings (not only those of genus 2),
although not always minimal crossing number
diagrams suffice to give such a surface. (Among the
Rolfsen knots examples are the genus 3 knots $10_{155}$,
$10_{157}$ and $10_{159}$ and the genus 2
knots $10_{162}$ and $10_{164}$, where I found only
11 crossing diagrams doing the job; many more such
examples exist.) 
\end{rem}

In \cite{gen1} I showed that the number of knot
diagrams of given genus $g$ is polynomially bounded
in the crossing number. One sees that the maximal
exponent in this polynomial is $d_g-1$, where \em{$d_g$ is the 
maximal number of $\sim$-equivalence classes
in all diagrams of genus $g$}. For genus one
we had $d_1-1=2$ and for genus 2 we obtain $d_2-1=8$ for this
maximal exponent. The numbers $d_n$ seem not less important
then $c_n$ and will occur several times later.

\begin{corr}\label{on9}
The number of diagrams of genus 2 and crossing
number $n$ is $O^\asymp(n^8)$. Hence there are $O^\asymp(n^8)$
alternating genus 2 knots of crossing number $n$
and $O(n^9)$ positive knots of genus 2 or unknotting
number 2 and crossing number at most $n$.
\end{corr}

\proof[(to be continued)]
For the alternating case the only non-obvious point is to
show that there are $O^\asymp(n^8)$ alternating knots and not only
$O(n^8)$. I will give an argument for this at the end of \S\ref{S2}. 

The positive case is somewhat more involved as we do
not have the result of \cite{Kauffman2,Murasugi,Thistle} of minimality
(in crossing number) of alternating diagrams. Therefore we have
a result only for bounded but not fixed crossing number.
We also need to use that a positive genus 2 knot has a
positive diagram of minimal crossing number. This is again
not straightforward and will be proved in theorem \reference{thpos}.
The result for the unknotting number and positive knots follows from
the inequality $u\ge g$ \cite[corollary 4.3]{pos}. \qed

\section{Alternating genus two knots\label{alt2}\label{S3}}

The $\bt$ twist sequences  of some of the
24 knots contain those of some others as a
subfamily. This happens when resolving a clasp.
The relations are given in figure
\ref{figrel}. Therein the knots are encircled,
whose  twist sequences are not contained in
any other (we will call them \em{main}),
and for the others not all (but at least one) 
of the sequences containing them is indicated.

\begin{figure}
\[
\diag{1cm}{10}{10}{
   \picfilledellipse{2 9}{1 0.5}{$13_{4233}$}
   \picputtext{1 7}{$11_{123}$}
   \picputtext{3 7}{$11_{329}$}
   \picfilledellipse{5 7}{0.9 0.5}{$11_{148}$}
   \picputtext{1 5}{$9_{23}$}
   \picputtext{3 5}{$9_{38}$}
   \picputtext{5 5}{$9_{25}$}  
   \picputtext{7 5}{$9_{39}$}
   \picfilledellipse{9 5}{0.8 0.5}{$9_{41}$}
   \picputtext{2 3}{$7_5$}
   \picputtext{5 3}{$7_6$}
   \picputtext{9 3}{$7_7$}
   \picputtext{2 1}{$5_1$}
   \picvecline{5.5 6.4}{7.0 5.5}
   \picvecline{7.0 4.4}{5.5 3.5}
   \picvecline{1.7 8.4}{1.3 7.6}
   \picvecline{2.3 8.4}{2.7 7.6}
   \picvecline{1 6.4}{1 5.6}
   \picvecline{3 6.4}{3 5.6}
   \picvecline{2.7 6.4}{1.3 5.6}
   \picvecline{1.3 6.4}{2.7 5.6}
   \picvecline{5 6.4}{5 5.6}
   \picvecline{1.3 4.4}{1.7 3.6}
   \picvecline{2.7 4.4}{2.3 3.6}
   \picvecline{2 2.4}{2 1.6}
   \picvecline{5 4.4}{5 3.6}
   \picvecline{9 4.4}{9 3.6}
}
\]
\twrule
\[
\diag{1cm}{10}{8}{
   \picfilledellipse{2 7}{1 0.5}{$12_{1097}$}
   \picfilledellipse{5 7}{1 0.5}{$12_{1202}$}
   \picputtext{1 5}{$10_{120}$}
   \picputtext{3 5}{$10_{101}$}
   \picputtext{5 5}{$10_{58}$}
   \picfilledellipse{7 5}{0.9 0.5}{$10_{97}$}
   \picputtext{3 3}{$8_{15}$}
   \picputtext{5 3}{$8_{12}$}
   \picputtext{7 3}{$8_{14}$}  
   \picputtext{7 1}{$6_2$} 
   \picfilledellipse{9 1}{0.6 0.5}{$6_3$}
   \picvecline{1.7 6.4}{1.3 5.6}
   \picvecline{2.3 6.4}{2.7 5.6}
   \picvecline{5 6.4}{5 5.6}
   \picvecline{1.3 4.6}{2.6 3.6}
   \picvecline{3 4.4}{3 3.6}
   \picvecline{5 4.4}{5 3.6}
   \picvecline{7 4.4}{7 3.6}
   \picvecline{7 2.4}{7 1.6}
}
\]
\caption{The inclusion relations, under resolving clasps,
between the twist sequences of the 24 generating knots and
the indication (by encircling) of the main twist sequences.%
\label{figrel}}
\end{figure}

\begin{rem}
It is striking and suggested by the figure that inclusions
of series occur only between generators of the same parity of the
crossing number. This will be so for higher genera diagrams, too.
As already remarked, whenever resolving a clasp simplifies
the diagram by more than the two crossings (by removing nugatory
crossings), the resulting diagram must have already smaller genus.
\end{rem}

We record two small consequences. First note that
$6_3$ is simple but main. Some reason for this is
that it is the only knot of the 24 where the numbers
of positive and negative crossings in the alternating
diagram are both odd. Therefore, we have

\begin{prop}
Let $K$ be an alternating genus 2 knot with
$\{c(K),w(K)\}\bmod 4=\{0,2\}$. Then $K$ is an
arborescent knot with Conway notation $(p,q)rs(t,u)$
with $p,q,r,s,t,u>0$ all odd. \qed
\end{prop}

Another interesting aspect is to consider the
achiral knots among the alternating genus 2 knots.
First we obtain

\begin{prop}
A prime alternating genus 2 knot $K$ has zero signature, if and
only if a diagram of $K$ can be obtained from a diagram of $6_3$,
$7_7$, $8_{12}$, $9_{41}$, $10_{58}$ or $12_{1202}$ by (repeated)
$\bt$ moves.
\end{prop}

\proof The one direction follows from computing the signatures of
the 24 knots and the fact, that a $\bt$ move in an 
alternating
diagram does not change the signature (which follows
from the Traczyk-Murasugi formula, see \cite{Traczyk}
or \cite[p 437]{Kauffman}). For the reverse direction
note that by a result
of Menasco \cite{Menasco} the primeness of an alternating
knot is equivalent to the primeness of (any)one of its 
alternating diagrams. \qed

\begin{corr}\label{cvg}
Let $K$ be a prime achiral alternating genus 2 knot.
Then a diagram of $K$ can be obtained 
from a diagram of $6_3$, $8_{12}$, $10_{58}$ or $12_{1202}$ by 
(repeated) $\bt$ moves.
\end{corr}

\proof This follows from the preceding proposition by
excluding the odd crossing number knots.  
\qed

It is, however, much more interesting to have
an exact classification of all such knots. This is
obtained by applying the flyping theorem of Menasco
and Thistlethwaite. (Here I for completeness include
the composite case.)

\begin{theo}
Let $K$ be an achiral alternating genus 2 knot.
Then a diagram of $K$ is either
\begin{enumerate}
\item a composite diagram 
\begin{enumerate}
\item $C(q,q)\# C(p,p)$ with $p,q>0$ even or
\item $K\# !K$ with $K\in\{\,C(p,q)\,|\,\mbox{%
$p,q>0$ even}\,\}\cup\{\,P(p,q,r)\,|\,\mbox{%
$p,q,r>0$ odd}\,\}$;
\end{enumerate}
\item an arborescent diagram with Conway notation
$(a,b)cc(a',b')$ with $a,b,c,a',b'>0$ odd and $\{a,b\}=
\{a',b'\}$ (in which case the knot is $+$achiral if
$a=a'$ and $-$achiral if $a=b'$),
\item a rational diagram $C(a,b,b,a)$ with $a,b>0$ even
(which is invertible so the knot is $+-$achiral) or
\item a diagram in the $\bt$ twist sequence of 
$12_{1202}$ with $a,b,c$ $\bt$ twists at the three positive
clasps and $a',b',c'$ twists at the three negative
clasps, such that $a,b,c\ge 0$ and $\{a,b,c\}=\{a',b',c'\}$
(in which case the knot is $+$achiral or $-$achiral depending
on whether the cyclic orderings of $(a,b,c)$ and $(a',b',c')$
along the knot are the same or reverse).
\end{enumerate}
\end{theo}

\proof In the case the knot $K$ is composite 
it must have two prime factors of genus one
and by a result of Menasco \cite{Menasco} both
are alternating. By the uniqueness of the
decomposition into prime factors, if $K$
is achiral both factors must be so or mutually obverse.
Now use the classification of alternating genus one knots
in \cite{gen1}. It is an easy consequence of this classification that
the only achiral knots among them are the rational knots $C(q,q)$ with
$q>0$ even. Then one obtains the above characterization.

In the case the knot $K$ is prime, using corollary \reference{cvg},
we need to discuss 4 cases.
\begin{itemize}
\item[$12_{1202}$:] It is easy to see (e.g., by looking at the
Gau\ss{} diagram \cite{ninv,VirPol} shown of figure \ref{ff}) 
that the diagram of $12_{1202}$
and any other diagram in its $\bt$ twist sequence does
not admit a flype. Hence the knot is achiral if and only
if the Gau\ss{} diagram is isomorphic to itself (or its
mirror image) with the signs of the crossings switched,
which happens exactly in the cases recorded above.
\item[$10_{58}$:] To show that we have no achiral knot here
we use the intersection graph of the Gau\ss{} diagram.
Its vertices correspond to the arrows in the Gau\ss{}
diagram and are equipped with the sign of the crossing
in the knot diagram. Two vertices $a$ and $b$ are connected by an
edge if and only if the arrows in the Gau\ss{} diagram
intersect (or the crossings are linked in the sense of
definition \reference{dffr}). A flype
preserves the intersection graph and hence the 
intersection graph of an achiral alternating knot diagram 
must have an automorphism reversing the signs of all
vertices. To see that any diagram in the $\bt$
twist sequence of $10_{58}$ does not have such
an automorphism, consider the equivalence relation
between vertices from definition \reference{dffr}.
Then the number of $\sim$- and $\ssim$-equivalence classes
of positive resp.\ negative crossings in each such
diagram is $2$ resp.\ $3$, and hence there cannot
be an automorphism of the desired kind.
\item [$8_{12}$:] Use again the intersection graph.
Looking at the number of positive and negative arrows
intersecting only one $\sim$- or $\ssim$-equivalence class
of arrows, we find that in the form $C(a,b,c,d)$
we must have $a=d$. Then $b=c$ follows from
looking at the number of positive and negative arrows
at all (or the writhe). This also follows from general
rational knot theory arguments.
\item[$6_3$:] The Gau\ss{} diagram is shown schematically
on figure \ref{ff}. Looking at the number of positive and 
negative arrows intersecting only ones of the same sign we find
$c=c'$, and hence by the writhe argument $a+b=a'+b'$.
Then counting the number of intersections between arrows
of the same sign we find $ab=a'b'$, whence 
$\{a,b\}=\{a',b'\}$. \qed
\end{itemize}

\begin{figure}
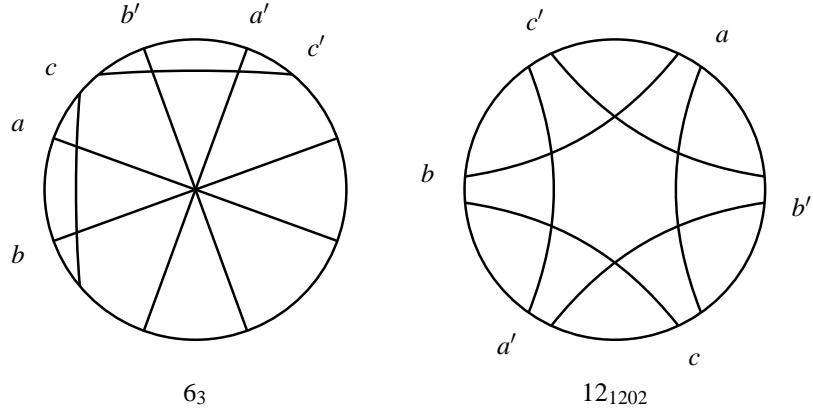

\[
\begin{array}{c@{\hskip1.5cm}c}
\diag{1cm}{4}{4}{
  \pictranslate{2 2}{
     \piccircle{0 0}{2}{}
     \picmultigraphics[rt]{2}{-90}{
        \picline{20 2 x polar}{200 2 x polar}
        \picline{-20 2 x polar}{160 2 x polar}
        \piccurve{140 2 x polar}{-1.6 0.5}{-1.6 -0.5}{220 2 x polar}
     }
     \picputtext{50 2.5 x polar}{$c'$}
     \picputtext{70 2.5 x polar}{$a'$}
     \picputtext{110 2.5 x polar}{$b'$}
     \picputtext{140 2.5 x polar}{$c$}
     \picputtext{160 2.5 x polar}{$a$}
     \picputtext{200 2.5 x polar}{$b$}
  }
} &
\diag{1cm}{4}{4}{
  \pictranslate{2 2}{
     \piccircle{0 0}{2}{}
     \picmultigraphics[rt]{6}{60}{
        \piccurve{55 2 x polar}{0.7 0.5}{0.7 -0.5}{-55 2 x polar}
     }
     \picputtext{140 25 - 2.5 x polar}{$c'$}
     \picputtext{260 25 - 2.5 x polar}{$a'$}
     \picputtext{20 25 - 2.5 x polar}{$b'$}
     \picputtext{320 25 - 2.5 x polar}{$c$}
     \picputtext{80 25 - 2.5 x polar}{$a$}
     \picputtext{200 25 - 2.5 x polar}{$b$}
  }
} \\[23mm]
6_3 & 12_{1202}
\end{array}
\]
\caption{Schematic drawing of the Gau\ss{} diagrams
in the $\bt$ twist sequence of $12_{1202}$ and $6_3$.
Orientation of the arrows is abused. A number like $a$ at
each chord denotes that it stands for a family of
$a$ neighbored non-intersecting chords. The crossings
are negative for the groups labeled by $a',b'$ and
$c'$ and positive for the groups labeled by $a,b$ and
$c$. For $12_{1202}$ all 6 numbers are even and for
$6_3$ odd.\label{ff}}
\end{figure}

\begin{rem} 
As far as orientation goes for the
composite case, the non-invertible genus 1 alternating
knots are $P(p,q,r)$ with $3\le p<q<r$ \cite{Trotter}.
So, taking one of these knots $K$, the knot $K\# !K$
is $+$achiral and $K\# -!K$ is $-$achiral. The rest
of the knots are invertible and so $+-$achiral.
\end{rem}

Using the intersection graph arguments we can now
easily complete the proof of corollary \ref{on9} in the
alternating case.

\proof[of corollary \ref{on9} (continued)] The only point is to
convince oneself that the $O^\asymp(n^8)$ alternating diagrams
remain at that quantity after modding out by flypes.
For this consider just diagrams, where the number of
$\bt$ moves applied to any $\sim$ equivalence class
of crossings in the diagram generating the series is
different and then there cannot be isomorphism of any 
two of the intersection graphs (just because the sets
of cardinalities of the $\sim$ equivalence classes
are never the same). But the number of compositions
of length $k$ of some number $n$ into strictly ascending
parts is the same as the number of compositions of
$n-\mybin{k}{2}$ into $k$ non-strictly ascending
parts (or the number of partitions of $n-\mybin{k}{2}$ 
of length $k$), which is $O^\asymp(n^{k-1})$. 

The proof of corollary \reference{on9} is now complete
modulo theorem \reference{thpos}. \qed

Considering the signature $\sg$, we mention a final consequence of theorem
\ref{thgen2} for positive knots, which also follows from \cite{PrzTan}.

\begin{corr}\label{crs4}
A positive genus two knot has $\sg=4$.
\end{corr}

\proof It is clear that $\sg\le 4$. To show $\sg=4$ it suffices to
check it on the (positively crossing switched) generating diagrams
on figure \ref{fig24}, as a $\bt$ move never reduces $\sg$. \qed

\section{Homogeneous genus two knots\label{S4}}

In \cite{Cromwell}, Cromwell introduced a certain class of
link diagrams he called homogeneous, which possess minimal
(genus) canonical Seifert surfaces. Roughly, a diagram is homogeneous,
if the connected components, called \em{blocks}, of the complement
of its Seifert picture (set of all Seifert circles lying in the
projection plane) contain only crossings of the same sign.
Letting this sign always remain the same or always change
when passing through a Seifert circle, we obtain the positive
(or negative) and alternating diagrams as special cases.
For five 10 crossing knots Cromwell could not
decide about the existence of a homogeneous diagram~-- $10_{144}$,  
$10_{151}$, $10_{158}$, $10_{160}$ and $10_{165}$. Two of
them have genus 2~-- $10_{144}$ and $10_{165}$. The present
discussion enables us to handle these cases.

\begin{theo}\label{tap}
Any homogeneous genus two knot $K$ is alternating or positive.
\end{theo}

Note, that this is no longer
true for genus three, as shows Cromwell's example $9_{43}$.

\begin{corr}\label{cor4.2}
The knots $10_{144}$ and $10_{165}$ are non-homogeneous.
\end{corr}

\proof The knots $10_{144}$ and $10_{165}$ violate obstructions 
to being positive (e.g. \cite[theorem 4(b)]{Cromwell} or \cite{pos})
or alternating (one edge coefficient of its Jones
polynomial is not $\pm 1$, see \cite{Kauffman2,Murasugi,Thistle}),
hence cannot be homogeneous. \qed

Before we start with the proof of theorem \reference{tap},
we need one more definition.

\begin{defi}
The \em{interior} of a Seifert circle is the bounded component
of its complement in the plane, and its \em{exterior} is the
unbounded one. The Seifert circle is called \em{separating},
if both its in- and exterior contain at least one other 
Seifert circle (or equivalently, at least one crossing),
and \em{non-separating} otherwise.
\end{defi}

First we record a statement we will use later to reduce the
number of cases to discuss.

\begin{lemma}\label{pq}
Let $D$ be an alternating diagram with (i) exactly three negative
crossings, all connecting a non-separating Seifert circle, or (ii)
with exactly two negative crossings. Assume furthermore that,
whatever case (i) or (ii) we are in, no flype can be
performed at any one of these two or three crossings. Then any
homogeneous diagram in the series of all its diagrams is
either positive or alternating.
\end{lemma}

\proof Assume a knot has in (all) its alternating diagram(s)
at most 3 negative (or positive) crossings. Then the fact
that alternating diagrams are homogeneous shows that any
Seifert circle must be connected from the same side 
by crossings of the same sign and then by the
non-existence of isthmus crossings any Seifert circle
is connected by either no or at least two negative crossings.
So, if they are at most three, all the negative crossings
connect the same pair of Seifert circles (there cannot be three
Seifert circles, each one connected with the other two, because of
orientation reasons). Then they belong to the same block.

If the crossings are three, one of the two Seifert circles 
to which they connect
has an empty interior (or exterior), and the diagram does not admit
a flype near one of these crossings, then 
(e.g. by looking at the chords
of the three crossings in the Gau\ss{} diagram) one convinces himself
that the triple of crossings is preserved by flypes,
and so the Seifert circle stays empty after any flype. 
Thus any alternating diagram of the knot has at most 
one separating Seifert circle,
and then each homogeneous diagram in
the series of this diagram is either positive or alternating. 

If the negative crossings are two and the diagram has two separating
Seifert circles, then these are exactly the Seifert circles connected by
the two negative crossings and both inside the inner one and outside
the outer one (or inside both if the one does not contain
the other) there are crossings. But then these negative
crossings admit a flype.
\qed
 
\proof[of theorem \reference{tap}] $g(K)=2$, so a homogeneous diagram
of $K$, if it exists, must lie in one of the 24 series (the composite
diagrams are connected sums of alternating pretzel diagrams, so
the claim is trivial for such diagrams). 


The series of $9_{38}$, $10_{101}$, $10_{120}$, $11_{123}$, $11_{329}$,
$12_{1097}$ and $13_{4233}$ one excludes by positivity (their
alternating diagrams are positive, and hence so is any homogeneous
diagram in their series).

Consider the series of $9_{39}$, $9_{41}$, $10_{97}$, $11_{148}$ and
$12_{1202}$. 
The diagram of $12_{1202}$ does not admit a flype (hence it is the only 
alternating diagram of $12_{1202}$) and it has exactly one separating
Seifert circle. $9_{39}$, $10_{97}$ and $11_{148}$ have two negative 
crossings which do not admit a flype. Finally, $9_{41}$ has three
negative crossings, all of which do not admit a flype and together
bound an empty Seifert circle. 
Then by the lemma each homogeneous diagram in
the series of all $5$ knots is either positive or alternating. 

There remain the 12 arborescent generating knots $5_1\,\dots 9_{25}$
and $10_{58}$. To handle these series, use the $\le 3$ negative
(or positive) crossing argument of lemma \reference{pq}. It works
except for $6_3$, $7_6$, $7_7$, $8_{12}$ and $10_{58}$. (Note, that
in most cases of two negative crossings they form a flypable clasp and
hence cannot admit a flype themselves.) 

$6_3$ one excludes because it has only three Seifert circles,
hence it cannot have two separating ones.

$8_{12}$ one excludes because it admits only
type B flypes and so the series of all its diagrams are equivalent,
but the one of $C(2,2,2,2)$ contains only rational knots, and
such knots are alternating.

Now, $10_{58}$ has an alternating diagram with $5$ clasps,
two of them negative (say, modulo mirroring). We find out
that the only possibility to flype is to flype the tangles
of these clasps, giving us (modulo symmetries) a total of
4 alternating diagrams of $10_{58}$. The only way to make
them homogeneous, but not positive and not alternating, is to
switch exactly one of the clasps in $3$ of these diagrams,
and then possibly to perform $\bt$ moves. As $10_{58}$'s
alternating diagrams differ only by type B flypes, it suffices
to consider one of these $3$ diagrams. But it is easy to see
that the diagram simplifies to an alternating one of one
crossing less.

$7_6$ one excludes similarly. We have the 2 negative crossings
admitting a flype, the flypable tangle being a positive clasp.
The proof of the lemma shows that the possibility to obtain modulo
flypes a homogeneous diagram is to switch or not the negative and/or
the flypable positive clasp. From the four cases only the two where
the flypable clasp is switched are neither alternating nor positive.
We end up with
\[
\diag{1cm}{4}{4}{
  \picputtext[dl]{0 0}{\epsfs{k-7-6.1}}
}
\quad\text{and}\quad
\diag{1cm}{4}{4}{
  \picputtext[dl]{0 0}{\epsfs{k-7-6.2}}
}\,.
\]
But in both cases one can see that after performing any series of
$\bt$ moves the diagram can be simplified to an alternating one.

$7_7$ one excludes the same way. The only way to obtain a
homogeneous non-positive and non-alternating diagram is to switch
exactly one of the two positive flypable clasps, but all diagrams
in this series simplify to an alternating one. 

In fact, one should be even a little more careful. Theorem 
\ref{thgen2} just said that one obtains a diagram in the series
modulo type B flypes (and a type B flype may change the
homogeneity of the diagram). But one can find out that the only
cases where the flype is necessary are to have 2 and 2
(for $7_6$ and $10_{58}$) and 2 and 1 (for $7_7$) flypable
crossings on both sides of the flypable negative clasp(s),
and these cases one handles exactly as above. \qed

\section{Classifying positive diagrams of some positive genus 2
knots\label{7_3_5}\label{S5}}

The strict increase of $v_2$ and $v_3$ under $\bt $ moves
at a positive diagram enables us to classify with reasonable effort
all positive diagrams of positive knots of genus 2 (or higher genera,
if an analogue of theorem \ref{thgen2} is worked out), if they are
not too complicated. We describe this procedure for the examples $7_3$
and $7_5$ (for which the use of $v_2$ suffices). The result is a
special case of a more general procedure, so the discussion aims
to show how in principle such a task can be solved.

Denote by $\ol{K}$ for an alternating knot $K$ the diagram obtained from
an alternating diagram of $K$ by making it positive by crossing changes
(this is defined up to flypes).

\begin{prop}\label{p7375}
The positive diagrams of $7_3$ are (up to flypes):
$\ol{7_3}$, $\ol{8_4}$, $\ol{8_{11}}$, $\ol{8_{13}}$, $\ol{9_{12}}$,
$\ol{9_{14}}$, $\ol{9_{21}}$, $\ol{9_{37}}$, and $\ol{10_{13}}$.
The positive diagrams of $7_5$ are: $\ol{7_5}$, $\ol{8_6}$, $\ol{8_8}$,
$\ol{8_{14}}$, $\ol{9_8}$, $\ol{9_{15}}$, $\ol{9_{19}}$, $\ol{10_{35}}$.
\end{prop}

\proof We have $v_2(7_5)=4$ and $v_2(7_3)=5$. 
Let $D$ be a positive diagram
of $7_3$ or $7_5$. Then $D$ belongs to the twist
sequence of one of the 24 knots above. In case of $8_{15}$,
$9_{23}$, $9_{38}$, $10_{101}$, $10_{120}$, $11_{123}$,
$11_{329}$, $12_{1097}$ and $13_{4233}$ the alternating
diagrams are positive and as doing $\bt$ moves does not
spoil alternation, all positive diagrams of their twist
sequence are alternating diagrams with at least $8$ crossings,
and hence by \cite{Kauffman2,Murasugi,Thistle} never belong to
$7_3$ or $7_5$. The same is true for the twist
sequence of $7_5$, with the exception that in it exactly
the diagram of $7_5$ belongs to itself and no one
belongs to $7_3$.

By an analogous argument the only diagram in the
twist sequence of $5_1$ belonging to $7_3$
is $7_3$'s usual $(1,1,1,1,3)$ pretzel diagram, and no
diagram belongs to $7_5$.

In the series of $9_{39}$, $9_{41}$, $10_{97}$, $11_{148}$ and
$12_{1202}$ the positive diagram obtained by crossing changes
from the alternating one has $v_2>5$ and as $v_2$ is
(strictly) augmented by applying $\bt$ moves to a
positive diagram (by Polyak-Viro formula, see \cite[exercise 4.3%
]{pos}), $7_3$ and $7_5$ do not occur here.

We are left with $9_{25}$, $10_{58}$, $8_{14}$, $8_{12}$,
$7_7$, $7_6$, $6_3$ and $6_2$. We discuss these series separately
in brief phrases.

\begin{itemize}
\item [$6_2$:] Making $6_2$'s diagram positive by crossing 
changes, we obtain $5_1$. The (positive) diagram has Dowker 
notation $4\ -8\ 10\ 12\ -2\ 6$, the alternating one the
same notation only without minus signs. 


By increase of $v_2$ under $\bt$ moves on a positive
diagram we need to apply twists on the positive
generator diagram only as long as $v_2 \le 5$.
 
Twisting at crossings 2 to 6, we obtain the diagrams
$\ol{8_4}$ and $\ol{8_{11}}$ of $7_3$ (the $P(1,-4,3)$ and
$P(1,-2,5)$ pretzel diagrams), and at crossing 1
the diagram $\ol{8_6}$ of $7_5$. In the case of the
diagrams of $7_3$, further twists can be excluded, since
$v_2=5$, but for $7_5$ (with $v_2=4$), we must also
consider a double twist at crossing 1. This gives
a diagram of $9_7$ (with $v_2=5$), which finishes
the case distinction for the series of $6_2$.

\item [$6_3$:] $-4\ -8\ 10\ -2\ 12\ 6$.
 
Since $v_2$ attains the value $5$, $\bt$ moves at crossings $2,3,4$
or $6$ cannot appear with another $\bt$ move. These twists yield
the diagram $\ol{8_{13}}$ of $7_3$. Twists at crossings $1$ and $5$
yield the diagram $\ol{8_8}$ of $7_5$. For two twists we
thus need to consider only these two crossings. Twisting twice
at one of them gives $9_7$, and once at each of both $9_{23}$.
Both $9_7$ and $9_{23}$ have $v_2=5$, and so we 
see that there are no more relevant diagrams.

\item [$7_6$:] $4\ 8\ 12\ 2\ -14\ 6\ -10$. 

To save work, note
that we have $5\sim 7$ in the sense of definition \reference{dffr}.
Thus crossing 7 can be excluded from twisting at. Without twists,
this is a diagram of $5_1$. The twists at crossings $2$, $3$,
$4$ or $6$ give the diagrams $\ol{9_{12}},\ \ol{9_{21}}$ of $7_3$.
The twists at crossings $1$ and $5$ result in the diagrams
$\ol{9_{8}}$ and $\ol{9_{15}}$ of $7_5$. Two twists at crossing
$1$ or $5$ give $9_7$, and one twist at each of both $9_{23}$,
and so we are done.

\item [$7_7$:] $-4\ 8\ -10\ 12\ 2\ -14\ 6$.

Without twists, this is a diagram of $5_1$. The twists at crossings 
$2$, $4$, $5$ or $7$ give the diagram $\ol{9_{14}}$ of $7_3$.
The twist at crossing $3$ gives its diagram $\ol{9_{37}}$.
The twists at crossings $1$ or $6$ give the diagram $\ol{9_{19}}$
of $7_5$. Two twists at latter crossings again give $9_7$ and
$9_{23}$.

\item [$8_{12}$:] $4\ -8\ 14\ 10\ -2\ -16\ 6\ -12$.

This is a diagram of $5_1$. We have three reverse
clasps, $(2,5)$, $(3,7)$ and $(6,8)$, and also $1\sim 4$.
Thus consider only crossings $1,2,3$ and $6$.
Twisting once at $1$ or $6$, we obtain $7_5\ (\ol{10_{35}})$ and
at $2$ or $3$, $7_3 \ (\ol{10_{13}})$ with $v_2=5$. For two twists
we need to consider only crossings $1$ and $6$. Then one
obtains diagrams of $9_7$ and $9_{23}$ with $v_2=5$.
%
%
%
%
Thus more twists cannot give any diagram of interest.

\item[$8_{14}$:] $4\ -8\ 10\ 14\ -2\ 16\ 6\ 12$.

Without a twist this is a diagram of $7_5$ ($\ol{8_{14}}$).
The alternating diagram has a negative clasp $(2,5)$. Not affecting
a crossing there by a $\bt$ move gives a diagram of an alternating
knot of $\ge 9$ crossings, which is excluded. Thus consider
twists at a crossing in the clasp (both crossings
are equivalent with respect to twists). A twist
gives $9_{18}$ with $v_2>5$, which is excluded, so there are no
more diagrams of $7_3$ and $7_5$.

\item[$9_{25}$:] $4\ 8\ 12\ 2\ -16\ 6\ 18\ -10\ 14$.

Again there is a negative clasp $(5,8)$. Use the above argument (for 
$8_{14}$). Without twists it is $8_{15}$, and with one twist
near a crossing in the clasp $9_{18}$ with $v_2>5$, so there are
no diagrams.

\item[$10_{58}$:] $4\ -8\ 14\ 10\ -2\ -18\ 6\ 20\ -12\ 16$.

This is a diagram of $8_{15}$. With one twist we obtain diagrams
of $10_{55}$ and $10_{63}$ with $v_2\ge 5$,
so there are no diagrams we seek.

\end{itemize}

By this exhaustive case distinction we have the desired description. 
\qed

Beside the diagrams we were interested in, we came across
many others used to exclude further possibilities.
This we use to remark also the following useful

\begin{exam}
The knot $!10_{145}$ is not positive.
It is then obviously almost positive as shows its Rolfsen diagram
\cite[appendix]{Rolfsen}. This is the reason for the difficulties
to show its non-positivity by obstructions based on skein arguments 
(see e.g. \cite{MorCro}), as skein arguments apply for almost
positive knots in the same way as for positive ones. The first
non-positivity proof is due to Cromwell \cite[corollary 5.1]{Cromwell} 
using the monicness of the Alexander polynomial. In our context 
the fact follows from the proof of proposition \reference{p7375}. 
We have $v_2(!10_{145})=5$ and $g(!10_{145})=2$, and so if
$!10_{145}$ were positive, we would have encountered it in the
above case distinction, but it did not.
\end{exam}

\section{Classifying all $2$-almost positive diagrams of a slice
or achiral knot\label{s5}\label{S6}}

In this section we give a proof of the classification, announced by
Przytycki and Taniyama in \cite{PrzTan}, of 2-almost positive achiral
and slice knots. Our proof will actually also describe all 2-almost
positive diagrams of such knots, in particular of the unknot (thus
extending the result announced by Przytycki and Taniyama and proved in
\cite{apos}, determining all almost positive unknot diagrams), although
for the unknot this result is not self-contained enough to be nicely
formulable in a closed statement.

\begin{prop}\label{ppaa}
The only non-trivial achiral 2-almost positive knot is $4_1$ (the 
figure eight knot), and the only non-trivial slice 2-almost positive 
knot is $6_1$ (stevedore's knot). Each one of them has only the two
obvious 2-almost positive diagrams.
\end{prop}

The procedure for this task is similar to the one in the previous
section, with the difference that it is better here to use the
signature instead of Vassiliev invariants.

\proof By the slice Bennequin
inequality (see \cite{Rudolph}), 2-almost positive diagrams
of achiral or slice knots have canonical genus $\tl g\le 2$ and $\sg=0$.
For simplicity we content ourselves only to the (interesting)
case, where the diagram is connected, as the composite
case reduces to it and to the almost positive diagram case.

\begin{itemize}
\item [$\tl g=0$:] A connected diagram of canonical genus zero
has one crossing and is hence not $2$-almost positive.
\item [$\tl g=1$:] If we have a subdiagram like
\[
\diag{7mm}{3}{2}{
  \picPSgraphics{0 setlinecap}
  \pictranslate{0.5 1}{
    \picrotate{-90}{
      \lbraid{0 -0.5}{1 1}
      \lbraid{0 0.5}{1 1}
      \rbraid{0 1.5}{1 1}
      \pictranslate{-0.5 0}{
    }
    }
    }
}
\]
then the diagram $D$ reduces to a connected almost positive
diagram and so $D$ belongs to a positive or almost positive
knot. If such a knot is slice or achiral, then it is the unknot.
Let $p$ and $q$ be odd and even positive integers.
All connected almost positive diagrams of the unknot are
unknotted twist knot diagrams \cite{apos} (that is, a twist knot 
diagram with one of the crossings in the clasp changed). 
Hence $D$ is either an unknotted twist knot diagram with
one of the crossings in the twist changed, a pretzel diagram
$P(3,-1,p)$ with of the crossings in the 3-crossing group
changed, or a rational diagram $C(4,-q)$ with two
of the crossings in the 4-crossing group changed.

If we do not have a subdiagram like the one above, then the 
classification of diagrams of canonical genus one
\cite{apos} shows that we have either a $C(-2,q)$ or $P(p,-1,-1)$
diagram, which are the even and odd crossing number diagrams
of the (negative clasp) even crossing number twist knots.
The only achiral twist knot is $4_1$ (a fact, which is almost
trivial to prove using knot polynomials) and the only
slice twist knot is $6_1$ (a fact, which is less trivial to prove,
and it was done by Casson and Gordon \cite{CassonGordon}, see
\cite[p. 215 bottom]{Kauffman3}). After discussing the case $\tl
g=2$ below, in which only the unknot occurs, we will conclude
that each of these knots has only the two $2$-almost positive
diagrams we just found.

\item [$\tl g=2$:] Again we discuss the 24 cases separately.
Consider all diagrams $D_0$ obtained by switching the crossings of the
generators so that exactly two are negative. Then apply $\bt$ moves
at some of the positive crossings of $D_0$. Using the fact that $\sg$
does not decrease when a $\bt$ move is applied to a positive
crossing in any diagram, we can exclude any diagrams obtained by
$\bt$ moves (at positive crossings) from $D$, if $\sg(D)>0$.
(Here $D$ will be obtained by some $\bt$ moves from $D_0$.)

Hereby some symmetries reduce the number of cases to be checked.
%
%
When fixing the crossings to be switched to become negative,
only one choice of crossing(s) in each $\sim$ and $\ssim$ 
equivalence class needs to be considered. The diagrams for the
other choices are obtained (even after $\bt$-twists) by flypes from
the choice made. Also, when applying twists, it needs to be done only at
one choice of crossing(s) in a $\sim$ equivalence class. (In a $\ssim$
equivalence class, \em{a priori} all crossings must be discussed,
if more than one crossing is involved in the twisting, and we would 
like to take care of mutations. However, signature and unknottedness
are invariant under mutations, so that the outcome of our calculation
\em{a posteriori} justifies also symmetry reduction in $\ssim$
equivalence classes.)

For $5_1$, using the signature and symmetry arguments, and
that $\sg(P(-1,-1,1,3,3))>0$, we see that the only diagrams
with $\sg=0$ are $P(-1,-1,1,1,p)$, with $p$ odd and up to
permutation of the entries, and they are all unknotted.
%
Considering the remaining 23 series, a complete distinction of the
cases was done using KnotScape and is given in the 3 tables on the
following pages. By explicit computation of $\sg$ we find
that $\sg(D_0)>0$ except for the choices of negative crossings
given in the tables below (where the aforementioned
symmetries have been discarded). To explain the notation
used in the tables, therein ``\ $\fbox{$6_2$}\ 1\,3$\,'' means: the
diagram obtained from this of $6_2$ (given by its Dowker notation
specified in \S\reference{S2}) by switching crossings so that
all crossings are positive except $1$ and $3$. It turns out that
in all cases of $\sg(D_0)=0$ the diagram $D_0$ is unknotted. Then we
start applying $\bt$ moves at (combinations of) positive crossings
of $D_0$, noticing that either all these moves do not change $\sg$,
or until some $\bt$ move gives a knot diagram with $\sg>0$. 
In latter case we exclude any further $\bt$ moves at that crossing.
It former case it turns out that we always obtain the unknot.
(That arbitrarily many twists at some specific knot diagram give the
unknot can be seen in each situation directly, but it also follows from
checking the first two diagrams in the sequence because of the result 
of \cite{ST}.) The twisting procedure is denoted, exemplarily,
in the following way:
\[
\begin{array}{c} \fbox{$6_2$}\quad 1\ 3\\ \\ \end{array}
\qquad \begin{array}{c@{\ \lra\ }c} 4 & 3_1 \\ 2*5* & 0_1 \end{array}
\]
The notation means: the diagram $\fbox{$6_2$}\ 1\,3$, described
above, with one twist at the crossing numbered by $4$ gives
the trefoil (with $\sg=2$, so we cannot have twists
at crossing $4$), and arbitrarily many twists at crossings
$2$ and $5$ give the unknot. Hereby
a `$,$' (comma) on the left of a term `$x\lra y$'
means `\em{or}', while `\em{and}' is written as a space.
Thus `$4\ 4,\, 1\ 5$' means double twist at crossing $4$ \em{or}
twists at crossings $1$ \em{and} $5$. The $\sim$ and $\ssim$
equivalences for each generator are denoted below it to justify why
certain crossings are not considered for symmetry reasons. \qed

\begin{table}[htb]
\twrule
\[
\begin{array}{r@{\quad\enspace}*4c}
\fbox{$6_2$} & 1\ 3 & 2\,\lra\,0_1 \\
    &      & 4\,\lra\,3_1 & 2*5*\,\lra\,0_1 \\
    &      & 5\,\lra\,0_1 \\
2\sim 5 &      & 6\,\lra\,3_1 \\
3\ssim 4\ssim 6 \\
    & 2\ 3 & 1\,\lra\,0_1 & 4\ 6\,\lra\,3_1 \\
    &      & 4\,\lra\,0_1 \\
    &      & 6\,\lra\,0_1 & 1*4*\,\lra\,0_1 \\
    &      &              & 1*6*\,\lra\,0_1 \\
\\
    & 3\ 4 & 1\,\lra\,3_1 \\
    &      & 2\,\lra\,0_1 & 2*\,\lra\,0_1 \\
    &      & 6\,\lra\,0_1 & 6*\,\lra\,0_1 \\
    &      &              & 2\ 6\,\lra\,3_1 \\
\\
\end{array}\left|
%
\begin{array}{r@{\quad\enspace}*4c}
\fbox{$7_6$} & 1\ 3 & 2\,\lra\,0_1 & 2\ 4\,\lra\,3_1 \\
    &      & 4\,\lra\,0_1 & 2\ 5\,\lra\,0_1 & 2*5*\,\lra\,0_1 \\
    &      & 5\,\lra\,0_1 & 4\ 5\,\lra\,0_1 & 4*5*\,\lra\,0_1 \\
2\ssim 4\\
3\sim 6 & 2\ 3 & 1\,\lra\,0_1 & 1\ 4\,\lra\,3_1 \\
5\sim 7 &      & 4\,\lra\,0_1 & 1\ 5\,\lra\,0_1 & 1*5*\,\lra\,0_1 \\
    &      & 5\,\lra\,0_1 & 4\ 5\,\lra\,0_1 & 4*5*\,\lra\,0_1 \\
\\
    & 2\ 4 & 1\,\lra\,6_2 \\
    &      & 3\,\lra\,0_1 & 3*\,\lra\,0_1 \\
    &      & 5\,\lra\,3_1 & \\
\\
    & 2\ 5 & 1\,\lra\,0_1 \\
    &      & 3\,\lra\,0_1 & 1*3*\,\lra\,0_1 \\
    &      & 4\,\lra\,3_1 & \\
\end{array}
\right.
\]
\twrule
\caption{Proof of proposition \reference{ppaa}: the series of $6_2$
and $7_6$.}
\end{table}

\begin{table}[htb]
\twrule
\[
\begin{array}{r@{\quad\ }*4c}
\fbox{$6_3$} & 1\ 3 & 2\,\lra\,0_1 & 2\ 2\,\lra\,0_1 & \\
    &      & 4\,\lra\,0_1 & 2\ 4\,\lra\,3_1 & 2*5*\,\lra\,0_1\\
    &      & 5\,\lra\,0_1 & 4\ 4\,\lra\,0_1 & 4*\,\lra\,0_1\\
    &      & 6\,\lra\,3_1 & 2\ 5\,\lra\,0_1 & \\
2\ssim 4 \\
3\ssim 6 & 2\ 3 & 1,\,4\,\lra\,0_1 & 1\ 4,\ 4\ 6,\ 5\ 6\,\lra\,3_1 & 4*5*\,\lra\,0_1 \\
    &      & 5,\,6\,\lra\,0_1 & 1\ 5,\ 1\ 6,\ 4\ 5\,\lra\,0_1 & 1*5*\,\lra\,0_1 \\
    &      &              &                               & 1*6*\,\lra\,0_1 \\
\\
    & 2\ 4 & 1\,\lra\,6_2 & \\
    &      & 3\,\lra\,0_1 & 3\ 6\,\lra\,3_1 & 3*\,\lra\,0_1 \\
    &      & 5\,\lra\,3_1 &                 & 6*\,\lra\,0_1 \\
    &      & 6\,\lra\,0_1 \\
\\
    & 2\ 5 & 1\,\lra\,0_1 \\
    &      & 3\,\lra\,0_1 & 3\ 6\,\lra\,3_1 & 1*3*\,\lra\,0_1 \\
    &      & 4\,\lra\,3_1 &                 & 1*6*\,\lra\,0_1 \\
    &      & 6\,\lra\,0_1 & \\
\\
    & 3\ 6 & 1\,\lra\,3_1 \\
    &      & 2\,\lra\,0_1 & 2\ 2\,\lra\,0_1 & 4*\,\lra\,0_1 \\
    &      & 4\,\lra\,0_1 & 2\ 4\,\lra\,3_1 & 2*\,\lra\,0_1 \\
    &      & 5\,\lra\,6_2 & 4\ 4\,\lra\,0_1 & 
\end{array}\ \left|\,\,
%
\begin{array}{r@{\quad\ }*4c}
\fbox{$7_7$} & 1\ 3 & 2\,\lra\,0_1 & 2*\,\lra\,0_1 \\
    &      & 4\,\lra\,3_1 &    &  \\
    &      & 6\,\lra\,7_6 &    \\
2\sim 5 & \\
4\sim 7 & 1\ 4 & 2\,\lra\,0_1 & 2*6*\,\lra\,0_1 \\
    &      & 3\,\lra\,3_1 \\
    &      & 6\,\lra\,0_1 \\
\\
    & 2\ 3 & 1\,\lra\,0_1 \\
    &      & 4\,\lra\,0_1 & 1*4*\,\lra\,0_1 \\
    &      & 6\,\lra\,3_1 \\
\\
    & 2\ 4 & 1\,\lra\,0_1 \\
    &      & 3\,\lra\,0_1 & 1*3*6*\,\lra\,0_1 \\
    &      & 6\,\lra\,0_1 \\
\\
    & 2\ 6 & 1\,\lra\,0_1 \\
    &      & 3\,\lra\,3_1 & 1*4*\,\lra\,0_1 \\
    &      & 4\,\lra\,0_1 \\
\\
    & 3\ 4 & 1\,\lra\,3_1 \\
    &      & 2\,\lra\,0_1 & 2*6*\,\lra\,0_1 \\
    &      & 6\,\lra\,0_1 \\
\\
    & 3\ 6 & 1\,\lra\,7_6 \\
    &      & 2\,\lra\,3_1 & 4*\,\lra\,0_1 \\
    &      & 4\,\lra\,0_1 
\end{array}\right.
\]
\twrule
\caption{Proof of proposition \reference{ppaa}: the series of $6_3$
and $7_7$.}
\end{table}

\begin{table}[htb]
\twrule
\[
\begin{array}{r*4c}
\fbox{$7_5$} & 2\ 3 & 1\,\lra\,0_1 & 1\ 1\,\lra\,0_1 \\
    &      & 4\,\lra\,0_1 & 1\ 4\,\lra\,0_1 & 1*4*\,\lra\,0_1 \\
1\sim 6,\ 2\ssim 5,\ 3\ssim 4\ssim 7 
    &      & 5\,\lra\,3_1 & 4\ 4\,\lra\,0_1 \\
\\
\fbox{$8_{12}$} & 1\ 3 & 7\,\lra\,5_2 & 4\,\lra\,3_1 & 2*6*\,\lra\,0_1 \\
       & 2\ 3 & 7\,\lra\,3_1 & 5\,\lra\,3_1 & 1*6*\,\lra\,0_1 \\
1\sim 4,\ 2\sim 5,\ 3\sim 7,\ 6\sim 8
       & 2\ 6 & 5\,\lra\,5_2 & 8\,\lra\,3_1 & 1*3*\,\lra\,0_1 \\
\\
\fbox{$8_{14}$} & 1\ 4 & 7\,\lra\,5_2 & 3\,\lra\,3_1 & 2*6*\,\lra\,0_1 \\
       & 2\ 4 & 5,\,7\,\lra\,3_1 & 1,3,6\,\lra\,0_1 & 1*3*6*\,\lra\,0_1 \\
2\sim 5,\ 4\ssim 7,\ 6\sim 8
       & 3\ 4 & 1,7\,\lra\,3_1 & 2,\,6\,\lra\,0_1 & 2*6*\,\lra\,0_1 \\
%
\\
\fbox{$8_{15}$} & 2\ 5 & 4\,\lra\,3_1 & 8\,\lra\,3_1 & 1*3*7*\,\lra\,0_1 \\
       & \\
2\ssim 4,\ 3\sim 6,\ 5\ssim 8
\\ \\
\fbox{$9_{23}$} & 2\ 4 & 5\,\lra\,3_1 & 8\,\lra\,3_1 & 1*3*7*\,\lra\,0_1 \\
       & \\
1\sim 6,\ 2\ssim 5,\ 4\ssim 8,\ 7\sim 9
\\ \\
\fbox{$9_{25}$} & 2\ 5 & 4\,\lra\,3_1 & 8\,\lra\,5_2 & 1*3*7*\,\lra\,0_1 \\
       & \\
2\ssim 4,\ 3\sim 6,\ 5\sim 8,\ 7\sim 9 
\\ \\
\fbox{$10_{58}$} & 2\ 6 & 5\,\lra\,5_2 & 9\,\lra\,5_2 & 1*3*8*\,\lra\,0_1 \\
       & \\
1\sim 4,\ 2\sim 5,\ 3\sim 7,\ 6\sim 9,\ 8\sim 10\\
\end{array}
\]
\twrule
\caption{Proof of proposition \reference{ppaa}: the series of $7_5$
and the 8 to 10 crossing generators.}
\end{table}

\end{itemize}

%
%

\begin{rem}
Looking more carefully at our arguments, we see that we only
needed the knot to be slice or achiral to ensure that the diagram has
genus at most two, then we only used that the signature is zero.
We could therefore hope to eliminate completely the condition
of achirality or sliceness by the condition of zero signature.
(This would reprove the result of Przytycki and Taniyama
\cite{PrzTan} that the only 2-almost positive zero signature
knots are twist and additionally show that they have only
the two obvious 2-almost positive diagrams.) For this we
would basically need a version of the ``slice Bennequin inequality''
of \cite{Rudolph} with signature replacing the slice genus.
But the inequality $\sg(D)\ge |w(D)|-n(D)+1$ is \em{not} true for
arbitrary diagrams. Lee Rudolph disappointed my hopes in this
regard, quoting the braid representation of the untwisted
Whitehead double of the trefoil in Bennequin's paper \cite[fig p.\ 121
bottom]{Bennequin}. It a is 7-string braid (so $n(D)=7$) 
consisting of 8 positive bands (so $w(D)=8$), but clearly $\sg=0$.
\end{rem}
 
\section{Almost positive knots\label{S7}}

Almost positive knots, although very intuitively defined,
are rather exotic~-- the simplest example $!10_{145}$ has 10 crossings.
Therefore, non-surprisingly, several properties of such knots
have been proved. For example, they have positive $\sg,v_2$ and
$v_3$ (see \cite{PrzTan} and \cite{apos}), so they are chiral
and non-slice, and are non-alternating
\cite{restr}. Here we add the following property:

\begin{theorem}
There is no almost positive knot of genus one.
\end{theorem}

\proof Assume there is an almost positive knot $K$ of genus one.
By the Bennequin-Vogel inequality (or ``slice Bennequin inequality''
of \cite{Rudolph}) an almost positive diagram $D$ of $K$ has
genus at most $2$. The classification of genus one diagrams relatively
easily excludes the cases where $\tl g(D)=1$ or $D$ is composite.
Thus again we need to consider the 24 series.

To have an almost positive diagram of an almost positive knot
we need to switch (exactly) one crossing in the generator diagram
to the negative, all others to the positive, and possibly apply $\bt$
moves at the positive crossings.

First note, that the negative crossing must have no $\sim$-equivalent
or $\ssim$-equivalent crossing. Otherwise, after possible flypes,
the negative crossing can be canceled by a Reidemeister II move
or a simple-to-see tangle isotopy, giving a positive diagram.
%

Then note, that the $\bt$ move at a positive crossing $p$ in an
almost positive diagram $D$ changes $\nabla$ (the Conway polynomial)
by a multiple of $\nabla_L$, where $L$ is the link resulting by 
smoothing out the crossing $p$ in $D$. Now, by Cromwell
\cite[corollary 2.2, p.\ 539]{Cromwell},
$\nabla_L$ has only non-negative coefficients,
hence such a $\bt$ move never reduces a coefficient in $\nabla$,
in particular not $[\nabla]_{z^4}$. Hence if at some point
$[\nabla]_{z^4}>0$, any further $\bt$ moves cannot produce a
genus one knot.

In many cases $[\nabla]_{z^4}>0$ already after the crossing switch
(without $\bt$ moves) and we can exclude such cases \em{a priori}.

Finally note, that $D$ must have at least 11 crossings, as the
only almost positive knot of at most 10 crossings is $!10_{145}$,
which has genus two.

There arguments exclude after some check all but 7 of the series.
We discuss these cases in more detail.

The argument we apply for these cases basically repeats itself
7 times and consists mainly in drawing and looking more carefully
at the corresponding pictures to see how to eliminate the
negative crossing by Reidemeister moves in most of the cases,
and to check that in the remaining cases $\md\Dl=2$. I list up
the cases, leaving the drawing of the pictures to the reader.

\begin{itemize}

\item[$6_2$:] We have $3\ssim 4\ssim 6$ and $2\sim 5$. The negative
crossing may be chosen to be $1$ or $3$. If it is $1$, the
diagram simplifies to a positive diagram unless at 3 is twisted.
However, if at any of 3, 4 and 6 it is not twisted, then by
flypes this crossing can be made to be 3, hence at all these
three crossings there must be $\bt$ moves. The resulting 12 crossing 
diagram has $\md\Dl=2$. 


\item[$6_3$:] $2\ssim 4$ and $3\ssim 6$. Then modulo
flypes and inversion the negative crossing only needs
to be chosen to be 1 or 3. In case it is crossing $1$,
the resulting diagram can be transformed into a positive one,
unless at both crossings 2 and 4 there are $\bt$ moves applied,
in which case $\md\Dl=2$. In case crossing 3 is changed to the
negative, the transformation into a positive diagram is always
possible.

\item[$7_6$:] $3\sim 6$, $5\sim 7$, $2\ssim 4$. This reduces
to checking the negative crossings to be 1. Then
the diagram can be transformed into a positive one, unless
at both 2 and 4 is twisted, in which case $\md\Dl=2$. 

\item[$7_7$:] $2\sim 5$, $4\sim 7$ and inversion symmetry
leave us with the negative crossing being
1 or 3. Former case simplifies to a positive diagram unless at crossings
3 and 6 is twisted, and so does latter case,
unless at crossings 1 and 6 is twisted. In both situations $\md\Dl=2$.

\item[$8_{14}$:] $2\sim 5$, $4\ssim 7$ and $6\sim 8$ leave us with
crossing 1 or 3. $1$ simplifies unless 3 is twisted,
in which case $\md\Dl=2$; 3 simplifies unless 1 
is twisted, in which case again $\md\Dl=2$. 

\item[$9_{39}$:] $1\sim 4$, $2\sim 6$, $3\sim 8$ leave us with 
crossings 5, 7 and 9 to be negative. When 5 is negative, then already
$\md\Dl=2$. When one of 7 and 9 are negative, the diagram
simplifies unless at the other one there is a $\bt$ move, 
in which case $\md\Dl=2$.
\end{itemize}

Finally, we have

\begin{itemize}
\item[$9_{41}$:] $2\sim 6$, $3\sim 8$, $5\sim 9$ leave 1, 4 and 7
to be negative. However, the diagram has (modulo $S^2$ moves) a
$\bZ_3$-symmetry (rotation around $2\pi/3$), hence we need to deal
just with crossing 1 switched to the negative. This simplifies to
a positive diagram unless at both $4$ and $7$ is twisted, in which case
$\md\Dl=2$. \qed
\end{itemize}

Similar properties to the one I proved remain still open.

\begin{question}
Is there an almost positive knot of 4-ball genus or unknotting
number one?
\end{question}

The expected answer to both is negative. (Note that in this case the
answer to the second part of the question is a consequence of the
answer to the first part.) To give a negative answer,
one could try to apply the argument excluding $10_{145}$~-- namely
that it has an almost positive genus three diagram~-- to
the other knots occurring in our proof whose diagrams
are not straightforwardly transformable into positive ones
(instead of showing $\md\Dl=2$ for them), but this appears to
require hard labor.

\section{Unique and minimal positive diagrams\label{S8}}

One of the achievements of the revolution initiated by the
Jones polynomial was the proof of the fact that an alternating knot
has an alternating diagram of minimal crossing number \cite{Murasugi}.
Unfortunately, such a sharp tool is yet missing to answer the problem
in the positive case. Hence the question whether there is a positive
knot with no positive minimal diagram is unanswered. In \cite{restr}
I managed to give the negative answer to this question in the case
the positive knot is alternating, and subsequently I received a paper
\cite{Nakamura}, where this result was proved independently.
Moreover, it follows from \cite{gen1}
that the answer is the same for (positive) knots of genus one (in
fact, a positive genus one knot is an alternating pretzel knot). Here
we extend this result to genus two.

\begin{theorem}\label{thpos}
Any positive genus two knot has a positive minimal diagram.
\end{theorem}

With this theorem we finish also the proof of corollary
\reference{on9}. The main tool we use to prove it is the $Q$ polynomial
of Brandt-Lickorish-Millett \cite{BLM} and Ho \cite{Ho}
(sometimes also called absolute polynomial) and some results about
its maximal degree obtained by Kidwell \cite{Kidwell}. (They were later
extended by Thistlethwaite to the Kauffman polynomial.)

Recall, that the $Q$ polynomial is a Laurent polynomial in one
variable $z$ for links without orientation, defined by being 1 on the
unknot and the relation
\begin{eqn}\label{Qrel}
A_{1}+A_{-1}=z(A_{0}+A_{\infty})\,,
\end{eqn}
where $A_i$ are the $Q$ polynomials of links $K_i$ and $K_i$ ($i\in\bZ
\cup\{\infty\}$) possess diagrams equal except in one room, where 
an $i$-tangle (in the Conway sense) is inserted, see figure
\ref{figtan}. (Orientation of any of the link components
is unimportant for this polynomial.)

\begin{figure}
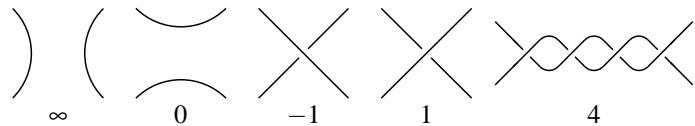

\[
\begin{array}{*6c}
\diag{6mm}{2}{2}{
  \piccirclearc{-1 1}{1.41}{-45 45}
  \piccirclearc{3 1}{1.41}{135 -135}
} & 
\diag{6mm}{2}{2}{
  \piccirclearc{1 -1}{1.41}{45 135}
  \piccirclearc{1 3}{1.41}{-135 -45}
} &
\diag{6mm}{2}{2}{
\picmultiline{-8.5 1 -1 0}{0 0}{2 2}
\picmultiline{-8.5 1 -1 0}{2 0}{0 2}
} &
\diag{6mm}{2}{2}{
\picmultiline{-8.5 1 -1 0}{2 0}{0 2}
\picmultiline{-8.5 1 -1 0}{0 0}{2 2}
}
&
\diag{6mm}{4}{1}{
  \picmultigraphics{4}{1 0}{
    \picline{0.2 0.8}{0.8 0.2}
    \picmultiline{0.22 1 -1.0 0}{0.2 0.2}{0.8 0.8}
  }
  \picmultigraphics{3}{1 0}{
    \piccirclearc{1 0.6}{0.28}{45 135}
    \piccirclearc{1 0.4}{0.28}{-135 -45}
  }
  \picline{-0.2 -0.2}{0.3 0.3}
  \picline{-0.2 1.2}{0.3 0.7}
  \picline{4.2 -0.2}{3.7 0.3}
  \picline{4.2 1.2}{3.7 0.7}
}
\\[2mm]
\infty & 0 & -1 & 1 & 4 
\end{array}
\]
\caption{\label{figtan}The Conway tangles.}
\end{figure}

The following result on $\md Q$ is the one that will 
be subsequently applied.

\begin{theorem}\label{thK}(Kidwell \cite{Kidwell})
Let $K$ be a knot. Then $\md Q(K)\le c(K)-1$ with equality
if and only if $K$ is alternating.
\end{theorem}

\begin{corr}\label{cr1}
Let $D$ be a positive diagram with $\md Q(D)=c(D)-2$. Then the
knot $K$ represented by $D$ has a positive minimal diagram.
\end{corr}

\proof By the theorem \reference{thK}, either $c(K)=c(D)$, in which case the claim
is trivial, or $c(K)=c(D)-1$ and $K$ is alternating, in which case the 
claim follows from the above mentioned result of \cite{restr}. \qed

\begin{lemma}\label{lm1}
With the above notation the $Q$ polynomial
satisfies the following property:
\begin{eqn}\label{Qrel2}
A_n\,=\,(z^2-1)(A_{n-2}-A_{n-4})+A_{n-6}\,.
\end{eqn}
\end{lemma}

\proof We have from \eqref{Qrel}
\begin{eqn}\label{Qrel3}
A_n+A_{n-2}\,=\,z(A_{n-1}+A_{\infty})\,.
\end{eqn}
Now, adding two copies of \eqref{Qrel3} for $n$ and $n-2$ we obtain
\[
A_n+2A_{n-2}+A_{n-4}=2zA_{\infty}+z(A_{n-1}+A_{n-3})=(z^2+2z)
A_{\infty}+z^2A_{n-2}\,.
\]
So
\[
A_n=(z^2+2z)A_{\infty}+(z^2-2)A_{n-2}-A_{n-4}\,.
\]
Therefore
\[
A_n-(z^2-2)A_{n-2}+A_{n-4}=A_{n-2}-(z^2-2)A_{n-4}+A_{n-6}\,,
\]
which is equivalent to the assertion. \qed

\proof[of theorem \reference{thpos}]
Take a positive diagram $D$ of a positive genus 2 knot $K$. If
$D$ is composite, the genus one case shows that $K$ is the
connected sum of two alternating pretzel knots, hence $K$ is
alternating. Thus consider the prime case. The series of $12_{1097}$
and $13_{4233}$ and their progeny in figure \reference{figrel}
contain only positive diagrams
which are alternating, so these cases are trivial. Considering
$11_{148}$, the diagram is made positive by switching the negative
clasp. But all diagrams arising by $\bt$ moves from this diagram can be
simplified near the switched (and possibly twisted) clasp by one
crossing, so as to become alternating. The same argument excludes
(the series of) $10_{97}$, $9_{25}$ and $8_{14}$. The case of $8_{12}$
is trivial, because it contains only rational knots, which are
alternating. For $10_{58}$ and $7_7$ apply the clasp argument
separately to the two negative clasps. For $7_6$ the two negative
clasps cannot be reduced independently, but still drawing a picture
one convinces himself, that a reducing to an alternating diagram by
1 crossing is always possible no matter how many $\bt$ moves
have been performed to the diagram. Here is a typical example:
\[
\diag{1cm}{4}{3}{
  \picputtext[dl]{0 0}{\epsfs[3cm]{k-mv7-6-1}}
}
\quad\lra\quad
\diag{1cm}{4}{3}{
  \picputtext[dl]{0 0}{\epsfs[3cm]{k-mv7-6-2}}
}
\quad\lra\quad
\diag{1cm}{4}{3}{
  \picputtext[dl]{0 0}{\epsfs[3cm]{k-mv7-6-3}}
}\,.
\]


This leaves us with $12_{1202}$, $9_{41}$, $9_{39}$ and $6_3$.
By corollary \reference{cr1} is suffices the check that for any
positive diagram $D$ in their series $\md Q(D)=c(D)-2$.
By the lemma \reference{lm1} and theorem \reference{thK} this reduces
to calculating $Q$ for at most one $\bt$ move applied near a 
crossing and a (reverse) clasp being positive or resolved. However, when the
clasp is resolved, the diagram reduces to one in the series of
some specialization, for which $\md Q(D)=c(D)-2$ or $\md Q(D)=c(D)-3$
by the above discussion. The formula in the lemma \reference{lm1}
then shows that we need to consider just positive clasps without $\bt$ moves.
This leaves a small number of diagrams. E.g., the diagram of $12_{1202}$
consists only of clasps, hence only one diagram needs to be checked.
Switching all crossings in the diagram of $12_{1202}$ to the positive,
we obtain a diagram of the knot $12_{2169}$, for which $\md Q=10$
is directly verified. $9_{39}$ and $9_{41}$ have 3 non-clasp
crossings, hence there are 8 diagrams to be checked, and for $6_3$
we have 64 diagrams. Using various symmetries one can further
reduce the work, but even that far I had no longer serious difficulty
checking the $8+8+64=80$ relevant diagrams by computer. \qed

Our proof actually also shows the following:

\begin{corr}
Any positive (reduced) diagram of a positive genus 2 knot $K$ has at
most $c(K)+3$ crossings. \qed
\end{corr}

This is, in this special case, a much better estimate than the
general bound $c(K)^2/2$ known from \cite{pos}.

The method used in the proof can also be used further. In \cite{gen1}
I exhibited the $(p,q,r)$-pretzel knots with $p,q,r>1$ odd as positive
knots with a unique positive diagram (up to inversion and moves in
$S^2$) and asked whether these are the only examples. The reason behind
this question was that (as I already expected at that point) the number
of series generators grows rapidly with the genus and hence so does
the number of diagram candidates for a positive knot of that genus.
Here we observe that at least for genus 2 the variety on generators
is not sufficiently large, so that such examples still exist.
We take one of our generators.

\begin{exam}
The knot $!10_{120}$ has a unique positive diagram. To see this,
first use that $10_{120}$ is non-arborescent ($\mc Q=6$). This
excludes the series of the knots up to $9_{25}$, and $10_{58}$.
The series of $9_{38}$, $10_{101}$, $11_{123}$, $11_{329}$, $12_{1097}$
and $13_{4233}$ are excluded because all positive diagrams in
these series are alternating (and the only $\bt$ irreducible diagrams
they contain are the generators themselves and $\bt$ (ir)reducibility
is preserved by flypes). $10_{97}$ is excluded, because by the
above discussion the maximal degree of $Q$ on positive
diagrams in its series
is equal to the crossing number minus 2, and hence all maximal degrees
are even (whereas clearly $\md Q(10_{120})=9$). The same argument
excludes $12_{1202}$ and reduces checking positive diagrams in the
series of $11_{148}$ only to the one with no $\bt$ moves applied,
which belongs to $10_{101}$, and the diagrams of $9_{39}$ and $9_{41}$
made positive by crossing changes and with exactly one
$\bt$ move applied. In all the latter cases $\md V=11$, whereas
$\md V(!10_{120})=12$. To finish the argument, it remains only
to notice that the alternating diagram of $10_{120}$ does not admit a
flype itself and because of the crossing number, any $\bt$ twisted
diagram in its series cannot belong to it.
\end{exam}

\section{Some evaluations of the Jones and HOMFLY polynomial\label{S9}}

\subsection{Roots of unity}

The first obstruction to particular values of $\tl g$ is an inequality
of Morton \cite{Morton}: $\max\deg_m\,P/2\,\le\,\tl g$, which shows
that $\tl g>g$ for the untwisted Whitehead double of the trefoil
\cite[remark 2]{Morton} and also for one of the two 11 crossing knots
with trivial Alexander polynomial, which according to \cite[fig. 5]
{Gabai} has genus 2 (I cannot identify which one). In this
section we will discuss an alternative approach to such an obstruction,
and apply it to exhibit the first examples of knots on which the weak
genus inequality of Morton is not sharp. 

\begin{theorem}\label{theoX}
There exist knots $K$ with $\tl g(K)>2=\max\deg_m\,P(K)/2$\,.
\end{theorem}

The present classification opens the question for alternative
criteria which can be applied to exclude a knot from belonging
to a given $\bt$ twist sequence. (We noted that some of the
$\bt$ twist sequences contain others, so we need to consider
only main $\bt$ twist sequences.) Such a criterion is the following
fact, which is a direct consequence of the skein relations
for the Jones \cite{Jones} and HOMFLY \cite{HOMFLY} polynomial 
and has been probably first noted by Przytycki \cite{Przytycki}.

\begin{theo}(Przytycki)
Let $a^{2k}=1$, $a\ne \pm 1$. Then $V(a)\in\bC$ and $P(ia,m)\in\bC[m^2]$
are $\bar t_{2k}$(-move) invariant.
\end{theo}

\begin{corr}\label{theval}
The sets
$\cP_{k,g}\,:=\,\bigl\{\,P_K\bmod\,\mbox{\small$\ds \frac{(il)^{2k}-1}{-l^2-1}$}\,
\big|\ \tl g(K)=g\,\bigr\}$ 
and 
$\cV_{k,g}\,:=\,\bigl\{\,V_K\bmod\,\mbox{\small$\ds \frac{t^{2k}-1}{t^2-1}$}\,
\big|\ \tl g(K)=g\,\bigr\}$ 
are finite for any $k$ and $g\in\bN$.
\end{corr}

\proof From the theorem it is obvious that for every generator $K$,
the sets of residui
\[
\begin{array}{l}
V_{k,K'}\,:=\,V_{K'}\bmod\,\mbox{\small$\ds\frac{t^{2k}-1}{t^2-1}$}
\mbox{\quad and\quad }P_{k,K'}\,:=\,P_{K'}\bmod\,\mbox{\small$\ds
\frac{(il)^{2k}-1}{-l^2-1}$}
\end{array}
\]
are finite on the series of $K$. $\cP_{k,g}$ and $\cV_{k,g}$ are a
finite union of such sets. \qed


\proof[of theorem \reference{theoX}]
We will explain how the knots have been found. The obvious idea is to
compute the sets in corollary \reference{theval} for some appropriate
$k$ in all 24 series and to hope not to find the value of some knot
therein, for which $\max\deg_m\,P\le 4$. Note, that the polynomials are
preserved by mutations, so we need to consider just one diagram for
any generating knot.


The cases $k=2$ and $k=3$ did not suggest themselves as particularly
interesting at least for $V$, because the corresponding evaluations
can be well controlled \cite{LickMil2,Lipson}. Thus I started with $k=4$.
In case of $V$, this is mainly the information given by its evaluations
at $e^{\pi i/4}$ and $e^{3\pi i/4}$ (modulo conjugation and the value
at $i$, which is equivalent to the Arf invariant \cite[\S 14]{Jones2}
and hence not very informative). Table \reference{tabev} summarizes the
number of evaluations of each series.

\begin{table}
{\small
\[
\begin{array}{c*{14}{|c}}
\mbox{\small series}\ry{1.3em}
& 5_1 & 6_2 & 6_3 & 7_5 & 7_6 & 7_7 & 8_{12} &
8_{14} & 8_{15} & 9_{23} & 9_{25} & 9_{38} & 9_{39} & 9_{41} \\[4pt]
\hline\ry{1.3em}
\#\,V_{4,K} &
47 & 121 & 202 & 226 & 136 & 119 & 52 & 302 & 702 & 418 & 479 & 1195 & 413 & 268 \\[4pt]
\hline\ry{1.3em}
\#\,P_{4,K} &
47 & 121 & 202 & 226 & 136 & 119 & 52 & 302 & 710 & 418 & 487 & 1231 & 413 & 268 \\[4pt]
\hline\ry{1.3em}
\#\,V_{5,K}\,=\,\#\,P_{5,K}
& 112 & 408 & 919 & 988 & 538 & 456 & 146 & 1610 & 4281 & 2634 & 2554 & 8588 & 2271 & 1270 \\[4pt]
\end{array}
\]
\[
\begin{array}{c*{11}{|c}}
\ry{1.3em}\mbox{\small series} & 
10_{58} & 10_{97} & 10_{101} & 10_{120} & 11_{123} & 11_{148} &
11_{329} & 12_{1097} & 12_{1202} & 13_{4233} & \mbox{\small total} \\[4pt]
\hline\ry{1.3em}
\#\,V_{4,K} &
157 & 980 & 2380 & 2587 & 2284 & 1041 & 2858 & 5791 & 197 & 5604 & 6645 \\[4pt]
\hline\ry{1.3em}
\#\,P_{4,K} &
163 & 1020 & 2429 & 2673 & 2349 & 1073 & 2970 & 6084 & 209 & 5915 & 6794 \\[4pt]
\hline\ry{1.3em}
\#\,V_{5,K}\,=\,\#\,P_{5,K}
& 624 & 8161 & 23,714 & 27,510 & 22,817 & 8489 & 34,905 & 104,620 & 938 & 102,940 & 128,898 \\[4pt]
\end{array}
\]
}
\caption{\label{tabev} The number of evaluations of $V$ and $P$
in the eighth and tenth roots of unity on each series, and in total.
(The number of evaluations for $V$ and $P$ coincide for
tenth roots of unity.)}
\end{table}

As established in remark \ref{rem1} and \cite{Gabai}, $\tl
g=g=\max\deg\Delta$ for $\le 10$ crossing knots, so we need to look at
more complicated examples. Examining Thistlethwaite's tables, I found
2010 non-alternating 11 to 15 crossing, for which $\max\deg_mP\le 4$.
(Among these 2010 knots only the expected 12 pretzel knots had
$\max\deg_mP\le 2$.) The unity root test for $V$ and $k=4$ does not
exclude any of these 2010 knots from having $\tl g=2$. The test for
$P$ with $k=4$ produced the same disappointing result. (The above
table shows that it does not bring much improvement compared to $V$.)

However, examining $V$ with $k=5$ exhibited four 15 crossing knots of
the type sought. These examples are shown on figure \ref{figg>2}.
One explanation of this outcome
may be that for $k=5$ all four relevant evaluations (at $e^{k\pi i/5}$
for $k=1,2,3,4$) admit very little control. The only known result about
them is Jones' norm bound for $k=1$ in terms of the braid index
and bridge number (cite \cite[propositions 15.3 and 15.6]{Jones2})
and the fact that this evaluation is
finite on closed 3-braids (see \cite[(12.8)]{Jones2}).
Experiments with $P$ and $k=5$, 
however, revealed significantly more time and memory consuming, and
all the values on all of the 24 series reported by my C++ program 
repeated those of $V$, so considering $P$ appears little rewarding.
\qed

\begin{figure}[htb]
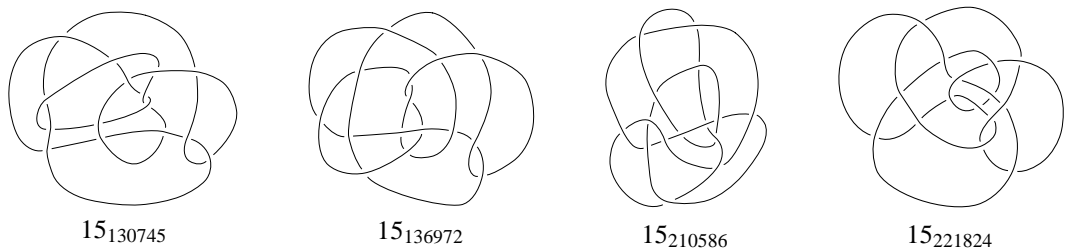

\[
\unitxsize2.8cm
\begin{array}{*4c}
\vis{15}{130745} &
\vis{15}{136972} &
\vis{15}{210586} &
\vis{15}{221824} \\
\end{array}
\]
\caption{\label{figg>2}The simplest examples of knots, for which 
$\tl g>2$ can be proved using Jones polynomial unity root evaluations,
but not using Morton's inequality.}
\end{figure}
 
\begin{rem}\label{remHS}
M.~Hirasawa has found that at least the last example on figure
\reference{figg>2}, $15_{221824}$, has a genus 2 Seifert surface,
so that the ordinary genus is not an applicable obstruction either
to weak genus two in this case.
\end{rem}

It would clearly be helpful to find some nice properties of the sets
occurring in corollary \reference{theval},
but such seem unlikely to exist or at least are obscured by the
electronic way of obtaining them.

Here is some more special example.

\begin{figure}
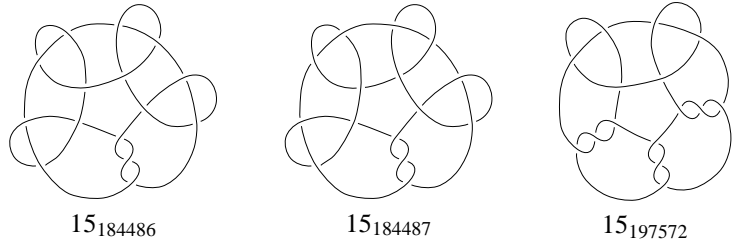

\[
\unitxsize2.8cm
\vis{15}{184486} \quad \vis{15}{184487} \quad
\unitxsize2.55cm
\vis{15}{197572}
\]
\caption{\label{fig15.2}The three 15 crossing pretzel knots
(two of them mutants), 
for which we can show at least that they have no diagram of
even crossing number of genus 2, but which have $\max\deg_mP=4$.}
\end{figure}

\begin{exam}
Consider the knots $15_{184486}$ and $15_{184487}$
on figure \reference{fig15.2}. These knots are slice (generalized)
pretzel knots, which are mutants. Their (common) Jones polynomial is
\begin{eqnarray*}
V(15_{184486}) \, = \, V(15_{184487}) & = &
(-1\ 3\ -6\ 13\ -16\ 17\ -19\ 12\ -7\ [4]\ 4\ -5\ 4\ -3\ 1)
\end{eqnarray*}
A check of the evaluation of 
$V\bmod\mbox{\small$\ds\frac{t^{10}-1}{t^2-1}$}$ shows that the
polynomial modulus is not realized in any main series of
even crossing number. 
Thus these knots do not have a (reduced) genus 2 diagram of even
crossing number (although clearly they have some in the series of
$5_1$). A similar situation occurs for $15_{197572}$.
\end{exam}

\subsection{The Jones polynomial on the unit circle}

While the unity root values of $V$ have been useful in a practical
purpose, we can continue the discussion of the polynomial evaluations 
in a more theoretical direction.

More generally than just in roots of unity, it is possible to say
something about the evaluations of the polynomials on the unit circle.
Here are two slightly weaker but hopefully also useful modifications of 
corollary \ref{theval}. They are also possible for $P$, but I content 
myself to $V$ for simplicity.

\begin{prop}\label{pr}
Let $z\in\bC$ with $|z|=1$ and $z\ne -1$. Then the set $\{\,V_K(z)\,|\,
\tl g(K)=g\,\}\subset\bC$ is bounded for any $g\in\bN$.
\end{prop}

\proof We use the Jones skein relation to expand the Jones polynomial
of a knot in the $\bt$ twist sequence of a diagram in terms of the
Jones polynomials of the diagram and all its crossing-changed
versions. We obtain a complex expression of partial sums of the 
Neumann series for $z^2$ and $z^{-2}$. Now we use the 
boundedness of these partial sums if $|z|=1$ and
$z\ne -1$. (The value $V(1)\equiv 1$ is of little interest.) \qed

\begin{prop}\label{pr2}
Let $z\in\bC$ with $|z|\le 1$ and $z\ne -1$. Then the set $\{\,V_K(z)\,|\,
\mbox{$K$ is positive and } g(K)=g\,\}\subset\bC$ is bounded for any $g\in\bN$.
\end{prop}

\proof In case of positive $\bt$ twists only, the Neumann series 
for $z^{-2}$ do not occur, and we are done as before. \qed

This result seems similar to the boundedness of some other sets
of evaluations of $V$ on closed braids of given strand number
considered by Jones \cite[\S 14]{Jones2}. The nature of our sets is 
quite different, though. Note, for example, that their closure is countable
(so in particular its set of norms has empty interior) for $|z|<1$,
while Jones showed that for the evaluations he considered, the closure
is an interval.

\begin{theo}
The map $f=f_g\,:\,S^1\to \bR$ defined for $g\in\bN_+$ by
\[
f_g(q)\,:=\,\sup\left\{\,\big|V_K(q)\big|\,\big|\,\tl g(K)=g\,\right\}
\]
has the following properties:
{
\def\labelenumi{\theenumi)}
\itemsep0.5pt\relax
\begin{enumerate}
\item $f(\bar q)=f(q)$, where bar denotes complex conjugation;
\item $f(1)=1$, $f(-1)=\infty$;
\item $f$ is upper-semicontinuous on $S^1\sm\{-1\}$, that is, for
$q\in S^1$ and $q\ne -1$ we have $\limsup\limits_{q_n\ne q, q_n\to q}
f_g(q_n)\ge f_g(q)$\,. 
\item $f_g$ satisfies the bound
\[
f_g(q)\,\le\,\max_L \big|V_L(q)\big|\,\cdot\,\left(\frac{2}{|1+q|}+1
\right)^{d_g}\,,
\]
where the maximum is taken over $L$ being a(n alternating) link diagram
obtained by smoothing out some sets of crossings in an alternating $\bt$
irreducible diagram of genus $g$. In particular, the order of the 
singularity of $f_g$ at $-1$ is at most $d_g$.
\end{enumerate}
}
The same properties hold if we modify the definition of $f_g$ by taking
the supremum only over positive or alternating knots.
\end{theo}

\proof The explicit estimate follows from the same argument as in
the proof of proposition 
\ref{pr}. If $V_n$ denote the Jones polynomials of $L_n$, where
$L_n$ are links with diagrams equal except in one room,
where $n$ antiparallel half-twist crossings are inserted,
then from the skein relation for the Jones polynomial we have
\[
V_{2n+1}(q)\,=\,q^{2n}V_1(q)+\frac{q^{2n}-1}{q^2-1}\,\bigl(q^{1/2}-
q^{-1/2}\bigr)V_\infty(q)\,,
\]
with $V_\infty$ denoting the Jones polynomial of $L_\infty$ and $L_\infty$ 
being the link obtained by smoothing out a(ny) crossing in the room.

Expand this relation with respect to any of the $d_k$ crossings,
at which $\bt$ moves can be applied, obtaining $2^{d_k}$ terms to
the right, and take the norm, applying the triangle inequality and
using $|q|=1$.

The upper-semicontinuity of $f_g$ is straightforward from its
definition and the continuity of $V$. Thus the only fact remaining to
prove is 
$f_g(-1)=\infty$.
For this first one easily observes that the determinant (even the
whole Alexander polynomial) depends linearly on the number of $\bt$
twists. Thus we could achieve arbitrarily high and low determinants
in the $\bt$ twist sequence (and at least one of both in alternating
or positive diagrams), unless all linear coefficients in this
dependency are zero. But the fact that the determinant never changes
sign by a $\bt$ twist implies that all knots in the series have
the same signature, and as any diagram can be unknotted by crossing
changes, it must be $0$. But for $\tl g>0$ for example the $(2,2\tl
g+1)$ torus knot has signature $2\tl g>0$, which is a contradiction. 
\qed

\subsection{\label{dense}Jones's denseness result for knots}

{\let\bt\beta

This subsection is unrelated to our discussion as far as weak genus
two knots are considered. However, it is interesting in connection
with (or rather contrast to) the properties of their Jones polynomial 
unity root evaluations.

In \cite[proposition 14.6]{Jones2}, Jones exhibited the denseness
of the norms of $V\left(e^{2\pi i/k}\right)$ on closed $3$-braids
in $[0,4\cos^2\pi/k]$, if $k\in\bN\sm\{1,2,3,4,6,10\}$.

Here we modify this result restricting our attention to \em{knots},
which are closed $3$-braids.

\begin{prop}\label{prd}
If $k\in\bN\sm\{1,2,3,4,6,10\}$, then
\begin{eqn}\label{dV}
\left.\begin{array}{cl} [0,4\cos^2\pi/5-1] & k=5 \\[2mm]
{}[4\cos^2\pi/k-3,4\cos^2\pi/k-1] & k\ge 7
\end{array}\right\}\,\subset\,
\overline{\left\{\,\big|\,V_K\left(e^{2\pi i/k}\right)\,\big|\,:\,
\mbox{ $K$ is a closed $3$-braid knot }\,\right\}} \,\subset\,
[0,4\cos^2\pi/k]\,.
\end{eqn}
\end{prop}

\proof We follow closely Jones's proof. The second inclusion is due
to him. The essential point is the first inclusion.

In the following by $\psi$ we denote the (reduced) Burau
representation. If $\bt$ is a braid, then $\psi_\bt=\psi(\bt)$
is its Burau matrix. We also write $\psi_n$ for the $n$-strand
Burau representation, when dealing with different strand numbers.
(Since numbers and braids are disjoint, the subscripts of $\psi$
cannot be interpreted ambiguously.)

By Jones's proof, we have for $\bt\in B_3$ with even exponent sum
$[\bt]$ (in particular when $\bt$'s closure $\hat\bt$ is a knot), that
\begin{eqn}\label{Vk}
\frac{1}{4\cos^2\pi/k}\,V_{\hat\bt}\left(e^{2\pi i/k}\right)\,=\,
f(\tr(\psi_\bt))\,:=\,
1-\frac{1}{2\cos^2\pi/k} + \frac{1}{4\cos^2\pi/k}\,\tr(\psi_\bt)\,,
\end{eqn}
with $\psi$ being the reduced Burau representation of $B_3$.

Now 
by \cite{Squier}, up to a conjugation (not affecting the trace),
$\psi(\bt)\in U(2)$, and hence, if additionally $k$ divides $[\bt]$,
then 
\[
\psi_\bt\left(e^{2\pi i/k}\right)\Big|_
{\{\,\bt\in B_3\,:\,k|[\bt]\,\}}\subset SU(2)\,,
\]
in particular $\tr(\psi(\bt))$ is real.

Now
\[
\Gm'\,:=\,\{\,\bt\in B_3\,:\,\mbox{ $\hat\bt$ is a knot and $k|[\bt]$ }\}
\]
is a coset in $B_3/\Gm$, where $\Gm$ is the kernel of $\bt\mapsto
\bigl(\sg(\bt),e^{2\pi i[\bt]/k}\bigr)\in S_3\times \bZ_k$. (Here $\sg$ 
is \em{not} the signature, but the induced permutation homomorphism
$B_3\to S_3$.) Again
$\Gm\subset B_3$ is normal and of finite index, hence the closure
of $\psi(\Gm)\subset SU(2)$ has non-trivial connected sets.
In particular the connected component of $1$ contains an $S^1\owns -1$.
Therefore, $\overline{\psi(\Gm')}$ with each $\psi'$ also contains
a coset of $S^1$ we call $G_{\psi'}$ (not necessarily a subgroup),
with $G_{\psi'}\owns -\psi'$.

If now for some $\psi'\in\psi(\Gm')$ we had $\tr(\psi')=\tau$ (where
$\tau\in\bR$), then
$|f|\Big|_{G_{\psi'}}$ would be a continuous function on $G_{\psi'}$
admitting the values $f(-\tau)$ and $f(\tau)$,
and for $\tau\ne 0$ we would apply Jones's argument.

Therefore, we are interested in some $\psi'$ where $|\tau|$ is maximal.
Now if $\xi_{1,2}$ are the eigenvalues of $\psi'$ (with $|\xi_{1,2}|=
1$), then because of $\Gm'^k\,:=\,\{\gm^k\,:\,\gm\in\Gm'\,\}\subset\Gm'$
for any $3\ndiv k$, we consider the maximal trace of $\psi'^k$
with $3\ndiv k$, which is
\[
\mu(\xi):=\sup_{3\ndiv k}\big|1+\xi^k\big|
\]
with $\xi:=\xi_1/\xi_2$. One sees that $\mu$ is minimized by $\xi=e^{\pm
2\pi i/3}$, where it is $1$. Therefore, $f$ ranges at least between 
$f(-1)$ and $f(1)$ on one of the $G_{\psi'^k}$, which implies
the assertion.
\qed


While this is likely not the maximum we can get in our
restricted situation for 3-braids, Jones's corollary specializes
completely to knots.

\begin{corr} \label{crd}
If $k\in\bN\sm\{1,2,3,4,6,10\}$, then
$ \overline{\left\{\,\big|\,V_K\left(e^{2\pi i/k}\right)\,\big|\,:\,
\mbox{ $K$ is a knot }\,\right\}}=[0,\infty)$\,.
\end{corr}

\proof Use that $1$ is always in the interior of the interval to
the left of \eqref{dV} and apply connected sums. \qed

Now we attempt to generalize corollary \reference{crd} to the
case $k=10$. According to Jones \cite[p.\ 263 top]{Jones3},
by the work of Coxeter and Moser \cite{CoxeterMoser}, the image of
$B_3$ in the Hecke algebra is finite, so we need to start with
4-braids, which makes the situation somewhat more subtle.

\begin{prop}\label{pk10}
$\overline{\left\{\,\big|\,V_K\left(e^{\pi i/5}\right)\,\big|\,:\,
\mbox{ $K$ is a knot }\,\right\}}=[0,\infty)$\,.
\end{prop}

\proof First we show that $
\overline{\left\{\,\big|\,V_K\left(e^{\pi i/5}\right)\,\big|\,:\,
\mbox{ $K$ is a $4$-braid knot }\,\right\}}$ contains an interval.
This argument starts along similar lines as the proof of
proposition \reference{prd}.

Consider $\Gm\subset B_4$, which is the kernel of
\[
B_4\owns \bt\,\mapsto\,\bigl(\,[\bt]\bmod 10,\sg(\bt),\psi_3(\ol{\bt})
\bigr) \,\in\,\bZ_{10}\times S_4\times H(e^{\pi i/5},3)\,,
\]
where $H(e^{\pi i/5},3)$ denotes the 3-strand Hecke algebra of parameter
$e^{\pi i/5}$, $\ol{\,\cdot\,}$ is the homomorphism $B_4\to B_3$ with
$\ol{\sg_{1,3}}=\sg_1$, $\ol{\sg_{2}}=\sg_2$, and all other notations
are as before. ($\psi_3=\psi$ is the reduced 3-strand Burau
representation.)
Again $\Gm\subset B_3$ is normal and of finite index, hence the
closure of $\psi_4(\Gm)\subset SU(3)$ is non-discrete.

All subgroups $S^1$ of $SU(3)$ can be conjugated to subgroups of the 
standard maximal toral subgroup, which are of the form
\[
u\in[0,1]\,\mapsto\,
\left(\begin{array}{ccc}
e^{2k\pi iu} & 0 & 0 \\
0 & e^{2l\pi iu} & 0 \\
0 & 0 & e^{-2(k+l)\pi iu}
\end{array}\right)
\]
for some $k,\,l\in\bZ$ with $(k,l)=1$. 
We will refer to these $S^1$s as \em{standard} $S^1$s and denote
them by $S^1_{k,l}$. (The case of $(k,l)>1$ gives no new subgroups,
at least as \em{subsets} of $SU(3)$.)

Therefore, $\ol{\psi_4(\Gm)}$ contains some $AS_{k,l}^1A^{-1}$ for some
$A\in SU(3)$.

Now, consider some $\bt\in B_4$ with $\sg(\bt)$ a 4-cycle,
and write down the weighted trace sum for 4-braids. The result is
\[
V_{\hat\bt}\left(e^{\pi i/5}\right)\,=\,\pi_0(\bt)\,:=\,
8c^3-6c+\frac{1}{2c}\,+\,\frac{1}{c}\tr(\psi_3(\ol{\bt}))\,+\,
\Bigl(6c-\frac{3}{2c}\Bigr)\,\tr(\psi_4(\bt))\,,
\]
with $c:=\cos \mbox{\small$\ds \frac{\pi}{10}$}$. (Keep in mind, that
$\psi_{3,4}$ denote Burau representations of different braid groups.)

If $|\pi_0|\Big|_{\bt\Gm}$ is not constant,
we would find the desired interval. 
Therefore, assume that in particular 
$|\pi_0|\Big|_{\bt \psi_4^{-1}(AS_{k,l}^1A^{-1})}$ is constant. 
Now, on any coset of $\myfrac{B_4}{\Gm}$, $\psi_3(\ol{\bt})$ is
constant, and $AS_{k,l}^1A^{-1}$
acts by multiplying by unit norm complex numbers the columns, so
in particular the 
diagonal entries $\xi_{i}$ of $A\psi_4(\bt)A^{-1}$ ($i=1,2,3$). 
Therefore, for these $\xi_{i}$,
\[
f(u)\,=\,f_{\xi_1,\xi_2,\xi_3}(u)\,=\,
e^{2\pi iku}\xi_1+e^{2\pi ilu}\xi_2+e^{-2\pi i(k+l)u}\xi_3
\]
must lie in some sphere (boundary of some ball) in $\bC$ 
for all $u\in[0,1]$.

That this happens only in exceptional cases follows by
holomorphicity arguments. Namely, if 
\[
\gm=\,-\left(
8c^3-6c+\frac{1}{2c}\,+\,\frac{1}{c}\tr(\psi_3(\ol{\bt}))\right)\Big/
\Bigl(6c-\frac{3}{2c}\Bigr)\,
\]
is the center of this sphere, i.e.
\[
u\,\mapsto\, f_{\xi_1,\xi_2,\xi_3} (u)\,-\,\gm\,
\]
is of constant norm on $[0,1]$, so is 
\[
|f(u)-\gm|^2=(f(u)-\gm)(\ol{f(u)}-\ol{\gm})\,,
\]
which is holomorphic, since $\ol{f(u)}=f_{\ol{\xi_1},\ol{\xi_2},\ol{
\xi_3}}(-u)$. Thus $|f(u)-\gm|^2$ is constant for any $u\in \bC$.

Assume now that $\tr(\psi_4(\bt))=\xi_1+\xi_2+\xi_3\ne 0$.
We claim that $\xi_i\lm_i=0$ for $i=1,2,3$, with $\lm_1:=k$,
$\lm_2:=l$, $\lm_3:=-(k+l)$. In particular, since
(at least) two of the $\lm_i$'s are non-zero, (at least) 
two of the $\xi_i$'s are zero.

Assume the contrary, that is, some $\xi_i\lm_i\ne 0$. Then,
since $(\lm_1,\lm_2,\lm_3)$ is completely characterized
by being a triple of relatively prime integers summing up to 0,
we can by symmetry assume that $\xi_1\ne 0\ne k$. Since any
$\ap\in\bC\sm\{0\}$ is of the form $e^{2\pi iu}$, we have that
\[
P(\ap)\,=\,\ap^k\xi_1\,+\,\ap^l\xi_2\,+\,\ap^{-k-l}\xi_3-\gm
\]
has constant norm for any $\ap\in\bC\sm\{0\}$. Letting $\ap\to 0$
or $\ap\to\infty$, we see that this is possible only if $P(\ap)\equiv
C\in \bC[\ap,\ap^{-1}]$ is a constant as Laurent polynomial in $\ap$.
This in turn is possible (up to interchange of $\lm_{2,3}$ and
$\xi_{2,3}$) only if
(i) $\xi_1=\xi_2=0$ and $k=-l=1$ or (ii) $\xi_1=-\xi_2$, $\xi_3=0$ and
$k=l=1$. Both cases contradict the assumptions $\tr(\psi_4(\bt))\ne 0$
or $\xi_1\ne 0$ resp. 

Thus we have shown that if 
$|\pi_0|\Big|_{\bt \psi_4^{-1}(AS_{k,l}^1A^{-1})}$ is constant,
then $A\psi_4(\bt)A^{-1}\in \cM$, where $\cM$ is the
(closed) subset of $U(3)$, consisting of matrices with zero trace or
least two zero diagonal entries. 

%

But if $\sg(\bt)$ is a 4-cycle, so is $\sg(\bt^{2k+1})$ for any
$k\in\bZ$, so that in particular by the same argument 
any odd power of $A\psi_4(\bt)A^{-1}$
must lie in $\cM$. Taking $\bt=\sg_1\sg_2\sg_3^{-1}$ and setting
$U:=e^{-\pi i/5}A\psi_4(\bt)A^{-1}$ we obtain an element of infinite
order in $SU(3)$, with all its odd powers lying in $\cM$.
But now, $\ol{U^{\bZ}}\subset SU(3)$ is an Abelian closed
non-discrete subgroup, and hence $\ol{U^{\bZ}}$ contains some
$S^1$. But $\ol{U^{\bZ}}$ contains the dense subset $U^{2\bZ+1}$, which
is also a subset of $\cM$, and hence $\ol{U^{\bZ}}$ is contained
itself in $\cM$. Therefore, $\cM\cap SU(3)$ contains an
$S^1=A'S^1_{m,n}A'^{-1}$. 

To show that this is impossible, consider again the trace. If $\tr\ne
0$, we have from the two zero entries and Cauchy-Schwarz for the
third, that $|\tr|\le 1$ on the whole $\cM$. But
integrating the (conjugacy invariant) squared
trace norm on the standard $S^1$, and using that for
any $X\in\bZ[t,t^{-1}]$,
\[
 [X(t)X(1/t)]_{t^0}\,=\,
 \int\limits_0^1\,\Big|\,X\left(e^{2\pi iu}\right)\Big|^2\,du\,,
\]
we obtain
\[
\int\limits_0^1\,\big|\,e^{2\pi imu}+e^{2\pi inu}+e^{-2(m+n)\pi iu}\,
\big|^2 du \,=\,\left\{\begin{array}{cc} 3 & \mbox{for\ }
\big|\{\,m,\,n,\,-m-n\,\}\big|=3
\\[2mm]
5 & \mbox{for\ }\big|\{\,m,\,n,\,-m-n\,\}\big|=2 \end{array} \right.\,.
\]
Thus we must have $|\tr|>1$ somewhere on the standard, and hence
on any other $S^1\subset SU(3)$, providing us with the desired
contradiction.

In summary, we showed that $\big|\,V_K\left(e^{\pi i/5}\right)\,\big|$
is dense in some interval when taking knots $K$ ranging over
closed $4$-braids. From this the proposition follows by taking
connected sums once we can show that there are knots $K_{1,2}$ with
$\big|\,V_{K_1}\left(e^{\pi i/5}\right)\,\big|>1$ and
$0\ne\big|\,V_{K_2}\left(e^{\pi i/5}\right)\,\big|<1$. Luckily, already
$K_1=3_1$ (trefoil) and $K_2=5_1$ ($(2,5)$-torus knot) do the job,
and we are done. \qed

\begin{rem}
V.~Jones pointed out, that for $l=0,\dots,n-1$, 
$\big|\,V\left(e^{l\pi i/n}\right)\,\big|$ is invariant under
a $n$-move (adding or deleting subwords $\sg_i^{\pm n}$).
Thus for $4\nmid k$ our result follows directly from
his, in particular for $k=10$. However, since no proof was given
in this case in \cite{Jones2}, it is worth including one here
anyway. For $4\mid k$, the invariance of $\big|\,V\left(e^{l\pi i/n}
\right)\,\big|$ simplifies the first proof only so far, that it suffices
to consider $\bt\in\Gm$, instead of some non-trivial coset in $B_3/\Gm$.
%
%
This does not alter the proof substantially.
\end{rem}

There is another way to prove the last two statements on norm denseness
in $[0,\infty)$, avoiding any braid group theory, and just applying
connected sums. It would pass via showing for every $k$ the existence 
of knots $K_{1,2}$, such that $\ln \big|\,V_{K_1}\left(e^{2\pi i/k}
\right)\,\big|\big/\ln \big|\,V_{K_2}\left(e^{2\pi i/k}\right)\,\big|$
is irrational. It is unclear how to find such knots for general $k$,
but for single values this is a matter of some calculation. The
following example deals with $k=10$, and thus indicates an alternative
(but much less insightful) proof of proposition \reference{pk10}.

\begin{exam}
Consider the knots $6_3$, $9_{42}$, $11_{391}$ and $15_{134298}$.
Writing $V_{[n]}(t):=\ffrac{t^{n+1/2}+t^{-n-1/2}}{t^{1/2}+t^{-1/2}}$,
note that $V(4_1)=V_{[2]}$ is (up to units) the minimal polynomial of
$e^{\pi i/5}$. The polynomials of our four knots are given by:
\[
V(9_{42})=V_{[3]}\,,\qquad V(6_3)=-V_{[3]}+V_{[2]}+1,\qquad 
V(11_{391})=2-V_{[2]}^2\,,\qquad V(15_{134298})=3-2V_{[2]}^2\,.
\]
Their evaluations at $e^{\pi i/5}$ are $\ffrac{1\pm \sqrt{5}}{2}$,
$2$ and $3$ resp. Then we use that the first two numbers are inverse
up to sign, and $\ln 3/\ln 2$ is irrational. (Except for $6_3$,
the knots are not amphicheiral, although they were chosen to be with
self-conjugate $V$ to make its evaluation at $e^{\pi i/5}$ as simple
as possible.)
\end{exam}

}

To apply our results in this subsection to the weak genus,
we obtain

\begin{corr}
For any even $k>6$ and any $g\in\bN_+$ there are infinitely 
many knots $K$ with braid index
\[
b(K)\es\le\es\left\{\begin{array}{cl} 3 & k\ne 10 \\ 4 & k=10
\end{array}\,\right.\,,
\]
which are not $k$-equivalent to a knot of canonical genus
$\le g$. \qed
\end{corr}

Note that, when replacing $k$-equivalence just by isotopy,
this is well-known, because of the result of Birman and Menasco
\cite[theorem 2]{BirMen} that there exist only finitely many knots
of given (Seifert) genus and given braid index. We will consider
the $k$-moves in more detail later.


\section{$k$-moves and the Brandt-Lickorish-Millett-Ho polynomial%
\label{SecQ}} 

\subsection{The minimal coefficients of $Q$}


It becomes clear from the previous discussion that the Jones polynomial
evaluations for themselves will unlikely give some significantly more
powerful and applicable criteria for showing $\tl g>2$
than Morton's inequality, so it is interesting to find
additional methods that sometimes provide an efficient amplification.
Here we study the $Q$ polynomial in this regard. This is where
the effort in examining the 8th roots of unity of $V$ came to use
in practice.

First, we have the following (not maximally sharp, but easy to apply)
criterion on the low degree coefficients of $Q$.

\begin{prop}\label{PQ}
Let $k$ be a prime. Then $Q\bmod (k,z^k)$ is $\bar t_{4k}$ invariant.
\end{prop}

\proof As in the proof of lemma \reference{lm1},
adding two copies of \eqref{Qrel3} for $n$ and $n-2$ we get
\[
A_n+2A_{n-2}+A_{n-4}=2zA_{\infty}+z(A_{n-1}+A_{n-3})=(z^2+2z)
A_{\infty}+z^2A_{n-2}\,.
\]
Now we iterate this procedure and obtain
\begin{eqn}\label{mq}
\sum_{i=0}^{k}\mybin{k}{i}A_{n-2i}\,=\,z^kA_{n-k}+
z\frac{z^k-2^k}{z-2}A_{\infty}\,.
\end{eqn}
Note that $K_{n-k}$ is a knot when orienting $K_n$ the twists
are antiparallel, even if $k$ is odd (so that $\min\deg A_{n-k}=0$).
Now use the primality of $k$, so that modulo $k$ the left
hand-side collapses to two terms, and we get modulo $k$ and $z^k$
\[
A_n+A_{n-2k}\,=\,\left(z\frac{z^k-2^k}{z-2}\right)\cdot A_{\infty}\,.
\]
Subtracting two copies of this equality for $n$ and $n-2k$
instead of $n$ gives the assertion. \qed

\begin{rem}\label{rqm}
The proof also shows that $Q\bmod (k,z^{k-1})$ is invariant under a
$t_{4k}$ move.
\end{rem}

Working with unity roots of $V$ of order $8$ and $10$ it turns out
useful to consider the criterion for $k=5$.
This criterion has some chance to give partial information as
long as the number of cases left over by the unity root evaluations
is sufficiently less than the total number of values of $Q\bmod
(5,z^5)$, which is very likely $5^5=3125$. 

Another criterion for the Kauffman polynomial $F(a,z+1/z)$ follows
again from Przytycki's work (see
\cite[corollary 1.17, p.~629]{Przytycki}).
%
%
The Kauffman polynomial is a powerful invariant, but, especially
when dealing with many and/or high crossing number diagrams,
too complex for practical computations. Hence, to make this result
more computationally manageable, we set again $a=1$ and use the
$Q$ polynomial. Then from
corollary 1.17 (b) of \cite{Przytycki} it follows that
$Q(z+1/z)$ is invariant under a $\bar t_{2k}$ move for $k$-th
roots of unity $z$. However, we need and prove this condition
in a slightly sharper form, replacing the order $k$ by $2k$.
Our proof is slightly different from (and less technical than)
Przytycki's, since it uses generating series.

\begin{prop}\label{crQ}
$Q(z+1/z)\bmod\mbox{\small$\ds \frac{z^{2k}-1}{z^{3+(-1)^k}-1}$}$ is
$\bar t_{2k}$ invariant, and in particular
$Q(z+1/z)\bmod\mbox{\small$\ds \frac{z^{4k}-1}{z^4-1}$}$ is
$\bar t_{4k}$ invariant. 
\end{prop}

\proof We use the formula in the proof of theorem 3.4 of
\cite{beha}. We observed there that the formula \eqref{Qrel3}
and lemma \reference{lm1} imply that the generating series
\[
f(z,x):=\,\sum_{n=0}^\infty\,A_{2n}(z)x^n
\]
is of the form
\[
f\,=\,\frac{P(z,x)}{(1-x)(1+(2-z^2)x+x^2)}
\]
for some $P\in\bZ[z,x]$. The invariance of $Q(z)$ under a $2k$-move
is equivalent to the denominator dividing $x^k-1$.
Thus we need to choose $z$ so that the zeros of $1+(2-z^2)x+x^2$
are distinct $k$-th roots of unity, different from $1$.
Now if $x_0$ and $x_1$ are these zeros, then $x_0x_1=1$.
Thus $x_{0,1}=e^{\pm 2l\pi i/k}$ for some $0\le l\le k-1$. We must
assume that $l\ne k/2$ (for even $k$) and $l\ne 0$, since then $x_0=x_1$
is a double zero. Then $2-z^2=-x_0-x_1=2\cos (2l\pi/k)$, hence
\[
z^2\,=\,2+2\cos \left(\frac{2\pi}{k}\cdot l\right)\,=\,
4\cos^2\left(\frac{\pi}{k}\cdot l\right),
\]
and $z=\pm 2\cos\left(\mbox{\small$\ds\frac{\pi}{k}$}\cdot l\right)$.
(Since $l$ can be replaced by $k-l$, the sign freedom is fictive.)

Thus $Q(z+1/z)$ is invariant, if $z+1/z=
2\cos\left(\mbox{\small$\ds\frac{\pi}{k}$}\cdot l\right)$,
with $1\le l\le k-1$ and $l\ne k/2$ for even $k$, which means
$z=e^{\pm l\pi i/k}$ for such $l$, and these are exactly the zeros
of the modulo-polynomials stated above. \qed

In the following we decide to use the second property
in proposition \reference{crQ} for $k=5$.
(One could also take $k=10$ for the first property.)

Clearly, the second (Przytycki type) criterion is more powerful,
already because
the number of values of the invariant is infinite. But our first
criterion is easier to compute, and at least it is not a consequence
of the second one, as shows the following


\begin{exam}\label{xf}
Consider $k=5$. The knots $11_{367}$ and $9_{1}$
have $Q$ polynomials that leave the same rest 
modulo $\frac{z^{20}-1}{z^4-1}$. But modulo 5 they differ
in the $z^4$-term, so $11_{367}$ and $9_{1}$ are not
$\bar t_{20}$ equivalent. 
\end{exam}

In this example, the difference of $Q(z)\bmod (k,z^k)$
comes out in the highest coefficient covered (this
of $z^{k-1}$). Surprisingly, this turns out to be the
case for \em{all} other examples I found, that is,
proposition \reference{crQ} implies the weaker version 
of proposition \reference{PQ} for $t_{4k}$ moves noted in
remark \reference{rqm}.

I tested all prime and composite knots of at most 16 crossings
for $k=3,5,7$; for $k=3$ there were about a million coincidences
of $Q(z+1/z)\bmod\mbox{\small$\ds \frac{z^{4k}-1}{z^4-1}$}$
with different $Q(z)\bmod (k,z^k)$, for $k=5$ they were about
3200, and for $k=7$ only one, so in this range of knots
for higher $k$ there are too few coincidences of
$Q(z+1/z)\bmod\mbox{\small$\ds \frac{z^{4k}-1}{z^{4}-1}$}$
to have an interesting picture.

%

\subsection{Excluding weak genus two with the $Q$ polynomial\label{WW}}

The original intention for the $Q$ polynomial
criteria was to exclude further knots
from the set of 2010 from having $\tl g=2$. Then I was fairly surprised
that the most promising candidates (that is, the knots,
whose $V$ moduli appeared the least number of times in the series)
showed up in (at least one of) the series of $12_{1097}$ and
$13_{4233}$. Thus in practice the above criteria have been useful
to reduce the number of diagrams in the series to be considered to
identify these knots. The identification was done using KnotScape.

First, I considered diagrams in the series of $13_{4233}$ and
$12_{1097}$ obtained by switching crossings and performing
at most one $\bt$ move at each crossing/clasp, that is,
with $\le 4$ crossings in each $\sim$-equivalence class.
(Resolving clasps gives diagrams in the subseries of
$13_{4233}$ and $12_{1097}$ in figure \reference{figrel}.)
Then I added all the (other) diagrams in these series
of at most 17, resp.\ 18, crossings. From the set of diagrams
thus obtained, I selected diagram candidates for any knot
with $\md_mP\le 4$ by calculating the Jones polynomial, and
tracking down coincidences. Finally, on the diagrams with
matching polynomials, Thistlethwaite's diagram transformation tool
{\tt knotfind} was applied to identify the knot. By this procedure
I managed to identify all the $\le 15$ crossing knots with
$\md_mP\le 4$ in genus two diagrams expect 6. We already know
four of them~-- they were given in figure \reference{figg>2},
and the other two are shown on figure \ref{figud}.

Thus these 2 knots deserved closer consideration under the
$Q$ polynomial criteria. These criteria proved the two knots to share
the status of those in figure \reference{figg>2}. We give some details
just for the first knot, the other one is examined in the same way.

\begin{figure}[htb]
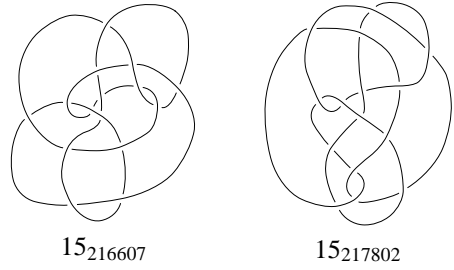

\[
\unitxsize2.8cm
\begin{array}{*2c}
\vis{15}{216607} &
\vis{15}{217802} \\
\end{array}
\]
\caption{\label{figud}
The two prime knots of at most 15 crossings,
for which one can use the $Q$ polynomial to show that
the lower bound $2$ for $\tl g$, coming from Morton's
inequality is not sharp. In all remaining (including composite) cases
it is sharp (if it is 2), except for the four knots on figure
\reference{figg>2}.}
\end{figure}

\begin{exam}
Consider $15_{216607}$ (figure \reference{figud}). We have 
\[
V(15_{216607})\bmod\frac{t^{10}-1}{t^2-1}\,=\,
(\ [-3]\ 0\ -2\ -3\ 1\ -5\ 2\ -4)\,.
\]
It turns out that in the series of $13_{4233}$ the modulus of
$V$ for $k=5$ appears 28 times. They can be encoded by the 
\em{twist vectors}:
\begin{center}
\begin{verbatim}
      {2,-1, 1, 1, 1,-2,-1, 0,-1},         { 1, 0,-1, 1, 1, 1, 1,-2,-2},
      {1, 1, 1, 2,-1, 1,-1, 2,-1},         { 1,-1, 1, 1, 1,-1,-1, 1,-2},
      {1, 1, 1, 0, 1, 1,-2,-2,-1},         { 1,-1, 1,-1, 1, 1,-1, 0,-1},
      {1, 1, 1, 0, 1,-1,-2,-2, 1},         { 0, 1, 0, 1,-2, 1, 0,-2, 0},
      {1, 1, 1, 0, 0, 0,-2,-2, 0},         { 0, 1,-2, 1, 0,-2, 0, 1, 0},
      {1, 1, 1,-1,-1, 1,-1, 2, 2},         { 0, 0, 0, 1, 1, 1, 0,-2,-2},
      {1, 1, 1,-1,-1, 1,-1,-2, 1},         { 0,-2, 0, 1, 1,-2, 0, 1, 0},
      {1, 1,-1, 2,-2, 1, 1,-2,-1},         {-1, 2, 1, 1, 1, 1,-1, 2,-1},
      {1, 1,-1,-1, 1, 1,-2, 1,-1},         {-1,-1, 1, 1, 1, 1, 2, 2,-1},
      {1, 1,-1,-1,-2, 1, 1,-1, 1},         {-1,-1, 1, 1, 1, 1, 1,-2,-1},
      {1, 1,-1,-1,-2, 1, 1,-2, 2},         {-2, 2, 1, 1, 1,-1,-1,-2, 1},
      {1, 1,-2,-1, 2, 1,-1, 0,-1},         {-2, 1, 1, 1, 0, 0, 0,-2, 0},
      {1, 1,-2,-2, 0, 0, 0, 1, 0},         {-2,-1, 1, 1, 1,-1, 2,-2, 1},
      {1, 0, 1, 1, 1, 1,-1,-2,-2},         {-2,-1, 1, 1, 1,-1, 1,-1, 1}
\end{verbatim}
\end{center}

We explain this notation. First,
the crossings are numbered as specified above in the order of the
Dowker notation of $13_{4233}$ given by 
\[
6\ 12\ 22\ \ul{26}\ 16\ \ul{4}\ 20\ 24\ \ul{8}\ \ul{14}\ 2\ 10\ 18\,.
\]
In this notation one skips an entry of a  crossing
appearing in a clasp with some crossing (entry) on its left.
For example, crossings denoted by `$6$' and `$26$' in the
notation form a clasp, so the fourth entry `$26$' is skipped,
and crossing number 4 in the list refers to the crossing represented
by the fifth integer `$16$' in the above Dowker notation.
To facilitate this renumbering, the integers of the crossings to
be skipped are underlined.
An entry $x_i$ at position $i$ ($1\le i\le 9$) in some list
denotes the switching and number of $\bt$ moves applied to the
crossing at number $i$. There are two possibilities.

If the crossing numbered as $i$ is a single element
in its $\sim$-equivalence class, then $x_i=-1$ means
a switched crossing in the alternating diagram, $x_i=0$
the crossing in the alternating diagram as it is, and for
$x_i\ge 1$ (resp. $x_i<-1$) the crossing in the alternating
diagram (resp. the switched one) with $x_i$ (resp. $-1-x_i$)
$\bt$ moves applied to it.

If the crossing $i$ builds (up to flype) a reverse clasp with another
crossing (that is, there are two elements in its $\sim$-equivalence
class), `$x_i>0$' means the clasp as it is with $x_i-1$
$\bt$ moves, `$x_i=0$' means the clasp resolved, and `$x_i<0$'
means the clasp switched with $-1-x_i$ twists applied.

Note, that all the values of $x_i$ need to be considered,
and hence are meant, only modulo $5$.\footnote{To avoid
confusion, let us remark that in a previous(ly cited) version
of the paper a \em{different} convention for the twist vectors
was used. There, for every crossing an entry $x_i=0$ meant
the crossing in the alternating diagram switched, $x_i=1$ the
crossing in the alternating diagram as it is, and $x_i\ge 2$
(resp. $x_i<0$) the crossing in the alternating diagram (resp.
the switched one) with $x_i-1$ (resp. $-x_i$) $\bt$ moves applied to
it. Thus if a crossing builds a (reversely oriented) clasp with
another one, as before `$1$' means the clasp as it is, `$0$' means
the clasp resolved, and `$-1$' the clasp switched.}

Similarly for the other main series (clearly only such
need to be considered) the modulus of $V$ appears 22 times for
the series of $12_{1097}$ and once for $10_{97}$. 

Checking the $51$ diagrams resulting from these vectors modulo 5,
we obtain the following values for $
Q(z+1/z)\bmod\mbox{\small$\ds \frac{z^{10}-1}{z^2-1}$}$:
\[
(\ [-13]\ 0\ 0\ 30\ 28\ 42\ 28\ 30)\,,\qquad
(\ [-23]\ 0\ 0\ 40\ 38\ 62\ 38\ 40)\,,\qquad
(\ [-25]\ 0\ 0\ 54\ 54\ 84\ 54\ 54)\,.
\]
But
\[
Q(15_{216607})\bmod\frac{z^{10}-1}{z^2-1}\,=\,
(\ [-17]\ 0\ 0\ 38\ 40\ 56\ 40\ 38)
\]
does not occur among them. Thus the $Q$ polynomial criterion in
proposition \reference{crQ} excludes all remaining possibilities,
and so $\tl g(15_{216607})>2$.
\end{exam}

\begin{rem}
It is striking that if we take, as above, the rest
$Q(z+1/z)\bmod \ffrac{z^{10}-1}
{z^2-1}$ to be an honest polynomial $P$ in $z$ of degree $\le 7$,
then always $[P]_1=[P]_2=[P]_4-[P]_6=[P]_3-[P]_7=0$ (with
$[P]_i=[P]_{z^i}$). This is in fact
true whatever polynomial $Q\in\bZ[z]$ may be, because the subalgebra of
$\bZ[z,1/z]\big/\ffrac{z^{10}-1}{z^2-1}$ generated by $z+1/z=-z^3-z^5-
z^7$ is the $\bZ$-module with basis $1,z^5,z^4+z^6$ and $z^3+z^7$, and
hence is a rank $4$ subalgebra of an algebra of rank $8$ over $\bZ$.
Therefore, Przytycki's criterion loses on power whenever this
subalgebra (considered also with $10$ replaced by other values $n$)
is small.

For $n$ divisible by $5$ an additional restriction comes from the
Jones-Rong result \cite{Jones4,Rong}, showing that (depending on the
parity of $\dim_{\bZ_5}H_1(D_K,\bZ_5)$) $Q(z+1/z)\bmod\ffrac{z^5-1}
{z-1}$ is always either of the form $\pm 5^k$ or $\pm 5^k(2z^3+2z^2+1)$
for some natural number $k$. 
\end{rem}
 
\begin{rem}
For both knots in figure \reference{figud}
not only the modulus of the Jones polynomial, but
the whole polynomial itself, and even the HOMFLY polynomial,
are realized by weak genus 2 knots ($14_{27627}$, $14_{34335}$ and
$15_{123857}$ for $15_{217802}$ and $14_{35025}$ for $15_{216607}$),
so that the HOMFLY polynomial cannot give complete information on
the weak genus 2 property.
\end{rem}


We obtain in summary that the 6 knots on figures
\reference{figg>2} and \reference{figud}
are indeed the only examples up to 15 crossings, revealing Morton's
inequality, despite these examples, as extremely effective, at least
for $\tl g=2$, even at that ``high'' (in comparison to Rolfsen's
classical tables) crossing numbers.


Beside the ones given above, this quest produced some further interesting
examples with no minimal crossing number diagram of weak genus two.
Contrarily, using a similar argument as in the proof of
theorem \reference{thpos} for the maximal degree of the $Q$ polynomial
on the non-alternating pretzel knots, one can show that for $\tl g=1$
any (weak genus one) knot has a genus one minimal diagram.


\subsection{16 crossing knots}

After the verification of 15 crossing knots, the
16 crossing knot tables were released by Thistlethwaite.
A check therein shows that there are 2249 non-alternating
16 crossing prime knots with $\md_mP\le 4$. (There were
no knots with $\md_mP\le 2$.)
Most of these knots again have weak genus 2. There are 19
knots, which can be excluded using $V\bmod\mbox{\small
$\ds\frac{t^{10}-1}{t^2-1}$}$ and 3 using additionally
$Q(z+1/z)\bmod\mbox{\small$\ds\frac{z^{10}-1}{z^2-1}$}$.
As a counterpart to the knots in figure \reference{fig15.2},
there is one knot, $16_{1265905}$, whose $V$ modulus occurs
only in (main) series of even crossing number, so this knot can
have no reduced weak genus two diagram of odd crossing number.

\begin{figure}
\begin{center}
\captionwidth0.4\textwidth\relax
\@captionmargin0.4\@captionmargin\relax
\begin{tabular}{c@{\qquad}c}
\diag{2.8cm}{1}{1}{
\picputtext[dl]{0 0}{\xepsfs{2.8cm}{t1-16_1265905}}} &
\diag{2.8cm}{1}{1}{
\picputtext[dl]{0 0}{\xepsfs{2.6cm}{t1-16_686716}}}
\\
\ry{6mm}$16_{1265905}$ & $16_{686716}$ \\
\vbox{\caption{\label{fig15.4}A weak genus two knot with 
no reduced weak genus two diagram of odd crossing number.}}
& \vbox{\caption{\label{fig15.5}Does this knot have weak genus two?}}
\end{tabular}
\end{center}
\end{figure}

As another novelty, there is one knot, $16_{686716}$, which
cannot be decided upon. It has the same $V$ and $Q$ moduli
(in fact the same $V$ and $P$, but not $Q$ polynomial) as two
weak genus two knots, $16_{619178}$ and $16_{733071}$. Thus
our criteria cannot exclude weak genus $2$. (Apparently
Przytycki's Kauffman polynomial criteria do not apply
either.) But, after testing all (still potentially relevant)
diagrams in the series of
$12_{1097}$ and $13_{4233}$ of $\le 49$ crossings
corresponding to twist vectors with all $|x_i|\le 9$,
I was unable to find a diagram of this knot. 

\begin{rem}\label{rem10.4}
It is interesting to remark that unusually many of
the above examples are slice, \em{inter alia} $15_{184486}$ and
$15_{184487}$, $15_{221824}$, and the knots on
figures \reference{figud}, \reference{fig15.4} and 
\reference{fig15.5}. It is so far quite unclear whether and what
a relation exists between sliceness and exceptional behaviour
regarding Morton's inequality.
\end{rem}

\subsection{Unknotting numbers and the 3-move conjecture\label{3move}}

Among the family of $k$-moves defined above,
3-moves are of particular interest because of their relation
to unknotting numbers.
An important conjecture of Nakanishi \cite{Nakanishi} is 

\begin{conjecture}(Nakanishi's 3-move conjecture)
Any link is 3-unlinked, that is, 3-equivalent to some (unique) unlink.
\end{conjecture}

This conjecture is by
trivial means true for rational and arborescent links and by non-trivial
work of Coxeter has been made checkable for closures of braids of at
most 5 strands, as he showed in \cite{Coxeter} that $B_n/<\,\sg_i^3\,>$
is finite for $n\le 5$, so proving the conjecture reduces to
verifying (a representative of) a finite number of classes. Qi Chen in
his thesis settled all of them except (the class of) the $5$-braid
$(\sg_1\sg_2\sg_3\sg_4)^{10}$.

As for our context, we get a finite case simplification
for the conjecture for knots of any given weak genus.
The weak genus one case is arborescent and hence trivial, and
we can now do by hand the proof of the 3-move conjecture for
weak genus two knots.

\begin{prop}
Any weak genus two knot is 3-unlinked.
\end{prop}

\proof Applying 3-moves near the $\bt$ twisted crossings in the
24 generators, we can simplify any genus 2 knot diagram to one of the
generators with possibly a crossing eliminated or switched,
and a clasp resolved or reduced to one crossing. We obtain
this way a link diagram of at most 9 crossings. These links
are easy to check directly, but this has previously also
been done by Qi Chen \cite{Chen}. \qed

\subsection{On the 4-move conjecture\label{4move}}

Similar arguments as for the 3-move conjecture allow us to give
a proof of Przytycki's 4-move conjecture for weak genus two knots.

\begin{conjecture} (Przytycki \cite{Przytycki})
Any knot is 4-equivalent to the unknot.
\end{conjecture}

Thus we have 

\begin{prop}\label{cor4}
Any weak genus two knot is 4-equivalent to the unknot.
\end{prop}

\proof  By 4-moves we can simplify any genus 2 knot diagram to one of
the generators of the 24 series with possibly crossings switched. As the
conjecture is verified by Nakanishi for knots of up to 10 crossings,
we need to consider just the diagrams of the 6 last generators
(with possibly crossings switched). In their diagrams
we still have the freedom to change clasps.

The 11 crossing generators and $13_{4233}$ have one of the tangles
\[
T_1\,=\quad
\diagram{1cm}{4}{2}{
  \picPSgraphics{0 setlinecap}
  \pictranslate{4 0}{\picrotate{90}{
  \pictranslate{1 2.8}{
    \picmultigraphics[rt]{4}{90}{
      \picline{0.5 -0.8}{0.5 0}
    }
    \picmultigraphics[rt]{4}{90}{
      \picmultiline{0.12 1 -1.0 0}{0.5 0.8}{0.5 0}
    }
  }
  \piccirclearc{1.0 3.6}{0.5}{0 180}
  \lbraid{1 1.5}{1 1}
  \lbraid{1 0.5}{1 1}
  }}
}
\mbox{\quad and \quad}
T_2\,=\,
\diag{1cm}{2}{1.6}{
  \picline{0.3 1}{1 1}
  \picmulticurve{0.12 1 -1.0 0}{1.0 1.5}{0.7 1.5}{0.4 1.1}{0.7 0.8}
  \picmulticurve{0.12 1 -1.0 0}{1.3 0.8}{1.6 1.1}{1.3 1.5}{1.0 1.5}
  \picmultiline{0.12 1 -1.0 0}{1 1}{1.7 1}
  \picmultiline{0.12 1 -1.0 0}{0.7 0.8}{1.3 0.2}
  \picmultiline{0.12 1 -1.0 0}{0.7 0.2}{1.3 0.8}
  \piccurve{1.0 1.5}{0.7 1.5}{0.4 1.1}{0.7 0.8}
}\,.
\]
It is easily observed that, in which way ever the non-clasp
crossings are changed, the clasps can be adjusted so as the diagram
to simplify by one crossing. Then for the 11 crossing generators we are
done, while for $13_{4233}$ we work inductively over the crossing
number.

$12_{1097}$ has the tangle
\[
\diagram{1cm}{3}{2}{
  \picPSgraphics{0 setlinecap}
  \pictranslate{4 0}{\picrotate{90}{
  \pictranslate{1 2.8}{
    \picmultigraphics[rt]{4}{90}{
      \picline{0.5 -0.8}{0.5 0}
    }
    \picmultigraphics[rt]{4}{90}{
      \picmultiline{0.12 1 -1.0 0}{0.5 0.8}{0.5 0}
    }
  }
  \piccirclearc{1.0 3.6}{0.5}{0 180}
  \lbraid{1 1.5}{1 1}
  }}
  \picputtext{2.5 0.55}{$a$}
  \picputtext{1.5 1.3}{$b$}
  \picputtext{1.5 0.7}{$c$}
}\,,
\]
and the same argument as for $T_1$ applies, unless none (or all) of
crossings $a$, $b$ and $c$ are switched. In this case, by switching the
lower clasp in the diagram of $12_{1097}$, one simplifies the diagram by
2 crossings independently of how the remaining crossings are switched:
\[
\eepsfs[3cm]{k-12-1097-m1}
\quad\lra\quad
\eepsfs[3cm]{k-12-1097-m2}
\quad\lra\quad
\eepsfs[3cm]{k-12-1097-m3}\ .
\]

Finally, the procedure for $12_{1202}$ (and it clasp-switched variants)
is shown below:
\begin{myeqn}{\qed}
\eepsfs[2.2cm]{k-12-1202-m1}
\enspace\lra\enspace
\eepsfs[2.2cm]{k-12-1202-m2}
\enspace\lra\enspace
\eepsfs[2.2cm]{k-12-1202-m3}
\enspace\lra\enspace
\eepsfs[2.2cm]{k-12-1202-m4}
\enspace\lra\enspace
\eepsfs[2.2cm]{k-12-1202-m5}
\end{myeqn}

\section{An asymptotical estimate for the Seifert algorithm\label{S11}}

The Seifert algorithm gives us the possibility to construct a lot of
Seifert surfaces for a knot, and although there is 
not always a minimal one, we may hope
that these cases are rather exceptional. Theorem 3.1 of \cite{gen1}
together with a property of the Alexander polynomial give us the tools
to confirm this in a way we make precise followingly.

\begin{theo}\label{thas}
Fix $g\in\bN_+$. Then
\begin{eqn}\label{conv}
\frac{\#\bigl\{\,D\,:\,\md\Dl(D)=g([D])=g(D)=g,\,c(D)\le n\,\bigr\}}
{\#\bigl\{\,D\,:\,g(D)=g,\,c(D)\le n\,\bigr\}}
\enspace\namedarrow{n\to\infty}\,1\,,
\end{eqn}
where $D$ is a knot diagram, $g(D)$ denotes its genus and $[D]$
the knot it represents. 
\end{theo}

This theorem says that for an arbitrary genus $g$ diagram with many
crossings, the probability the canonical Seifert surface to be of
minimal genus is very high. For the proof we use the Alexander
polynomial.

\begin{rem}\label{remG}
There is a purely topological result due to Gabai, which also can
be applied (see corollary 2.4 of \cite{Gabai2}), as a $\bt$ move
corresponds to change of the Dehn filling of a torus in the knot
complement. However, Gabai needs the condition the manifold
(obtained from the knot complement by cutting out this torus) to
be Haken, and I don't know how to encode this condition in the diagram.
\end{rem}

%

The proof of our theorem bases on the following lemma.

\begin{lemma}\label{lem1}
Let $S$ be s subset of $\bZ^n$ with the following property: if
$(x_1,\dots,x_{k-1},$\ $a,x_{k+1}\dots,x_{n})\in S$ and
$(x_1,\dots,x_{k-1},$\ $b,x_{k+1}\dots,x_{n})\in S$ for some $a\ne b$, then
$(x_1,\dots,x_{k-1},x_k,x_{k+1}\dots,x_{n})\in S$ for all $x_k\in\bZ$.
Then $\fa n\,\ex \eps_n,\,k_n\,\fa k\ge k_n\,:$
\begin{eqn}\label{(l)}
\frac{\is}{(2k+1)^n}\,\ge\,\eps_n\,\Lra\,S\ps [-k,k]^{n}\,.
\end{eqn}
\end{lemma}

\proof Fix some parameter $p\in\bN$ and use induction on $n$. For $n=1$
the claim is evident: set $\ds\eps_1=\frac{1}{2p}$ and $k_1=p$.
Assume now the assertion holds for $n-1$. Let $S\subset\bZ^n$ and set
\[
n_{i,k}\,:=\,\# \bigl(
S\cap ([-k,k]^{n-1}\times \{i\})\,\bigr), \quad |i|\le k\,.
\]
Set $\eps_n:=1-(1-\eps_{n-1})^2$. If now $\x k_0'\,\fa k\ge k_0'$ have 
maximally one $i_0$ such that
\[
\frac{n_{i_0,k}}{(2k+1)^{n-1}}\,\ge\,\eps_{n-1}\,,
\]
then for each such $k$
\[
\sum_{i=-k}^k\,\frac{n_{i,k}}{(2k+1)^{n}}\,<\,\frac{1}{2k+1}+
\eps_{n-1}\,\namedarrow{k\to\infty}\,\eps_{n-1}\,<\,\eps_n\,.
\]
Therefore, $\ex k_0''$ such that $\fa k\ge k_0''$
\[
\frac{\big|\,S\cap [-k,k]^{n}\,\big|}{(2k+1)^{n}}\,<\,\eps_n\,,
\]
and, choosing $k_n$ large enough, there is nothing to prove,
as the premise of \eqref{(l)} does not hold. Therefore, assume
that $\fa k_0':\,\x k\ge k_0'\,\ex i_0\ne i_1:$
\[
\frac{n_{i_0,k}}{(2k+1)^{n-1}}\,\ge\,\eps_{n-1}\,,\qquad
\frac{n_{i_1,k}}{(2k+1)^{n-1}}\,\ge\,\eps_{n-1}\,.
\]
Set $k_n:=k_{n-1}$. Then for $k\ge k_n$ $\x k'\ge k:\,$ $S\ps 
[-k',k']^{n-1}\times \{i_0,i_1\}$. Then $S\ps [-k',k']^{n}$, so
$S\ps [-k,k]^{n}$. \qed

Note, that yet we have the freedom to vary the parameter $p$.
This we need now.

\begin{lemma} 
Lemma \ref{lem1} can be modified by
replacing ``Then $\fa n\,\ex \eps_n,k_n\,:\dots$'' by
``Then $\fa n\,\fa \eps\,\ex k_{n,\eps}\,:\dots$''.
\end{lemma}

\proof Let $p\to\infty$ in the proof of lemma \ref{lem1}. \qed

\proof[of theorem \ref{thas}] Clearly (even taking care of possible
flypes) it suffices to prove the assertion for the $\bt$ twist 
sequence of one fixed diagram $D$, which we parametrize 
using the twist vectors $(n_1,\dots,n_l)$ introduced in \S\reference{WW}
by $\{D(x_1,\dots,x_n)\}_{x_i=-\infty}^{\infty}$, so that a positive
parameter corresponds to a $\bt$ twisted positive crossing. 

Then we apply the previous lemma to
\[
S\,:=\,\{\,(x_1,\dots,x_n)\,:\,\tl g(D(x_1,\dots,x_n))\,>\,
\md \Dl(D(x_1,\dots,x_n))\,\}\,.
\]
The property needed for $S$ is established by the simple fact that
the Alexander polynomials of knots in a 1-parameter $\bt$ twist sequence
form an arithmetic progression.

Denoting $c_{g,n}$ the fraction on the left of \eqref{conv}, assume
$\liminf\limits_{n\to\infty}c_{g,n}<1$. It is equivalent to use 
the $k$-ball around $0$ in the $||\,.\,||_{1}$ or $||\,.\,||_{\infty}$ 
norm, so this means to assume $\ex\eps>0:\,\fa k_0\,\ex k\ge k_0:
\is\,>\,\eps(2k+1)^n$. Then by the second lemma $S\ps [-k,k]^n$ for 
$k\ge k_{n,\eps}$, hence $S=\bZ^n$. But this is clearly impossible as
for example by the canonical Seifert surface minimality of positive 
diagrams $S\cap \bN_+^n=\varnothing$. Hence $\liminf\limits_
{n\to\infty}c_{g,n}=1$. Therefore
$\lim\limits_{n\to\infty}c_{g,n}$ exists, and it is 1. \qed

\section{Estimates and applications
of the hyperbolic volume\label{Shyp}}

We conclude the discussion of the weak genus in general, and
weak genus two in particular, by some remarks concerning the
hyperbolic volume. Surprisingly, it turned out that with
regard to the hyperbolic volume, the setting of \cite{gen1}
has been previously considered in a preprint of Brittenham
\cite{Brittenham2}, of which I learned only with great delay.
Parts of the material in this section (for example, the reference
to \cite{Adams2}) have been completed using Brittenham's work.

\begin{defi}
For an alternating knot $K$ define a link $\tl K$ by adding a circle
with linking number $lk=0$ (i.e. disjoint from the
canonical Seifert surface) around a crossing in each
$\sim$-equivalence class of an alternating diagram of $K$.
\begin{eqn}\label{repl}
\diag{1cm}{2}{2}{
 \piccirclearc{1 1}{0.8}{90 270}
 \picmultivecline{-8 1 -1.0 0}{2 0}{0 2}
 \picmultivecline{-8 1 -1.0 0}{0 0}{2 2}
 \picmulticirclearc{-8 1 -1.0 0}{1 1}{0.8}{270 90}
}\rx{2cm}
\diag{1cm}{2}{2}{
 \piccirclearc{1 1}{0.8}{90 270}
 \picmultivecline{-8 1 -1.0 0}{0 0}{2 2}
 \picmultivecline{-8 1 -1.0 0}{2 0}{0 2}
 \picmulticirclearc{-8 1 -1.0 0}{1 1}{0.8}{270 90}
}
\end{eqn}%
(The orientation of the circles is not important.)
\end{defi}

In this language one can obtain all weak genus $g$ knots
by $1/n_i$-Dehn surgery along the unknotted components of
$\tl K$ for the genus $g$ generators $K$. (In fact, the
main generators suffice, and the cases of composite generators
can be discarded.)

In this situation we can apply a result of
Thurston (see \cite{NZ}). 
To state it, here and below $\vol(K)$ denotes the hyperbolic volume
of (the complement of) $K$, or $0$ if $K$ is not hyperbolic. 
$K(n_1,\dots,n_l)$ denotes, as in \cite{Bseq}, the knot
in the series of $K$ with twist vector $(n_1,\dots,n_l)$,
as explained in \S\reference{WW}.

\begin{theorem}(Thurston)
If $\vol(\tl K)>0$, then for all vectors $(n_1,\dots,n_l)\in\bZ^l$,
\[
\vol\bigl(\,K(n_1,\dots,n_l)\,\bigr)\,<\,
\vol(\tl K)\,,
\]
and
\[
\vol\bigl(\,K(n_1,\dots,n_l)\,\bigr)\,\lra\,
\vol(\tl K)\,,
\]
as $\ds\min_{i=1}^l\,n_i\to\infty$.
\end{theorem}

As a consequence, we obtain the following theorem.

\begin{theorem}\label{thhy}
Let
\[
S_g\,:=\,\{\,\vol(\tl K)\,:\,\mbox{$K$ main generator of genus
$g$}\,\}\,.
\]
Then
\begin{myeqn}{\qed}
\sup\,\{\,\vol(K)\,:\,\tl g(K)\,=\,g\,\}\quad=\quad
\max S_g\,.
\end{myeqn}
\end{theorem}

\proof The $\tl K$ are augmented alternating links in the sense of
Adams \cite{Adams2}, and hence by his result are hyperbolic,
if $K$ is a prime alternating knot different from a torus
knot. Applying Thurston's result, it remains to prove that
the alternating torus knot is never a main generator. This is an
easy exercise. \qed

This theorem shows in particular that the hyperbolic volume of
knots of bounded weak genus is bounded, with an explicitly
computable exact upper estimate. 
%
%
%

In particular, we obtain from theorem \reference{thhy}
by explicit calculation:

\begin{corr}
\[
\sup\,\{\,\vol(K)\,:\,\tl g(K)\,=\,1\,\}\quad=\quad\vol(\tl{3}_1)\approx
14.6554495068355\,,
\]
and
\begin{myeqn}{\qed}
\sup\,\{\,\vol(K)\,:\,\tl g(K)\,=\,2\,\}\quad=\quad\vol(\wt{13}_{4233})
\approx 58.6217980273420\,.
\end{myeqn}
\end{corr}

The (approximate) volumes of $\tl K$ for the main generating knots 
$K$ of genus 2 are given as follows:
\[
\begin{array}{l|c}
\rx{1.0em}K & \vol(\tl K) \\[1mm]
\hline\ry{1.4em}
      6_3\ &\ 36.6386237671 \\
   9_{41}\ &\ 38.7476335870 \\
  10_{97}\ &\ 43.9663485205 \\
 11_{148}\ &\ 43.9663485205 \\
12_{1097}\ &\ 58.6217980273 \\
12_{1202}\ &\ 38.7476335870 \\
13_{4233}\ &\ 58.6217980273
\end{array}
\]

There is a further application of the hyperbolic volume.

\begin{prop}\label{thisp}
If $\vol(\tl K)>\vol(\tl K')$ for two generators
$K$ and $K'$, then a generic alternating knot in the series of
$K$ has no diagram in the series of $K'$. \qed
\end{prop}

%
%
%

%
%
To make precise what `generic' means we make a definition:

\begin{defi}
A subclass $\cB\subset\cC$ in a class $\cC$ of links is called
\em{asymptotically dense} or \em{generic}, if
\[
\lim_{n\to\infty}\,\frac{\big|\,\{\,K\in\cB\,:\,c(K)=n\,\}\,\big|}
{\big|\,\{\,K\in\cC\,:\,c(K)=n\,\}\,\big|}\,=\,1\,.
\]
\end{defi}

For example, in \cite{Thistle2} Thistlethwaite showed that the
non-alternating links are generic in the class of all links.
Similarly, a result of \cite{nlpol} is that any generic subclass
of the class of alternating links contains mutants.

The proof of  proposition \reference{thisp} is similar to the
arguments in \S\reference{S11}, but simpler, and is hence omitted.
(Again avoiding $K'$ to be a torus knot is easy.)

\begin{exam}
We have 
\[
\vol(\tl 9_{38})\approx 47.2069898171\,>\,
\vol(\wt{10}_{97})\approx 43.9663485205\,,
\]
so that a generic alternating knot in the series of $9_{38}$
will not have a diagram in the series of $10_{97}$. (Note that
both series have seven $\sim$-equivalence classes and thus
the number of diagrams in them grows comparably.)
\end{exam}

The fact that $\vol(\wt{13}_{4233})$ and $\vol(\wt{12}_{1097})$
are equal is unfortunate, as otherwise we would be able to 
conclude that a generic genus two alternating knot of one of the
crossing number parities has no genus two diagrams of the other
crossing number parity (as we did for specific examples before
using the values of the Jones polynomial at roots of unity).
Also, this value is much higher than the volume of any non-alternating
$\le 16$ crossing knot. (The maximal volume of such a knot is about
$32.9$, and the maximal volume among those knots with $\md_mP\le 4$ is
about $22.9$.) Thus the volume does not seem to have much
practical significance as an obstruction to $\tl g=2$.
%
On the other hand, we can use the fact that
$\vol(\wt{13}_{4233})=\vol(\wt{12}_{1097})$ is higher 
than $\vol(\tl K)$ for the other main generators $K$. From
this, and proposition \reference{thisp}, we obtain

\begin{corr}
A generic alternating genus two knot has no non-special
genus two diagrams (i.e. such diagrams with a separating Seifert
circle.) \qed
\end{corr}

This is not true for weak genus one, because of the
alternating knots of even crossing number. For odd crossing
number genus one alternating knots it is, contrarily, trivial.
However, being such a narrow class, genus one diagrams are not
interesting anyway.


To estimate $\max S_g$, Brittenham uses a remark of Thurston that
for any link $L$, $\vol(L)\le 4V_0c(L)$, with $V_0$ being the
volume of the ideal tetrahedron. Then he studies
\[
C_g\,:=\,\{\,c(\tl K')\,:\,\mbox{$K$ main generator of genus
$g$}\,\}\,,
\]
where $\tl K'$ is obtained from $\tl K$ by resolving in $K$ clasps
of $\sim$-equivalence classes with two crossings. (This move preserves
the link complement.) Brittenham shows that $\max\,C_g\le 30g-3$. 

We conclude this section by giving an estimate for $\max\,C_g$,
which is the best possible for $g\ge 6$.

\begin{prop}\label{pve}
$\max\,C_g\le 30g-15$, and this inequality is sharp for $g\ge 6$.
\end{prop}

In particular, we have a slight improvement of Brittenham's volume
estimate:

\begin{corr}
\begin{myeqn}{\qed}
\sup\,\{\,\vol(K)\,:\,\tl g(K)\,=\,g\,\}\quad\le\quad
(120g-60)\,V_0\,.
\end{myeqn}
\end{corr}

However, we also know now that a significant further improvement of
Brittenham's volume estimate is possible only by studying the volume
of the $\tl K$ directly, and not via their crossing number.

\proof[of proposition \reference{pve}]
We know from \cite{STV} that $d_g\le 6g-3$, and in each
$\sim$-equivalence class we need 4 crossings for the trivial
loop, and at most one crossing for the generating knot. (Recall that
$d_g$ are the numbers introduced at the end of \S\reference{S2}.)
If it some $\sim$-equivalence class of the generating diagram 
has two $\sim$-equivalent crossings, their clasp can be resolved,
since this preserves the link complement. Thus each $\sim$-equivalence
class contributes at most 5 crossings to $c(\tl K')$, showing the
estimate claimed.

To show that the estimate is sharp, we need to construct a prime
alternating knot $K=K_g$ of genus $g\ge 6$ with $6g-3$
$\sim$-equivalence classes, all consisting of a single crossing.

Once this is done, it is easy to show that $c(\tl K')=c(\tl K)=30g-15$.
Let $L_1,\dots,L_n$ be the trivial components of $\tl K$.
Then $K\sqcup L_i$ is non-split for any $i$, since $1/n_i$
surgery on $L_i$ changes $K$, as it may give an alternating knot
of higher crossing number. Also, as this knot is prime
(by \cite{Menasco} and the primality of the diagram),
$L_i$ cannot be enclosed in a sphere intersecting $K$ in
an unknotted arc (otherwise the result from $K$ after
$1/n_i$ surgery on $L_i$ will always have $K$ as prime factor).
Thus $L_i$ and $K$ have at least 4 mixed crossings in any diagram
of $\tl K$. Since $K$ appears in a reduced alternating diagram
in the diagram of $\tl K$ obtained by the replacements \eqref{repl},
it is also of minimal crossing number.

We give the $K_g$ in terms of their Seifert graphs; since
all $K_g$ are special alternating, these graphs determine
uniquely a special alternating diagram of $K_g$ (see e.g.
\cite{Cromwell}; these graphs are trivalent and bipartite). We include
the graphs only for $g=6$ and $g=7$. (The genus can be determined
easily, since the number of regions of the graph is $2g+1$.)
Given a graph of a knot of genus $g$, one can obtain a graph of
a knot of genus $g+2$ by the replacement
\[
\diag{1cm}{1}{1}{
    \picline{0 0}{1 0}
    \picline{0 1}{1 1}
}\quad\lra\quad
\diag{1cm}{1.8}{1}{
    \picline{0 0}{1.8 0}
    \picline{0 1}{1.8 1}
    \picmultigraphics{4}{0.4 0}{
      \edge{0.3 0}{0.3 1}
    }
}\,,
\]
performed so that the number of edges in each face remains even.
\begin{myeqn}{\qed}
\diag{1cm}{6}{3}{
  \pictranslate{3 1.5}{
    \picellipse{0 0}{3 1.5}{}
    \edge{60 3 x polar 2 :}{60 3 x polar -2 :}
    \edge{90 3 x polar 2 :}{90 3 x polar -2 :}
    \edge{120 3 x polar 2 :}{120 3 x polar -2 :}
    \edge{1.1 0 x}{1.1 1.5 x}
    \edge{0.8 0 x}{0.8 1.5 x}
    \edge{-0.8 0 x}{-0.8 1.5 x}
    \edge{-1.1 0 x}{-1.1 1.5 x}
    \edge{0 0.5}{-1.5 0.5}
    \edge{0 -0.5}{-1.5 -0.5}
    \edge{-0.9 -0.5}{-0.9 -1.41}
    \edge{-1.2 -0.5}{-1.2 -1.37}
  }
}\rx{1.5cm}
\diag{1cm}{6}{3}{
  \pictranslate{3 1.5}{
    \picellipse{0 0}{3 1.5}{}
    \edge{60 3 x polar 2 :}{60 3 x polar -2 :}
    \edge{90 3 x polar 2 :}{90 3 x polar -2 :}
    \edge{120 3 x polar 2 :}{120 3 x polar -2 :}
    \edge{1.1 0 x}{1.1 1.5 x}
    \edge{0.8 0 x}{0.8 1.5 x}
    \edge{-0.8 0 x}{-0.8 1.5 x}
    \edge{-1.1 0 x}{-1.1 1.5 x}
    \edge{0 0.5}{-1.5 0.5}
    \edge{0 -0.5}{-1.5 -0.5}
    \edge{-0.3 0.5}{-0.3 -0.5}
    \edge{-0.6 0.5}{-0.6 -0.5}
    \edge{-0.9 -0.5}{-0.9 -1.41}
    \edge{-1.2 -0.5}{-1.2 -1.37}
  }
}
\end{myeqn}

%

\begin{rem}\label{rBr}
Brittenham uses his proof that weak genus bounds the volume to show in
\cite{Brittenham} that
there are (hyperbolic) knots of genus one and arbitrarily large weak
genus. Because of his use of Thurston's theorem, however, he cannot give
any particular examples. Such examples, although not hyperbolic,
have been previously given in \cite{Moriah,gwg}.
\end{rem}

\begin{rem}
One can check that for $g\le 5$ the knots $K_g$ do not exist.
This follows for $g=1$ from \cite{gen1}, for $g=2$
from our discussion, and for $g=3$ from the calculation given later
in \S\reference{gen3}. For $g=4,\,5$, one can establish this
in the following way. It follows from the results of \cite{MS} and
\cite{SV} that the Seifert graphs of the alternating diagrams
of $K_g$ are exactly the planar, 3-connected, bipartite,
3-valent graphs with $4g-2$ vertices and an odd number of
spanning trees. A list of candidates for such graphs was generated
and then examined with MATHEMATICA. It showed that for $g\le 5$
no such graphs exist.
\end{rem}

\section{Genus three\label{S12}}

\subsection{The homogeneity of $10_{151}$, $10_{158}$ and $10_{160}$}

After having some success with $\tl g=2$, I was encouraged to face
the combinatorial explosion and to try to obtain at least some
partial results about $\tl g=3$. One motivation for this attempt were
the 3 undecided genus 3 knots in \cite[appendix]{Cromwell}. They
can now be settled, and thus, together with corollary
\reference{cor4.2}, Cromwell's table completed.

\begin{prop}
The knots $10_{151}$, $10_{158}$ and $10_{160}$ are non-homogeneous.
\end{prop}

\proof These knots all have monic Alexander polynomial, and hence a
homogeneous diagram must be a genus 3 diagram of at most 12 crossings
\cite[corollary 5.1]{Cromwell} with no $\bt$ move applied (see proof
of \cite[theorem 4]{Cromwell}). As crossing changes commute with
flypes, deciding about homogeneity reduces to looking for homogeneous
diagrams obtained by flypes and crossing changes from a $\bt$-%
irreducible alternating diagram of between 10 and 12 crossings. We can
exclude special alternating series generators, as homogeneous diagrams
therein are alternating (and positive). Since the leading coefficient
of $\Dl$ is multiplicative under Murasugi sum, and invariant up to
sign under mirroring, the monicness of the Alexander polynomial is 
preserved under passing from the homogeneous to the alternating
diagram. Therefore, it suffices to consider only (alternating)
generating knots, whose Alexander polynomial is itself monic.
There are 37 such knots.

Unfortunately, (non-)homogeneity of a diagram, unlike alternation and
positivity, is a condition, which is not necessarily preserved by
flypes. Thus we must apply flypes on the 37 generators, obtaining 275
(alternating) generating diagrams. 

\begin{figure}[htb]
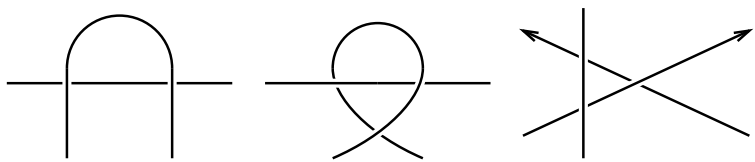

\[
\begin{array}{*3c}
\diag{1cm}{3}{2}{
   \picPSgraphics{0 setlinejoin}
   \picline{0 1}{3 1}
   \picmultiline{0.12 1 -1.0 0}{0.8 0}{0.8 1.2}
   \picmultiline{0.12 1 -1.0 0}{2.2 0}{2.2 1.2}
   \piccirclearc{1.5 1.2}{0.7}{0 180}
}
&
\diag{1cm}{3}{2}{
   \picPSgraphics{0 setlinejoin}
   \pictranslate{1.5 1}{
       \picmultigraphics[S]{2}{-1 1}{
           \picmulticurve{0.12 1 -1.0 0}{0.6 -1}{-0.1 -0.7}
                         {-0.6 -0.2}{-0.6 0.2}
           \picmultiline{0.12 1 -1.0 0}{0 0}{1.5 0}
       }
   }
   \piccirclearc{1.5 1.2}{0.6}{0 180}
}
& 
\diag{1cm}{3}{2}{
   \pictranslate{1.5 1}{
       \picmultigraphics[S]{2}{-1 1}{
           \picmultivecline{0.12 1 -1.0 0}{1.5 -0.7}{-1.5 0.7}
       }
   }
   \picmultiline{0.12 1 -1.0 0}{0.8 0}{0.8 2}
}
\end{array}
\]
\caption{\label{fignh}
Fragments to exclude, together with their obverses, in a homogeneity
test. Unoriented lines may have both orientations. The first and
third fragments above make the diagram non-homogeneous even after
flypes. The second one may or may not do so (depending on the
orientation) but if the diagram is homogeneous, then this property
is not spoiled by reducing the fragment to a clasp (so there is a
simpler homogeneous diagram).}
\end{figure}

We must now consider the diagrams obtained from these 275 by
crossing changes, and then test homogeneity. However, it is useful
to make a pre-selection. There are several simple fragments in a
diagram, which either render it non-homogeneous,
or reveal a simpler homogeneous diagram. See figure \ref{fignh}.
Thus it suffices to consider diagrams without such fragments.
More generally than (excluding) the first fragment, if $p\sim q$ or
$p\ssim q$, then $p$ and $q$ must have equal sign. We also apply
theorem \reference{thK} to discard diagrams with long bridges
(we have $\max\deg Q=8$ for all three knots). There remain to be
considered 1430 diagrams.

Homogeneity test on these 1430 diagrams gives 430 homogeneous ones.
It is easy to check that none of them matches the Alexander
polynomial of any of the 3 knots we seek, and so we are done. \qed

\subsection{The complete classification\label{gen3}}

The above three knots were a motivation to find the $\bt$
irreducible alternating genus 3 knots at least up to 12 crossings.
However, a complete classification of the $\bt$ irreducible
genus three alternating knots is considerably more difficult.

Theorem 3.1 of \cite{gen1} shows that at least at $c_3\le 8c_2+6=110$
crossings the series will terminate. The situation becomes then
more optimistic, though. If one repeats the discussion at the
end of \S\ref{S2} for a $\bt$ irreducible alternating genus 3
diagram, this leads to expect $c_3$ to be around $23$.
Then we found in \cite{STV}, that it is indeed equal to $23$. The
method there used the list of maximal Wicks forms compiled as described
in \cite{BV}. This method becomes increasingly efficient when the 
crossing number grows beyond 15. After some optimization, I was able to
process with it also the crossing numbers below 23, finally reaching
17 crossings. For fewer crossings, one can select generators
directly from the alternating knot tables. (I also processed 16
crossings by both methods to check that the results are consistent.)

The number of generating knots is shown in table \reference{tab1}.
In particular $d_3=15$. These data show that there is a huge number of
generators, which render
discussions by hand, or with moderately reasonable electronic
calculation, as for $\tl g=2$, practically impossible in most cases.

\begin{table}[ptb]
\captionwidth0.8\vsize\relax
\newpage
\vbox to \textheight{\vfil
\rottab{%
\def\hh{\\[0.5mm]\hline[1.2mm]}%
\hbox to \vsize{\footnotesize\hss
\begin{mytab}{|c||r|r|r|r|r|r|r|r|r|r|r|r|r|r|r|r|r||c|}%
  { & \multicolumn{17}{|r||}{ } & }%
  \hline [2mm]%
  \rx{-0.8em}%
  \diag{1cm}{1.5}{1}{
    \piclinewidth{90}
    \picline{0 1}{1.5 0}
    \picputtext{1.2 0.7}{$c$}
    \picputtext{0.5 0.3}{$\#\sim$}
  }\rx{-0.8em}%
 & $7$ &   $8$ &    $9$ &    $10$ &   $11$ &   $12$ &    $13$ &    $14$ &    $15$ &     $16$ &  $17$ &     $18$ &      $19$ &      $20$ &      $21$ &      $22$ &       $23$  & total \\[2mm]%
\hline
\hline [2mm]%
6 &   &   &   &   &   & 4 &   &   &   &   &   &   &   &   &   &   &   & 4 \hh
7 & 1 & 2 & 5 & 8 & 11&   &   & 9 &   &   &   &   &   &   &   &   &   & 36\hh
8 &   & 6 & 10& 21& 22& 30& 44&   &   & 13&   &   &   &   &   &   &   &146\hh
9 &   &   & 4 & 16& 42& 72& 64& 55& 68&   &   &  7&   &   &   &   &   &328\hh
10&   &   &   & 2 & 15& 51&104&159&119& 52& 45&   &   & 2 &   &   &   &549\hh
11&   &   &   &   & 1 & 10& 49&120&194&211&130& 20& 14&   &   &   &   &749\hh
12&   &   &   &   &   & 1 & 5 & 32&112&220&229&154& 75& 2 & 1 &   &   &831\hh
13&   &   &   &   &   &   & 1 & 2 & 17& 63&170&252&178& 48& 18&   &   &749\hh
14&   &   &   &   &   &   &   &   & 1 & 4 & 22& 63&132&163& 82&   &   &467\hh
15&   &   &   &   &   &   &   &   &   &   & 2 & 3 & 12& 25& 47& 46& 23&158\\
[2mm]
\hline
\hline [2mm]%
total& 1 & 8 & 19 & 47 & 91 & 168 & 267 & 377 & 511 & 563 & 598 & 499 & 411 & 240 & 148 & 46 & 23 & $4017$ \\[1.2mm]
\hline
\end{mytab}%
\hss}%
}{The number of $\bt$ irreducible prime genus 3 alternating knots
tabulated by crossing number $c$ and number of $\sim$-equivalence classes
($\#\sim$). \protect\label{tab1}}
\vss}
\newpage
\end{table}


Nonetheless, one can obtain some interesting information already
from the data in the table, for example:

\begin{prop}
The number of alternating genus 3 knots of odd and even
crossing number grows in the ratio $\myfrac{42}{37}$. \qed
\end{prop}

This is certainly not a fact one would expect from considering the
genus two case.

\subsection{The achiral alternating knots}

As the condition a knot to be achiral is relatively
restrictive, I tried, similarly as for genus 2, to
consider the achiral alternating knots of genus three,
hoping to reduce significantly the number of cases and
to obtain an interesting collection of knots.
As we saw, a knot to generate a series with an
achiral alternating knot, it must be in particular of
even crossing number, zero signature and even number of
$\sim$-equivalence classes of crossings. (In fact, among these
classes there must be equally many of both signs for the same number
1 or 2 of elements.) From the generators compiled
above, 68 passed these tests. 

To deal with these 68 cases more conveniently, it is
worth mentioning a further simple criterion which can be often
useful. It uses Gau\ss{} sums (see for the definitions 
\cite{pos,ninv,VirPol}).

\begin{prop}\label{pachir}
Let $K$ be the alternating generator of a series
containing an alternating achiral knot $K'$. Then the following
Gau\ss{} sums vanish on an(y) alternating diagram of $K$:
\[
\begin{array}{c}
\CD{\arrow{-90}{90}}w_p\,\mbox{(writhe)}\,,\quad
\CD{\arrow{-120}{60}\arrow{-60}{120}}\frac{w_p+w_q}{2}\,,\quad
\CD{\chrd{-120}{60}\chrd{-60}{120}\chrd{0}{180}}w_pw_qw_r\,,\\[6mm]
\CD{\chrd{-120}{60}\chrd{-60}{120}\chrd{0}{180}}w_p+w_q+w_r\,,\quad
\CD{\chrd{-70}{70}\chrd{-110}{110}\chrd{-180}{0}}w_pw_qw_r\,,\quad
\CD{\chrd{-70}{70}\chrd{-110}{110}\chrd{-180}{0}}w_p+w_q+w_r\,,\quad
\CD{\labch{-70}{70}{p}\labch{-110}{110}{q}\chrd{-180}{0}}\frac{w_p+w_q}{2}\,.
\end{array}
\]
\end{prop}

\proof The the intersection graph of the Gau\ss{} diagram (IGGD)
of $K'$ has an automorphism taking each vertex to one with the
opposite sign. But building $K$ out of $K'$ means reducing the
number of elements in a $\sim$-equivalence class in the IGGD
to 1 or 2 according to their parity, and hence the above automorphism
carries over to (the IGGD of) $K$. But the above Gau\ss{} sums
are clearly invariants of the intersection graph (and not only
of the Gau\ss{} diagram). They change sign under mirroring the
diagram, and hence the result follows. \qed


The proof suggests that more is likely.

\begin{conj}
If $K$ is the alternating generator of a series
containing an alternating achiral knot, then 
{\def\theenumi{(\@roman\c@enumi)}
\def\labelenumi{\theenumi}
\begin{enumerate}
\item\label{itema} $K$ is achiral, or
\item\label{itemb} $K$ is an iterated mutant of its obverse, or
\item\label{itemc} $K$ has self-conjugate HOMFLY and/or Kauffman polynomial.
\end{enumerate}}
\end{conj}

Clearly \reference{itema} and \reference{itemb} are stronger
than our result. But beware that \reference{itemc} is not.
Remarkably some the above simple Gau\ss{} sums can sometimes
do better in distinguishing an alternating knot from its obverse
than the HOMFLY and/or Kauffman polynomial, as one can see from
the examples $10_{48}$ and $10_{71}$.

It is a good exercise to apply the above criteria by hand in some
simple examples. However, for many and/or more complicated diagrams
it is easier and safer to use computer.

Applying proposition \reference{pachir} on the 68 knots
only the 30 achiral (\em{without} regard of orientation) knots
remained. Up to 14 crossing the list is $8_9$, $8_{17}$, $8_{18}$,
$10_{43}$, $10_{45}$, $10_{81}$, $10_{88}$, $10_{115}$, $12_{125}$,
$12_{273}$, $12_{477}$, $12_{510}$, $12_{960}$, $12_{1124}$,
$12_{1251}$, $14_{1202}$, $14_{5678}$, $14_{15366}$, $14_{16078}$,
$14_{16857}$ and $14_{17247}$. There are 6 knots of 16 crossings,
two of 18 and one of 20 crossings.

Again one can study more detailedly their series as for
$\tl g=2$. For example, we have

\begin{prop}
The fibered achiral alternating genus 3 knots are: $8_9$, $8_{17}$, $8_{18}$,
$10_{43}$, $10_{45}$, $10_{81}$, $10_{88}$, $10_{115}$, $12_{125}$, $12_{477}$
and $12_{1124}$. \qed
\end{prop}

Since the maximal number of $\sim$-equivalence classes of these
30 knots is 12 ($16_{277679}$, $16_{309640}$ and the two 18 crossing
knots have that many), we have

\begin{prop}
The number of prime achiral alternating genus three knots of
$n$ crossings is $O^\asymp(n^{5})$. \qed\vspace{5mm}
\end{prop}

\section{Questions\label{S13}}

\begin{question}
Are there any composite (other than the obvious ones)
or satellite knots of $\tl g=2$?
\end{question}


The lack of ``exotic'' composite $\tl g=2$ knots is suggested
by a conjecture of Cromwell:

\begin{conjecture}(Cromwell \cite{Cromwell2})
If $D$ is a diagram of a composite knot $K=K_1\#K_2$ and
$g(D)=\tl g(K)$, then $D$ is composite.
\end{conjecture}

The conjecture is true by Cromwell's work if $D$ is a
diagram of a closed positive braid and my Menasco's work \cite{Menasco}
if $D$ is alternating. However, the conjecture in general turns
out wrong, as shows the example of figure \ref{figCr}, discovered
in the course of the work previously described here.

\begin{figure}[htb]
\[
\begin{array}{*1c}
\eepsfs[6cm]{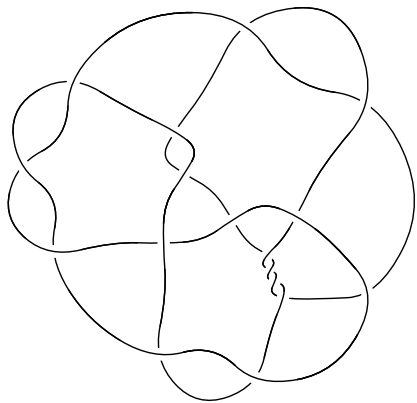}
\end{array}
\]
\caption{\label{figCr}A counterexample to a conjecture of Cromwell:
a prime genus two diagram of the knot $5_2\#!5_2$.}
\end{figure}

We can pose, however, a different problem:

\begin{question}
Does any knot have only finitely many reduced diagrams of
minimal (weak) genus?
\end{question}

It is an easy observation (similar to the proof of proposition
\reference{cor4}) that there are infinitely many slice knots of
$\tl g=2$. (See also remark \reference{rem10.4}.) Take $13_{4233}$.
Then switching 2 of the clasps we obtain a knot bounding a
ribbon disc with two singularities, and can change by twists
the half-twist crossings: 
\[
\eepsfs[3cm]{k-13-4233-m1}\,.
\]

\begin{question}
Can one decide more exactly which weak genus two knots are
slice?
\end{question}

Finally, we point out two general problems:

\begin{question}
Is $\tl g$ always additive under connected sum?
\end{question}

In this case the combinatorial nature of $\tl g$ seems to make
the problem much more involved than for $g$ (for which there is
an easy cut-and-paste argument, see \cite{Adams}). Note again that
the answer would be positive if Cromwell's conjecture had been true.

\begin{question}
Is $\tl g$ invariant under mutation?
\end{question}

\noindent{\bf Acknowledgement.} This paper benefited from
the contributions of several people. Original motivation for
it was L.\ Rudolph's paper 
\cite{Rudolph}, which made me realize how the Seifert algorithm
can be used together with the ``slice Bennequin equality''
to classify $k$-almost positive unknot diagrams. I am also
grateful to Y.\ Nakanishi, Y.\ Rong and V.~Jones for their
references and remarks on $k$-moves, to U.\ Kaiser and D.\ Gabai for
their hint in remark \reference{remG}, and to T.~Nakamura for
sending me a copy of his paper.

For the experimental part, the program KnotScape \cite{KnotScape}
did a huge amount of useful work. (See the survey article \cite{HTW}
for an introduction to this program.) In particular I acknowledge
M.\ Thistlethwaite's collaborative assistance in doing a part of
the polynomial calculations, and also the help of M.\ Heusener in
computing some of the hyperbolic volumes needed in \S\reference{Shyp}.

Finally,
I would wish to thank to the organizers of the ``Knots in Hellas '98''
conference in Delphi, Greece, and among them especially to Jozef
Przytycki and Sofia Lambropoulou for the invitation and giving me the 
possibility to give a talk and for their valuable help and support on 
both mathematical and non-mathematical matters during the conference.

%
%
%
%
%
%

\end{document}

%% file: myeqn.tex
\newenvironment{myeqn*}[1]{\begingroup\def\@eqnnum{\reset@font\rm#1}%
\xdef\@tempk{\arabic{equation}}\begin{equation}\edef\@currentlabel{#1}}
{\end{equation}\endgroup\setcounter{equation}{\@tempk}\ignorespaces}

\newenvironment{myeqn}[1]{\begingroup\let\eq@num\@eqnnum
\def\@eqnnum{\bgroup\let\r@fn\normalcolor 
\def\normalcolor####1(####2){\r@fn####1#1}%
\eq@num\egroup}%
\xdef\@tempk{\arabic{equation}}\begin{equation}\edef\@currentlabel{#1}}
{\end{equation}\endgroup\setcounter{equation}{\@tempk}\ignorespaces}

\newenvironment{myeqn**}{\begin{myeqn}{(\arabic{equation})\es\es\mbox{\qed}}\edef\@currentlabel{\arabic{equation}}}
{\end{myeqn}\stepcounter{equation}}